# General Fractional Vector Calculus

Vasily E. Tarasov [1]

[1] Skobeltsyn Institute of Nuclear Physics,
Lomonosov Moscow State University, Moscow 119991, Russia
E-mail: tarasov@theory.sinp.msu.ru

A generalization of fractional vector calculus (FVC) as a self-consistent mathematical theory is proposed to take into account a general form of non-locality in kernels of fractional vector differential and integral operators. Self-consistency involves proving generalizations of all fundamental theorems of vector calculus for generalized kernels of operators. In the generalization of FVC from power-law nonlocality to the general form of nonlocality in space, we use the general fractional calculus (GFC) in the Luchko approach, which was published in 2021. This paper proposed the following: (I) Self-consistent definitions of general fractional differential vector operators: the regional and line general fractional gradients, the regional and surface general fractional curl operators, the general fractional divergence are proposed. (II) Self-consistent definitions of general fractional integral vector operators: the general fractional circulation, general fractional flux and general fractional volume integral are proposed. (III) The general fractional gradient, Green's, Stokes' and Gauss's theorems as fundamental theorems of general fractional vector calculus are proved for simple and complex regions. The fundamental theorems (Gradient, Green, Stokes, Gauss theorems) of the proposed general FVC are proved for a wider class of domains, surfaces and curves. All these three parts allow us to state that we proposed a calculus, which is a general fractional vector calculus (General FVC). The difficulties and problems of defining general fractional integral and differential vector operators are discussed to the nonlocal case, caused by the violation of standard product rule (Leibniz rule), chain rule (rule of differentiation of function composition), and semigroup property. General FVC for orthogonal curvilinear coordinates, which includes general fractional vector operators for the spherical and cylindrical coordinates, is also proposed.

## 1   Introduction

Vector differential and integral operators are important in various fields of mechanics and physics. In standard vector calculus, integrals and derivatives of integer order are used, and therefore this mathematical tool cannot be used to describe systems, media and fields with nonlocality in space. In this regard, it is important to generalize the vector calculus for applications to the description of non-local media and systems. Nonlocality in space means a dependence of the process, states and variable at the current point of space on the changes on other points of the space.

To describe nonlocality in space, we can use derivatives and integrals of non-integer orders. There operators form a calculus if these operators of non-integer orders satisfy nonlocal analogs of fundamental theorems of standard calculus. These theorems connect the integral and differential operators of non-integer orders. Such calculus is called fractional calculus (FC), and operators are called fractional derivatives (FD) and fractional integrals (FI) [1, 2, 3, 4, 5] and [6, 7]. The FD and FI have different nonstandard properties. For example, the standard product rule, the standard chain rule, and standard semigroup rule are violated for FD of non-integer order [8, 9, 10, 11]. Fractional derivatives and integrals are actively used to describe non-standard properties of systems, media and fields with nonlocality in space and time in varios subjects of mechanics and physics [12, 13, 14], economics [15, 16], and biology [17].

Attempts to generalize some differential operators of vector calculus have been made since the beginning of the 21st century. The history of the development of fractional vector calculus can be conditionally divided into three stages.

1) At the first stage, definitions of various fractional generalizations of differential vector operators (gradient, divergence, rotor, Laplace operator) were constructed. This stage began with the work of Riesz published in 1949 [18, 19, 20] (see also [21, 22]). The definition of fractional generalizations of the gradient was proposed in the works of Adda in 1998 [23], which is actually based on the Sonin-Letnikov fractional derivative [24], and the works by Tarasov in 2005 [25, 26], which are based on the Caputo fractional derivative. The definition of a fractional curl operators was proposed by Engheta in 1998 [27, 28] and Meerschaert, Mortensen, and Wheatcraft in 2006 [29]. The definition of a fractional divergence was proposed by Meerschaert, Mortensen, and Wheatcraft in 2006 [29].

At the first stage, the proposed definitions of fractional vector differential operators were usually not consistent with each other. Fractional generalizations of integral vector operators (fractional circulation, fractional flux and fractional volume integral) have not been proposed. Fractional generalizations of fundamental theorems of vector calculus (such as the Green's, Stokes' and Gauss's theorems) have not been suggested.

2) At the second stage, definitions of fractional generalizations of differential and integral vector operators, which are consistent with each other, were suggested. Fractional generalizations of fundamental theorems of vector calculus were proved. This phase began with the work of Tarasov published in 2008 [30] and book [31], where the power-law spatial non-locality is considered by using the Caputo fractional derivatives and Riemann-Liouville fractional integrals. The fractional generalizations of the Green's, Stokes' and Gauss's theorems are formulated and proved in [30, 31].

After 2008, other articles began to appear, in which special aspects of self-consistent formulations of the fractional vector calculus are discussed. Let us note these aspects: (a) the product rule for FVC is discussed by Bolster, Meerschaert, and Sikorskii in 2012 [32]; (b) applications of fractional gradient to the fractional advection are considered by D'Ovidio, and Garra in 2014 [33]; (c) the discrete fractional vector calculus on lattices is proposed by Tarasov in 2014 [34]; (d) fractional generalizations of the Helmholtz decomposition are proposed by Ortigueira, Rivero and Trujillo in 2015 [35]; (e) the fractional vector operators are considered on convex domain by Agrawal and Xu in 2015 [36]; (f) the fractional vector calculus that is based on the Grunwald-Letnikov derivatives is discussed by Tarasov in 2015 [37], and then by Ortigueira and Machado in 2018 [38]; (g) the fractional Green and Gauss formulas are considered by Cheng and Dai in 2018 [39].

3) The third stages of the development of fractional vector calculus actually began in 2021. At this stage, the fractional vector calculus as a self-consistent mathematical theory is generalized for general form of non-locality and general form of kernels of fractional vector differential and integral operators. Self-consistent mathematical theory involves proving generalizations of all fundamental theorems of vector calculus for generalized kernels of operators. This stage began with the work of D'Elia, Gulian, Olson, and Karniadakis [40] published in 2021 and based on the generalization of the Meerschaert, Mortensen, and Wheatcraft approach to FVC [29]. We can also state that this stage began with our proposed work, which generalizes the approach to formulation of FVC that is proposed in basic paper [30] (see also Chapter 11 in book [14, pp.241-264]) and gave first self-consistent formulation of FVC. In the generalization of FVC from power-law nonlocality to the general form of nonlocality in space, we proposed to use the general fractional calculus (GFC) in the Luchko approach [50, 51, 52]. In our paper, we proposed the following:

(A) Self-consistent definitions of general fractional differential vector operators: the regional and line general fractional gradients, the regional and surface general fractional curl operators, the general fractional divergence are proposed.

(B) Self-consistent definitions of general fractional integral vector operators: the general fractional circulation, general fractional flux and general fractional volume integral are proposed.

(C) The general fractional gradient, Green's, Stokes' and Gauss's theorems as fundamental theorems of general fractional vector calculus are proved for simple and complex regions. The fundamental theorems (Gradient, Green, Stokes, Gauss theorems) of the proposed general FVC are proved for a wider class of domains, surfaces and curves.

All these three parts allow us to state that in this paper, for the first time, a calculus is proposed, which is a general fractional vector calculus (General FVC).

Let us describe approach, which is used in this paper, with some details.

In fractional calculus, nonlocality is described by the kernel of the operators, which are fractional integrals (FI) and fractional derivatives (FD) of non-integer orders. To take into account various types of nonlocality in space, we can use operators with various types of kernels. It is important to have a general fractional calculus that allows us to describe nonlocality in the most general form. We proposed [41, 42] to use general fractional calculus (GFC) to describe systems, media and fields with general form of nonlocality in space.

In this paper, we proposed to use the general fractional calculus (GFC) to formulate a general fractional vector calculus (General FVC). The term "general fractional calculus" (GFC) was proposed by Kochubei in work [43] in 2011 (see also [44] and [45, 46]). In papers [43, 44], the general fractional derivatives (GFDs) and general fractional integral (GFIs) are defined, and the fundamental theorems of the GFC are proved. The GFC is based on the concept of kernel pairs, which was proposed by N.Ya. Sonin (1849-1915) in 1884 article [47] (see also [48]). ("Sonin" is more correct name of the Russian sientist [49] instead of "Sonine" that is use in French [47]). The very important form of the GFC was proposed by Luchko in 2021 [50, 51, 52]. In works [50, 51], GFD and GFI of arbitrary order are suggested, and the general fundamental theorems for the GFI and GFDs are proved. Operational calculus for equations with general fractional derivatives is proposed in [52]. The GFC is also developed and applied in physics in works [53, 54, 55, 56, 57, 58, 59, 60, 61, 62, 63, 41, 42].

In this paper, we use general fractional derivatives (GFDs), general fractional integrals (GFIs) and fundamental theorems of GFC as mathematical tools to formulate General FVC. The proposed General FVC is built on the basis of the results in the GFC obtained in Luchko's work in 2021 [50, 51, 52].

In this paper, the proofs are detailed. This is due to the fact that many standard rules are violated for fractional integrals and derivatives, including such as Leibniz's rule, chain rule, and semigroup rule. These rules are often used in standard vector analysis to prove theorems. For example, the Stokes' theorem is usually proved by using the chain rule, which cannot be applied to fractional derivatives and general fractional derivatives. Moreover, vector differential operators in fractional calculus become nonlocal, which creates additional difficulties for the accurate formulation and proof of theorems. Nonlocality also leads to the possibility of defining different vector differential operators, instead of one operator in the standard vector calculus. For example, we can define Regional GF Gradient, the Surface GF Gradient and Line GF Gradient. The situation is similar for the GF Divergence and GF Curl operators. Moreover, fractional vector analogs of fundamental theorems are not fulfilled for all these general fractional vector operators. In addition, violation of the chain rule leads to the fact that operators defined in different coordinate systems (Cartesian, cylindrical and spherical) are not related to each other by coordinate transformations.

In Section 2, we consider and proposed general fractional integrals and derivatives. In Section 3, we consider and proposed general fractional integral and derivative for $[a, b]$. In Section 4, we consider and proposed line general fractional integral (Line GFI). In Section 5, we give formulation and proof for general fractional gradient theorem. In Section 6, we proof general fractional Green theorem. In Section 7, we proposed double and surface general fractional integral, and flux. In Section 8, we give formulation and proof general fractional Stokes theorem. In Section 9, the general fractional Gauss theorem is formulated and proved. In Section 10, the equalities for general fractional differential vector operators are proposed. Basic interpretations of general fractional differential vector operators are described. In Section 11, we consider and proposed General FVC for orthogonal curvilinear coordinates (OCC), which includes general fractional vector operators for the spherical and cylindrical coordinates.

## 2  General Fractional Integrals and Derivatives

### 2.1  Definitions of GFI and GFD

Let us assume that the functions $M(x)$ belongs to the space $C_{-1,0}(0, \infty)$, and suppose that there exists a function $K(x) \in C_{-1,0}(0, \infty)$, such that the Laplace convolution of these functions is equal to one for all $x \in (0, \infty)$. The function $f(x)$ belongs to the space $C_{-1,0}(0, \infty)$, if this function can be represented in the form $f(x) = x^p g(x)$, where $-1 < p < 0$, and $g(x) \in C[0, \infty)$.

**Definition 2.1** *Let the functions $(M(x), K(x))$ satisfy the following conditions.*
*1) The Sonin condition for the kernels $M(x)$ and $K(x)$ requires that the equation*

$$\int_0^x M(x - x') K(x') \, dx' = 1 \qquad (1)$$

*holds for all $x \in (0, \infty)$.*
*2) The functions $M(x), K(x)$ belong to the space $C_{-1,0}(0, \infty)$,*

$$M(x), K(x) \in C_{-1,0}(0, \infty), \qquad (2)$$

*where*

$$C_{-1,0}(0, \infty) = \{f(x): f(x) = x^p g(x), \ x > 0, \ -1 < p < 0, \ g(x) \in C[0, \infty)\}. \qquad (3)$$

*The set of the pairs $(M(x), K(x))$ that satisfy condition (1) is called the Sonin set and is denoted by $\mathbb{S}$.*

*The set of the pairs $(M(x), K(x))$ that satisfy conditions (1) and (2) is called the Luchko set and is denoted by $\mathbb{L}$.*

The definition of the set $\mathbb{L}$ was proposed in [50] to build a self-consistent formulation of a calculus of general fractional integrals and derivatives.

To define GFI and GFD, we use the Luchko's approach to general fractional calculus, which is proposed in [50, 51].

**Definition 2.2** *Let $M(x) \in \mathbb{L}$ and $f(x) \in C_{-1}(0, \infty) = C_{-1,\infty}(0, \infty)$. The general fractional integral (GFI) with the kernel $M(x) \in C_{-1,0}(0, \infty)$ is the operator on the space $C_{-1}(0, \infty)$:*

$$I^x_{(M)}: C_{-1}(0, \infty) \to C_{-1}(0, \infty), \tag{4}$$

*that is defined by the equation*

$$I^x_{(M)}[x']f(x') = (M * X)(x) = \int_0^x dx'\, M(x - x')\, f(x). \tag{5}$$

If the functions $M(x)$ and $K(x)$ belong to the Luchko set, then we can define general fractional derivatives $D^x_{(K)}$ and $D^{x,*}_{(K)}$ that are associated with GFI $I^x_{(M)}$.

**Definition 2.3** *Let $M(x), K(x) \in \mathbb{L}$ and $f(x) \in C^1_{-1}(0, \infty)$, i.e. $f^{(1)}(x) \in C_{-1}(0, \infty)$. The general fractional derivatives (GFD) with kernel $K(x) \in C_{-1,0}(0, \infty)$, which is associated with GFI (5), is defined as*

$$D^{x,*}_{(K)}[x']f(x') = (K * f^{(1)})(x) = \int_0^x dx'\, K(x - x')\, f^{(1)}(x') \tag{6}$$

*for $x \in (0, \infty)$. The GFD $D^x_{(K)}$ is defined by the equation*

$$D^x_{(K)}[x']f(x') = \frac{d}{dx}(K * f^{(1)})(x) = \frac{d}{dx}\int_0^x dx'\, K(x - x')\, f(x) \tag{7}$$

*for $x \in (0, \infty)$.*

*Operator (6) is called the GFD of the Caputo type, and operator (7) is called the GFD of the Riemann-Liouville type.*

**Remark 2.1** If the kernel pair $(M(x), K(x))$ belongs to the Luchko set $\mathbb{L}$ [51], the kernel $K(x)$ is called associated kernel to $M(x)$. Note that if $K(x)$ is associated kernel to $M(x)$, then $M(x)$ is associated kernel to $K(x)$. Therefore, if $(M(x), K(x))$ belongs to the set $\mathbb{L}$, then the both $I^x_{(M)}[x']$, $D^x_{(K)}[x']$ and $I^x_{(K)}[x']$, $D^x_{(M)}[x']$ can be used as the general fractional integrals (GFI) and general fractional derivatives (GFD).

**Remark 2.2** As it was proved in [50, 51], operators (6) and (7) are connected (see equation 47 in Definition 4 of [50, p.8]) by the equation

$$D^{x,*}_{(K)}[x']f(x') = D^x_{(K)}[x']f(x') - K(x)\, f(0). \tag{8}$$

The GFI and GFD are connected by the fundamental theorems of general fractional calculus (FT of GFC).

**Theorem 2.1** *(Firts fundamental theorem of GFC). If a pair of kernels $(M(x), K(x))$ belongs to the Luchko set $\mathbb{L}$, then the equality*

$$I^x_{(M)}[x']D^{x',*}_{(K)}[x'']f(x'') = f(x) - f(0) \tag{9}$$

holds for $f(t) \in C_{-1,(K)}(0,\infty)$, where
$$C_{-1,(K)}(0,\infty) := \{f: f(x) = I^x_{(K)}[x']g(x'), \ g(x) \in C_{-1}(0,\infty)\}.$$

*Proof.* This theorem is proved as Theorem 3 in [50, p.9], (see also Theorem 1 in [51, p.6]).

**Theorem 2.2** *(Second fundamental theorem of GFC). If a pair of kernels $(M(x), K(x))$ belongs to the Luchko set $\mathbb{L}$, then the equality*
$$I^x_{(M)}[x']D^{x',*}_{(K)}[x'']f(x'') = f(x) - f(0) \tag{10}$$
*holds for $f(x) \in C^1_{-1}(0,\infty)$, where*
$$C^1_{-1}(0,\infty) := \{f: f^{(1)}(x) \in C_{-1}(0,\infty)\}.$$

*Proof.* This theorem is proved as Theorem 4 in [50, p.11], (see also Theorem 2 in [51, p.7]).

### 2.2  Notations for GFI and GFD operators

To define fractional vector operations, we introduce the operators that correspond to the fractional differentiation and general fractional integration.

The general fractional integral operator (GFI operator) is
$$I^x_{(M)}[x'] := \int_0^x dx' \, M(x - x'). \tag{11}$$
The operator (11) acts on functions $f(x) \in C_{-1}(0,\infty)$ by
$$I^x_{(M)}[x']f(x') = \int_0^x dx' \, M(x - x') f(x'). \tag{12}$$
We define the general fractional differential operator (GFD operator) in the form
$$D^{x,*}_{(K)}[x'] := \int_0^x dx' \, K(x - x') \frac{d}{dx'}. \tag{13}$$
The Caputo operator (13) acts on functions on $f(x) \in C^1_{-1}(0,\infty)$ by
$$D^{x,*}_{(K)}[x']f(x') = \int_0^x dx' \, K(x - x') f^{(1)}(x'). \tag{14}$$
We note that operator (13) can be represented in the form
$$D^{x,*}_{(K)}[x'] := I^x_{(K)}[x'] D^1[x'].$$
The proposed notations of GFI and GFD are more convenient, since it allows us to take into account the variables of integration and differentiations, the range of change for these variables.

### 2.3  Examples of Kernel Pairs from Sonin Set and Luchko Set

Let us give examples of the pair of kernels that belongs to the Sonin set $\mathbb{S}$ and the Luchko set $\mathbb{S}$.

Note that if the kernel $M(t)$ is associated to kernel $K(t)$, then the kernel $K(t)$ is associated to $M(t)$. Therefore, if we have the operators
$$I^t_{(M)}[\tau]X(\tau) = (M * X)(t), \quad D^{t,*}_{(K)}[\tau]X(\tau) = (K * X^{(1)})(t), \tag{15}$$
where the kernel pair $(M(t), K(t))$ belong to the Sonin set $\mathbb{S}$, then we can use the operators
$$I^t_{(K)}[\tau]X(\tau) = (K * X)(t), \quad D^{t,*}_{(M)}[\tau]X(\tau) = (M * X^{(1)})(t). \tag{16}$$
A similar situation with the kernel pairs $(M(t), K(t))$ that belong to the Luchko set $\mathbb{L}$.

**Example 2.1** The pair of the operator kernel
$$M(t) = h_\alpha(t) = \frac{t^{\alpha-1}}{\Gamma(\alpha)}, \qquad (17)$$

$$K(t) = h_{1-\alpha}(t) = \frac{t^{-\alpha}}{\Gamma(1-\alpha)}, \qquad (18)$$

where $0 < \alpha < 1$, belongs to the Sonin set $\mathbb{S}$ and the Luchko set $\mathbb{L}$. These kernels define the well-known Riemann-Liouville fractional integral and the Caputo fractional derivative [4].

**Example 2.2** The pair of the kernels [57, p.3628]:
$$M(t) = h_{\alpha,\lambda}(t) = \frac{t^{\alpha-1}}{\Gamma(\alpha)} e^{-\lambda t} \qquad (19)$$

$$K(t) = h_{1-\alpha,\lambda}(t) + \frac{\lambda^\alpha}{\Gamma(1-\alpha)} \gamma(1-\alpha, \lambda t), \qquad (20)$$

belongs to the Sonin set $\mathbb{S}$ and the Luchko set $\mathbb{L}$, where $0 < \alpha < 1$, and $\lambda \geq 0$, $t > 0$, and $\gamma(\beta, t)$ is the incomplete gamma function
$$\gamma(\beta, t) = \int_0^t \tau^{\beta-1} e^{-\tau} d\tau. \qquad (21)$$

**Example 2.3** The pair of kernels
$$M(t) = (\sqrt{t})^{\alpha-1} J_{\alpha-1}(2\sqrt{t}), \qquad (22)$$

$$K(t) = (\sqrt{t})^{-\alpha} I_{-\alpha}(2\sqrt{t}), \qquad (23)$$

belongs to the Sonin set $\mathbb{S}$ and the Luchko set $\mathbb{L}$ (see [47, 48],[57, p.3627]), if $0 < \alpha < 1$, where
$$J_\nu(t) = \sum_{k=0}^\infty \frac{(-1)^k (t/2)^{2k+\nu}}{k!\Gamma(k+\nu+1)}, \quad I_\nu(t) = \sum_{k=0}^\infty \frac{(t/2)^{2k+\nu}}{k!\Gamma(k+\nu+1)} \qquad (24)$$

are the Bessel and the modified Bessel functions, respectively.

For $\alpha = 0.5$, kernels (22) and (23), takes the form
$$M(t) = \frac{\cos(2\sqrt{t})}{\sqrt{\pi t}}, \qquad (25)$$

$$K(t) = \frac{\cosh(2\sqrt{t})}{\sqrt{\pi t}}. \qquad (26)$$

**Example 2.4** The pair of the kernels (see Eq. 7.15 in [57, p.3629]):
$$M(t) = t^{\alpha-1} \Phi(\beta, \alpha; -\lambda t), \qquad (27)$$

$$K(t) = \frac{\sin(\pi\alpha)}{\pi} t^{-\alpha} \Phi(-\beta, 1-\alpha; -\lambda t), \qquad (28)$$

belongs to the Sonin set $\mathbb{S}$ and the Luchko set $\mathbb{L}$, if $0 < \alpha < 1$, and $\Phi(\beta, \alpha; z)$ is the Kummer function
$$\Phi(\beta, \alpha; z) = \sum_{k=0}^\infty \frac{(\beta)_k}{(\alpha)_k} \frac{z^k}{k!}. \qquad (29)$$

**Example 2.5** The pair of the kernels (see Eq. 7.16, 7.18 in [57, p.3629]):

$$M(t) = 1 + \frac{\lambda}{\Gamma(\alpha)\sqrt{t}}, \tag{30}$$

$$K(t) = \frac{1}{\sqrt{\pi t}} - \lambda\, e^{\lambda^2 t}\, \mathrm{erfc}(\lambda\sqrt{(t)}), \tag{31}$$

belongs to the Sonin set $\mathbb{S}$ and the Luchko set $\mathbb{L}$, if $\lambda > 0$, and $\mathrm{erfc}(z)$ is the complementary error function

$$\mathrm{erfc}(z) = 1 - \mathrm{erf}(z) = 1 - \frac{2}{\sqrt{\pi}}\int_0^t e^{-z^2}\, dz. \tag{32}$$

**Example 2.6** The pair of the kernels (see Eq. 7.17, 7.19 in [57, pp.3629-3630]):
$$M(t) = 1 - \frac{\lambda}{\Gamma(\alpha)} t^{\alpha-1}, \tag{33}$$

$$K(t) = \lambda\, t^{-\alpha}\, E_{1-\alpha, 1-\alpha}[\lambda\, t^{1-\alpha}], \tag{34}$$

belongs to the Sonin set $\mathbb{S}$ and the Luchko set $\mathbb{L}$, where $\lambda > 0$.

**Example 2.7** An important subset $\mathbb{H}$ of the Sonine set was proposed by Andrzej Hanyga in work [61]. In this paper, it was proved that any singular (unbounded in a neighborhood of the point zero) locally integrable completely monotone function is a Sonine kernel. If $M(x)$ belongs to the set $\mathcal{H}$ then its associate kernel $K(x)$ also does belong to $\mathcal{H}$. As an example of a kernel pair from $\mathcal{H}$, we can give the kernel pair [61]

$$M(t) = h_{1-\beta+\alpha}(t) + h_{1-\beta}(t), \tag{35}$$

$$K(t) = t^{\beta-1}\, E_{\alpha,\beta}[-t^\alpha], \tag{36}$$

that belongs to the Sonin set $\mathbb{S}$, if $0 < \alpha < \beta < 1$, and $E_{\alpha,\beta}[z]$ is the two-parameters Mittag-Leffler function

$$E_{\alpha,\beta}[z] = \sum_{k=0}^{\infty} \frac{z^k}{\Gamma(\alpha k + \beta)}, \tag{37}$$

where $\alpha > 0$, and $\beta, z \in \mathbb{C}$. The function (35) is a singular locally integrable completely monotone function.

**Example 2.8** Some functions with the power-logarithmic singularities at the origin can belongs to the Sonin set $\mathbb{S}$ [57, pp.3627-3630]. An example of such pair of the kernels (see Eq. 7.22-7.24 in [57, p.3630]) with the power-logarithmic singularities are the following

$$M(t) = \frac{A - \ln(t)}{\Gamma(\alpha)} t^{\alpha-1}, \tag{38}$$

$$K(t) = \mu_{\alpha,h}(t) = \int_0^\infty \frac{t^{z-\alpha}\, e^{h z}}{\Gamma(z+1-\alpha)}\, dz, \tag{39}$$

where $\mu_{\alpha,h}(t)$ is the Volterra function with
$$h = \frac{\Gamma'(\alpha)}{\Gamma(\alpha)} - A.$$

Such kernels do not belong to the set $C_{-1}(0, \infty)$.

**Example 2.9** An important subset $\mathbb{K}$ of the Sonine set was proposed by Anatoly N. Kochubei in works [43, 44]. The kernels, which belong to the subset $\mathcal{K}$, were defined in terms of their Laplace transforms and using the Sonine condition in the Laplace domain.

# 3 General Fractional Integral and Derivative for $[a,b]$

## 3.1 Definition of GFI and GFD for $[a,b]$

Let us define the general fractional integration and differentiation on $[a,b]$ with $0 \leq a < b$.

**Definition 3.1** *Let function $F(x')$ belongs to the funtion space $C_{-1}(0,\infty)$ and $(M(x), K(x)) \in \mathbb{L}$.*
*Then the general fractional integration on $[a,x]$ with $0 \leq a < x$ is defined by the equation*

$$I^{(M)}_{[a,x]}[x]\, F(x) := \begin{cases} I^{x}_{(M)}[x']\, F(x') - I^{a}_{(M)}[x']\, F(x') & x > a > 0, \\ I^{x}_{(M)}[x']\, F(x') & x > a = 0. \end{cases}$$

For $b > a > 0$ and functions $F(x) \in C_{-1}(0,\infty)$, we define the general fractional integration on $[a,b]$ by the equation

$$I^{(M)}_{[a,b]}[x]\, F(x) := I^{b}_{(M)}[x]\, F(x) - I^{a}_{(M)}[x]\, F(x) =$$

$$\int_0^b dx\, M(b-x)\, F(x) - \int_0^a dx\, M(a-x)\, F(x).$$

If $a = 0$ and $b > 0$, then

$$I^{(M)}_{[a,b]}[x]\, F(x) := I^{b}_{(M)}[x]\, F(x).$$

The FGD for $[a,b]$ are defined similarly.

**Definition 3.2** *Let function $F(x)$ belongs to the funtion space $C^{1}_{-1}(0,\infty)$ and $(M(x), K(x)) \in \mathbb{L}$.*
*Then the general fractional differenciation on $[a,x]$ with $0 \leq a < x$ is defined by the equation*

$$D^{(K)}_{[a,x]}[x]\, F(x) := \begin{cases} D^{x,*}_{(K)}[x']\, F(x') - D^{a,*}_{(K)}[x']\, F(x') & x > a > 0, \\ D^{x,*}_{(K)}[x']\, F(x') & x > a = 0. \end{cases}$$

For $b > a > 0$ and functions $F(x) \in C^{1}_{-1}(0,\infty)$, we define the general fractional derivative on $[a,b]$ by the equation

$$D^{(K)}_{[a,b]}[x]\, F(x) := D^{b,*}_{(K)}[x]\, F(x) - D^{a,*}_{(K)}[x]\, F(x) =$$

$$\int_0^b dx\, K(b-x)\, F^{(1)}(x) - \int_0^a dx\, M(a-x)\, F^{(1)}(x).$$

If $a = 0$ and $b > 0$, then

$$D^{(K)}_{[a,b]}[x]\, F(x) := D^{b,*}_{(K)}[x]\, F(x),$$

where $F^{(1)}(x) = dF(x)/dx$.

The GFI operators are the map

$$I^{(M)}_{[a,b]}[x]:\quad C_{-1}(0,\infty) \to C_{-1}(0,\infty).$$

**Remark 3.1** For arbitrary $a, b \geq 0$, we can define the GFI operator as
$$I^{(M)}_{[a,b]}[x] F(x) := sgn(b - a) I^{(M)}_{\omega[a,b]}[x] F(x)$$

$$D^{(K)}_{[a,b]}[x] F(x) := sgn(b - a) D^{(K)}_{\omega[a,b]}[x] F(x)$$

with the signum function
$$sgn(x) := \begin{cases} -1 & (x < 0), \\ 0 & (x = 0), \\ 1 & (x > 0), \end{cases}$$

and
$$\omega[a,b] := \begin{cases} [a,b] & (b > a), \\ \{a\} & (b = a), \\ [b,a] & (b < a). \end{cases}$$

We will use the notation
$$\mathbb{R}^n_+ = \{(x_1, ..., x_n): x_1 > 0, ..., x_n > 0\},$$

$$\mathbb{R}^n_{0,+} = \{(x_1, ..., x_n): x_1 \geq 0, ..., x_n \geq 0\},$$

where $\mathbb{R}^1_+ = (0, \infty)$.

## 3.2 Properties of GFI and GFD on $[a, b]$

Let us prove the property of the additivity for the proposed GFI and GFD on $[a, b]$.

**Theorem 3.1** Let $f(x)$ belongs to the space $C_{-1}(0, \infty)$ and $c > b > a \geq 0$. Then, the property is satisfied
$$I^{(M)}_{[a,b]}[x] f(x) + I^{(M)}_{[b,c]}[x] f(x) = I^{(M)}_{[a,c]}[x] f(x). \tag{40}$$

Let $f(x)$ belongs to the space $C^1_{-1}(0, \infty)$ and $c > b > a \geq 0$. Then. the equations is satisfied
$$D^{(K)}_{[a,b]}[x] f(x) + D^{(K)}_{[b,c]}[x] f(x) = D^{(K)}_{[a,c]}[x] f(x). \tag{41}$$

*Proof.* If $f(x) \in C_{-1}(0, \infty)$, then there are $I^{(M)}_{[a,b]}[x] f(x)$ and $I^{(M)}_{[b,c]}[x] f(x)$.
For $c > b > a > 0$, using the definition, we get
$$I^{(M)}_{[a,b]}[x] f(x) + I^{(M)}_{[b,c]}[x] f(x) =$$

$$(I^b_{(M)}[x] f(x) - I^a_{(M)}[x] f(x)) + (I^c_{(M)}[x] f(x) - I^b_{(M)}[x] f(x)) =$$

$$I^c_{(M)}[x] f(x) - I^a_{(M)}[x] f(x) = I^{(M)}_{[a,c]}[x] f(x).$$

For $c > b > a = 0$, using the definition, we have

$$I^{(M)}_{[0,b]}[x]\,f(x) + I^{(M)}_{[b,c]}[x]\,f(x) =$$

$$I^{b}_{(M)}[x]\,f(x) + (I^{c}_{(M)}[x]\,f(x) - I^{b}_{(M)}[x]\,f(x)) =$$

$$I^{c}_{(M)}[x]\,f(x) = I^{(M)}_{[0,c]}[x]\,f(x).$$

If $f(x) \in C^1_{-1}(0,\infty)$, then there are $D^{(K)}_{[a,b]}[x]\,f(x)$ and $D^{(K)}_{[b,c]}[x]\,f(x)$. The proof of the additivity of the operators is proved in a similar way.

### 3.3  Fundamental Theorems for GFI and GFD on $[a,b]$

Let us formulate the second FT of GFC for the operator $I^{(M)}_{[a,x]}[s]$.

**Theorem 3.2** *Let $f(x)$ belongs to the space $C^1_{-1}(0,\infty)$, and the pair of kernels $(M(t), K(t))$ the Luchko set $\mathbb{L}$. Then*

$$I^{(M)}_{[a,x]}[s]\,D^{S,*}_{(K)}[x']\,f(x') = f(x) - f(a), \tag{42}$$

*where $x > a \geq 0$.*

*Proof.* Let us assume that $f(s) \in C^1_{-1}(0,\infty)$ and
$$F(s) = D^{S,*}_{(K)}[x']f(x') \in C_{-1}(0,\infty).$$
Then using the second fundamental theorem GFC in the form
$$I^{x}_{(M)}[s]\,D^{S,*}_{(K)}[x']\,f(x') = f(x) - f(0),$$
which holds for all $x > 0$, we obtain
$$I^{(M)}_{[a,x]}[s]\,D^{S,*}_{(K)}[s']\,f(s') := I^{x}_{(M)}[s]\,D^{S,*}_{(K)}[x']\,f(x') - I^{a}_{(M)}[s]\,D^{S,*}_{(K)}[x']\,f(x') = .$$

$$(f(x) - f(0)) - (f(a) - f(0)) = f(x) - f(a).$$

**Corollary 3.1** *For $x = b$, equation (42) takes the form*
$$I^{(M)}_{[a,b]}[x]\,D^{x,*}_{(K)}[x']\,f(x') = f(b) - f(a).$$

Let us formulate the first FT of GFC for the operator $I^{(M)}_{[a,x]}[s]$.

**Theorem 3.3** *Let us assume that the conditions of the first FT of GFC for the operators $I^{x}_{(M)}$ and GFD $D^{x,*}_{(K)}$ are satisfied. Then*
$$D^{x,*}_{(K)}[s]\,I^{(M)}_{[a,s]}[x']\,F(x') = F(x),$$
*where we use*
$$D^{x,*}_{(K)}[s]\,1 = 0.$$

*Proof.* Then
$$D^{x,*}_{(K)}[s]\, I^{(M)}_{[a,s]}[x']\, F(x') := D^{x,*}_{(K)}[s]\, I^{s}_{(M)}[x']\, F(x') - D^{x,*}_{(K)}[s]\, I^{a}_{(M)}[x']\, F(x') =$$

$$F(x) - \left(I^{a}_{(M)}[x']\, F(x')\right)\left(D^{x,*}_{(K)}[s]\, 1\right) = F(x),$$

where we use
$$D^{x,*}_{(K)}[s]\, 1 = 0.$$

For functions $F(x) \in C^1_{-1}(0, \infty)$, we define the general fractional derivative on $[a,b]$ by the equation
$$D^{(K)}_{[a,x]}[x']\, F(x') := D^{x,*}_{(K)}[x']\, F(x') - D^{a,*}_{(K)}[x']\, F(x').$$

For the operator $D^{(K)}_{[a,x]}$ and GFI $I^{x}_{(M)}$, the first fundamental theorem GFC is satisfied in the following form.

**Theorem 3.4** *Let us assume that the conditions of the second FT of GFC for the operators $I^{x}_{(M)}$ and GFD $D^{x,*}_{(K)}$ are satisfied. Let us assume that $f(s) \in C_{-1}(0, \infty)$ and*
$$F(s) = I^{s}_{(M)}[x']f(x') \in C^1_{-1}(0, \infty).$$

*Then*
$$D^{(K)}_{[a,x]}[s]\, I^{s}_{(M)}[x']\, f(x') = f(x) - f(a).$$

*Proof.* Using the first fundamental theorem GFC in the form
$$D^{x,*}_{(K)}[x']\, I^{x',*}_{(M)}[s]\, f(s) = f(x),$$
which holds for all $x > 0$, we obtain
$$D^{(K)}_{[a,x]}[s]\, I^{s}_{(M)}[x']\, f(x') := D^{x,*}_{(K)}[s]\, I^{s}_{(M)}[x']\, f(x') - D^{a,*}_{(K)}[s]\, I^{s}_{(M)}[x']\, f(x') =$$

$$f(x) - f(a),$$

where $x > a \geq 0$.

**Remark 3.2** The second fundamental theorem GFC does not hold for the operator $D^{(K)}_{[a,x]}$ with $a > 0$, since $I^{x}_{(M)}[x]\, 1 \neq 0$.

# 4 Line General Fractional Integral (Line GFI)

## 4.1 Simple Line in $\mathbb{R}^2_{0,+}$

Let us define concept of simple line in $\mathbb{R}^2_{0,+}$ of the $XY$-plane.

**Definition 4.1** Let a line $L \subset \mathbb{R}^2_{0,+}$ be described by the function
$$y = y(x) \geq 0, \quad x \in [a,b] \subset \mathbb{R}_{0,+}, \tag{43}$$
which is a continuously differentiable function for all $x \in (a,b) \subset \mathbb{R}_{0,+}$, i.e. $y(x) \in C^1(a,b)$. Then the line $L$ will be called $Y$-simple line on the $XY$ plane.

**Remark 4.1** It is possible to weaken the condition of continuous differentiability of the function $y(x)$ at all points of the interval $(a,b)$, i.e. $y^{(1)}(x) \in C(a,b)$. We can consider the $Y$-simple lines, for which
$$y^{(1)}(x_k - 0) \neq y^{(1)}(x_k + 0)$$
for finite number of points $x_k \in [a,b]$ with $y(x_k - 0) = y(x_k + 0)$. For example, these conditions can be used for broken lines (polygonal chanins).

If $L$ is $Y$-simple line, then every line, which is parallel to the $Y$-axis, intersects the line $L$ at most once for $x \in [a,b]$.

Similarly, line $L \subset \mathbb{R}^2_{0,+}$ is called $X$-simple line on the $XY$ plane, if $L$ can be described by the function
$$x = x(y) \geq 0, \tag{44}$$
which is a continuously differentiable function for $y \in (c,d) \subset \mathbb{R}_{0,+}$, $x(y) \in C^1(c,d)$.

**Definition 4.2** The line $L$ in $\mathbb{R}^2_{0,+}$ of the $XY$ plane is called simple line in $\mathbb{R}^2_{0,+}$ of the $XY$ plane, if $L$ is $X$-simple and $Y$-simple line.

Let $L \subset \mathbb{R}^2_{0,+}$ is $X$-simple line on the $XY$ plane that is described by the single valued function $y = y(x)$ for $x \in [a,b]$. If the functions has inverse function, then $L$ is simple line on the $XY$ plane. It is known that if the function $y = f(x)$ has the derivative $f^{(1)}(x) = df(x)/dx$ for $x \in [a,b]$ such that $f^{(1)}(x) > 0$ (or $f^{(1)}(x) < 0$) for all $x \in [a,b]$, then there exists an inverse function $x = f^{-1}(y)$ for $y \in [c,d]$, where $c = y(a)$ and $d = y(b)$.

Therefore we can formulate the following theorem.

**Theorem 4.1** Let $L \subset \mathbb{R}^2_{0,+}$ is $Y$-simple line on the $XY$ plane, where $y = y(x) \in C^1(a,b)$ and $y^{(1)}(x) > 0$ (or $y^{(1)}(x) < 0$) for all $x \in [a,b]$.

Then there is an inverse function $x = x(y) \in C^1(c,d)$, the line $L$ is $X$-simple line with $y \in [c,d]$, and $L$ is simple line.

**Remark 4.2** In general FVC, we can consider smooth lines on $\mathbb{R}_{0,+}$ that consist of lines, which are simple lines with respect to one of the axes and lines parallel to one of the axes (X, Y, Z).

## 4.2 Simple Line in $\mathbb{R}^3_{0,+}$

Let us define concepts of simple line in $\mathbb{R}^3_{0,+}$.

**Definition 4.3** *Let $L \subset \mathbb{R}^3_{0,+}$ be a line that is described by the functions*
$$y = y(x) \geq 0, \quad z = z(x) \geq 0, \quad x \in [a,b] \subset \mathbb{R}_{0,+}, \tag{45}$$
*which are continuously differentiable functions for $x \in (a,b) \subset \mathbb{R}_{0,+}$, i.e. $y(x), z(x) \in C^1(a,b)$. Then this line will be called YZ-simple line.*

If $L$ is $YZ$-simple line, then every plane, which is parallel to the $YZ$-plane, intersects the line $L$ at most once for $[a,b]$.

**Definition 4.4** *The line $L$ is called simple line in $\mathbb{R}^3_{0,+}$, if $L$ is $XY$-, $XZ$- and $YZ$-simple line.*

If the derivatives $y_x^{(1)}(x)$ and $z_x^{(1)}(x)$ are nonzero and do not change sign on the interval $x \in [a,b]$, then there exist the inverse functions $x = x(y) \in C^1(c,d)$ and $x = x(z) \in C^1(e,f)$. The derivatives of the functions $x(y)$ and $x(z)$ are defined by the equations
$$x_y^{(1)} = \frac{1}{y_x^{(1)}}: \quad x_y^{(1)}(y) = \left(\frac{1}{y_x^{(1)}(x)}\right)_{x=x(y)}.$$
It is obvious that these derivatives of the functions are nonzero and do not change sign also.

We can also state that there exist the functions
$$y = y(z) := y(x(z)) \in C^1(e,f), \quad z = z(y) := z(x(y)) \in C^1(c,d),$$
and the derivatives of these functions that are are nonzero and do not change sign.

As a result, we proved the following theorem.

**Theorem 4.2** *The $YZ$-simple line $L \subset \mathbb{R}^3_{0,+}$, for which the derivatives of the functions $y = y(x)$ and $z = z(x)$ are nonzero and do not change sign on the interval $x \in [a,b]$, is $XY$-simple, $XZ$-simple line, and therefore simple line in $\mathbb{R}^3_{0,+}$.*

Let $L$ be a simple line in $\mathbb{R}^3_{0,+}$ and the lines $L_{xy}$, $L_{xz}$ and $L_{yz}$ are projections of $L$ on the $XY$-, $XZ$-, $YZ$ planes. Then $L_{xy}$, $L_{xz}$ and $L_{yz}$ are simple lines in $\mathbb{R}^2_{0,+}$ of the $XY$-, $XZ$-, $YZ$ planes.

**Remark 4.3** The simple line $L \subset \mathbb{R}^3_{0,+}$, which connects the points $A(a,c,e)$ and $B(b,d,f)$, can be defined by three equivalent forms
$$L = \{(x,y,z): \ x \in [a,b], \ y = y(x) \in C^1[a,b], \ z = z(x) \in C^1[a,b]\}, \tag{46}$$

$$L = \{(x,y,z): \ y \in [c,d], \ x = x(y) \in C^1[c,d], \ z = z(y) \in C^1[c,d]\}, \tag{47}$$

$$L = \{(x,y,z): \ z \in [e,f], \ x = x(z) \in C^1[e,f], \ y = y(z) \in C^1[e,f]\}, \tag{48}$$
where $y(x=a) = c$, $y(x=b) = d$, $z(x=a) = e$, $z(x=b) = f$.

**Remark 4.4** We can also consider the simple lines, for which the function that are not differentiable at a finite number of points. These conditions allows us to use for broken lines (polygonal chanins).

## 4.3 Problems with Definition of Line GFI of Vector Field

Let us consider the $Y$-simple line $L_1 \subset \mathbb{R}^2_{0,+}$ in the $XY$-plane that is described by the equation
$$y = y(x) \in C^1[a,b],$$
and the vector field
$$F(x,y) = e_x F_x(x,y) + e_y F_y(x,y),$$
where $(x,y) \in \mathbb{R}^2_{0,+}$. Then the standard line (curve) integral of second kind can be defined by the equation
$$\int_{L_1} (F_x(x,y)\, dx + F_y(x,y)\, dy) := \int_a^b (F_x(x,y(x))\, dx + F_y(x,y(x))\, y^{(1)}(x)\, dx). \quad (49)$$

If the $X$-simple line $L_2 \subset \mathbb{R}^2_{0,+}$ in the $XY$-plane is described by the equation
$$x = x(y) \in C^1[c,d],$$
then the line integral can be defined by the equation
$$\int_{L_2} (F_x(x,y)\, dx + F_y(x,y)\, dy) := \int_c^d (F_x(x(y),y)\, x^{(1)}(y)\, dy + F_y(x(y),y)\, dy). \quad (50)$$

If $y = y(x) \in C^1(a,b)$ and $y^{(1)}(x) > 0$ (or $y^{(1)}(x) < 0$) for all $x \in [a,b]$. Then there is an inverse function $x = x(y) \in C^1(c,d)$, and the $Y$-simple line $L_1$ is $X$-simple line. Therefore the line integral (49) can be represented as (50), and we have the equality
$$\int_a^b (F_x(x,y(x))\, dx + F_y(x,y(x))\, y^{(1)}(x)\, dx) =$$

$$\int_c^d (F_x(x(y),y)\, x^{(1)}(y)\, dy + F_y(x(y),y)\, dy). \quad (51)$$

This equality is based on the property
$$\int_c^d F_y(x(y),y)\, dy = \int_a^b F_y(x,y(x))\, y^{(1)}(x)\, dx \quad (52)$$

$$\int_a^b F_x(x,y(x))\, dx = \int_c^d F_x(x(y),y)\, x^{(1)}(y)\, dy. \quad (53)$$

Equations (49) and (50) cannot be used to define the line general fractional integral (line GFI) since property (52) and (53) does not hold for general fractional integrals. In general fractional calculus, we have the inequality
$$I^{(M_2)}_{[c,d]}[y]\, F_y(x(y),y) \neq I^{(M_1)}_{[a,b]}[x]\, (F_y(x,y(x))\, y^{(1)}(x))$$
that has the form
$$\int_0^d M_2(d-y)\, F_y(x(y),y)\, dy - \int_0^c M_2(c-y)\, F_y(x(y),y)\, dy \neq$$

$$\int_0^b M_1(b-x)\, F_y(x,y(x))\, y^{(1)}(x)\, dx - \int_0^a M_1(a-x)\, F_y(x,y(x))\, y^{(1)}(x)\, dx.$$

To solve this problem of definition of line GFI, we can use the fact that the line integral of second kind over simple line can be defined by the equation

$$\int_{L_1} (F_x(x,y)\,dx + F_y(x,y)\,dy) := \int_a^b F_x(x, y(x))\,dx + \int_c^d F_y(x(y), y)\,dy. \qquad (54)$$

Therefore the line GFI in $\mathbb{R}^2_{0,+}$ can be defined by the following definition.

**Definition 4.5** Let $L$ be a simple line in $\mathbb{R}^2_{0,+}$ of the $XY$-plane. Let the functions $f_x(x) := F_x(x, y(x))$ and $g_y(y) := F_y(x(y), y)$ belong to the function space $C_{-1}(0, \infty)$. Then line GFI for the line $L$ is defined by the equation
$$(I_L^{(M)}, F) = I_{[a,b]}^{(M_1)}[x]\,F_x(x, y(x)) + I_{[c,d]}^{(M_2)}[y]\,F_y(x(y), y) =$$

$$\int_0^b M_1(b - x)\,F_x(x, y(x))\,dx - \int_0^a M_1(a - x)\,F_x(x, y(x))\,dx +$$

$$\int_0^d M_2(d - y)\,F_y(x(y), y)\,dy - \int_0^c M_2(c - y)\,F_y(x(y), y)\,dy.$$

This line GFI exists, if the kernels $(M_1(x), K_1(x))$ and $(M_2(y), K_2(y))$ belong to the Luchko set $\mathbb{L}$.

The proposed approach to define line GFI for lines in $\mathbb{R}^2_{0,+}$ can be used to define line GFI for lines in $\mathbb{R}^3_{0,+}$. Let $L$ be a simple line in $\mathbb{R}^3_{0,+}$, which is defined in form (46), (47), (48), and the vector field
$$F(x, y, z) = e_x F_x(x, y, z) + e_y F_y(x, y, z) + e_z F_z(x, y, z),$$
where $(x, y, z) \in \mathbb{R}^3_{0,+}$. The standard line integral of second kind over a simple line $L$ in $\mathbb{R}^3_{0,+}$ can be defined by the equation
$$\int_L (F_x(x, y, z)\,dx + F_y(x, y, z)\,dy + F_z(x, y, z)\,dz) :=$$

$$\int_a^b F_x(x, y(x), z(x))\,dx + \int_c^d F_y(x(y), y, z(y))\,dy + \int_e^f F_z(x(z), y(z), z)\,dz. \qquad (55)$$

The line GFI over a simple line $L$ in $\mathbb{R}^3_{0,+}$ can be defined by the equation
$$(I_L^{(M)}, F) = I_{[a,b]}^{(M_1)}[x]\,F_x(x, y(x), z(x)) +$$

$$I_{[c,d]}^{(M_2)}[y]\,F_y(x(y), y, z(y)) + I_{[e,f]}^{(M_3)}[x]\,F_z(x(z), y(z), z). \qquad (56)$$

This line GFI exists, if the kernels $(M_1(x), K_1(x))$, $(M_2(y), K_2(y))$ and $(M_3(z), K_3(z))$ belong to the Luchko set $\mathbb{L}$, and the functions $F_1(x) := F_x(x, y(x), z(x))$ $F_2(y) := F_y(x(y), y, z(y))$, and $F_3(z) := F_z(x(z), y(z), z)$ belong to the function space $C_{-1}(0, \infty)$.

**Remark 4.5** The line integral can also be defined for wide class of lines in $\mathbb{R}^3_{0,+}$ that are not simple lines, if the lines $L$ can be split into finite number of simple lines in $\mathbb{R}^3_{0,+}$. These lines will be called piecewese simple lines.

### 4.4 Definition of Line GFI for Vector Field in $\mathbb{R}^3_{0,+}$

Let us define some conditions on the vector field. We will assume that the vector field
$$F = F(x, y, z) = e_x F_x(x, y, z) + e_y F_y(x, y, z) + e_z F_z(x, y, z)$$
on the simple line $L \subset \mathbb{R}^3_{0,+}$ is described by functions that belong to the space $C_{-1}(0, \infty)$:

$$F_1(x) := F_x(x, y(x), z(x)) \in C_{-1}(0, \infty),$$

$$F_2(y) := F_y(x(y), y, z(y)) \in C_{-1}(0, \infty),$$

$$F_3(z) := F_z(x(z), y(z), z) \in C_{-1}(0, \infty).$$

If these conditions are satisfied, then we will write $F(x, y, z) \in \mathbb{F}_{-1}(L)$.

In the case $F \in \mathbb{F}_{-1}(\mathbb{R}_+^3)$, the line general fractional integral for the vector field $F$ and the line $L = AB$ with endpoints $A(a_1, a_2, a_3)$ and $B(b_1, b_1, b_3)$ with all $b_k > a_k \geq 0$ is defined by the equation

$$(I_L^{(M)}, F) = \sum_{k=1}^{3} I_{[a_k, b_k]}^{(M_k)}[x_k] F_k(x_k) =$$

$$I_{[a_1, b_1]}^{(M_1)}[x] F_1(x) + I_{[a_2, b_2]}^{(M_2)}[y] F_2(y) + I_{[a_3, b_3]}^{(M_3)}[z] F_3(z) =$$

$$(I_{(M_1)}^{b_1}[x] - I_{(M_1)}^{a_1}[x]) F_x(x, y(x), z(x)) +$$

$$(I_{(M_2)}^{b_2}[x] - I_{(M_2)}^{a_2}[y]) F_y(x(y), y, z(y)) +$$

$$(I_{(M_3)}^{b_3}[z] - I_{(M_3)}^{a_3}[z]) F_z(x(z), y(z), z) =$$

$$\int_0^{b_1} dx \, M_1(b_1 - x) F_x(x, y(x), z(x)) - \int_0^{a_1} dx \, M_1(a_1 - x) F_x(x, y(x), z(x)) +$$

$$\int_0^{b_2} dy \, M_2(b_2 - y) F_y(x(y), y, z(y)) - \int_0^{a_2} dy \, M_2(a_2 - y) F_y(x(y), y, z(y)) +$$

$$\int_0^{b_3} dz \, M_3(b_3 - z) F_z(x(z), y(z), z) - \int_0^{a_3} dz \, M_3(a_3 - z) F_z(x(z), y(z), z).$$

For the proposed definition of a line general fractional integral, the kernels of the integral operator remain dependent on the difference of variables

$$M_k = M_k(b_k - x_k), \quad M_k = M_k(a_k - x_k),$$

and the line GFI itself is expressed through the Laplace convolution as a product in the ring $\mathbb{R}_{-1}$. Using the Laplace convolution, the line GFI can be written as

$$(I_L^{(M)}, F) = \sum_{k=1}^{3} ((M_k * F_k)(b_k) - (M_k * F_k)(a_k)),$$

where $b_k > a_k \geq 0$, and

$$(M_k * F_k)(c_k) = \int_0^{c_k} dx_k \, M(c_k - x_k) F_k(x_k),$$

where $c_k = a_k$ or $c_k = b_k$ with $M_k(x_k) \in C_{-1,0}(0, \infty)$ and $F_k(x_k) \in C_{-1}(0, \infty)$. Using the property of GFI, we have $(M_k * F_k)(c_k) \in C_{-1}(0, \infty)$.

We can consider the variables $x_k > a_k$ instead of the numbers $b_k$. Then

$$(I_L^{(M)}, F) = \sum_{k=1}^{3} ((M_k * F_k)(x_k) - (M_k * F_k)(a_k)),$$

where $x_k > a_k \geq 0$.

**Remark 4.6** In general FVC, we can consider lines that consist of lines, which are simple lines and lines parallel to one of the axes (X, Y, Z). As an example of this type of lines, we can consider polygonal chains (broken lines).

## 4.5 Line General Fractional Integral for Polygonal Chains

Let us define a line GFI for the polygonal chains. A polygonal chain (broken line) is a geometric figure consisting of line segments connected in series by their endpoints. A polygonal chain $L$ is a sequence of points $A_1, A_2, \ldots, A_n$ that forms the successively connected line segments $[A_1, A_2], [A_2, A_3], \ldots, [A_{n-1}, A_n]$. The points $A_1, A_2, \ldots, A_n$ are called the vertices of polygonal chain.

Let us consider the polygonal chain
$$L := \bigcup_{k=1}^{n-1} [A_k, A_{k+1}], \tag{57}$$
where the vertices $A_k$ have coordinates $(x_k, y_k, z_k)$, such that $x_k, y_k, z_k \geq 0$ for all $k = 1, 2, \ldots, n-1$.

For polygonal chain (57), the line GFI is defined by the equations
$$I_L^{(M)}[x, y, z] = \sum_{k=1}^{n-1} I_{[A_k, A_{k+1}]}^{(M)}[x, y, z], \tag{58}$$
where
$$I_{[A_k, A_{k+1}]}^{(M)}[x, y, z] := sgn(x_{k+1} - x_k) \, e_x \, I_{\omega[x_k, x_{k+1}]}^{(M_1)}[x] +$$

$$sgn(y_{k+1} - y_k) \, e_y \, I_{\omega[y_k, y_{k+1}]}^{(M_2)}[y] + sgn(z_{k+1} - z_k) \, e_z \, I_{\omega[z_k, z_{k+1}]}^{(M_3)}[z],$$

with the signum function
$$sgn(x) := \begin{cases} -1 & (x < 0), \\ 0 & (x = 0), \\ 1 & (x > 0), \end{cases} \tag{59}$$

and
$$\omega[x_k, x_{k+1}] := \begin{cases} [x_k, x_{k+1}] & (x_{k+1} > x_k), \\ \{x_k\} & (x_{k+1} = x_k), \\ [x_{k+1}, x_k] & (x_{k+1} < x_k). \end{cases} \tag{60}$$

If $x_k < x_{k+1}$, $y_k < y_{k+1}$, $z_k < z_{k+1}$, then
$$I_{[A_k, A_{k+1}]}^{(M)}[x, y, z] := e_x \, I_{[x_k, x_{k+1}]}^{(M_1)}[x] + e_y \, I_{[y_k, y_{k+1}]}^{(M_2)}[y] + e_z \, I_{[z_k, z_{k+1}]}^{(M_3)}[z],$$

If the opposite inequality holds, then a minus sign is put in front of the integral, and in the interval the numbers are set in ascending order. If equality ($x_k = x_{k+1}$ or $y_k = y_{k+1}$ or $z_k = z_{k+1}$) holds, then the integral is considered equal to zero.

Let vector field $F$ be defined by the equation
$$F = F(x, y, z) = e_x F_x(x, y, z) + e_y F_y(x, y, z) + e_z F_z(x, y, z),$$
where
$$F_x(x, y_k(x), z_k(x)), F_y(x_k(y), y, z_k(y)), F_y(x_k(y), y, z_k(y)) \in C_{-1}(0, \infty)$$
for all $k = 1, \ldots, n$.

Then the line GFI for the polygonal chain (broken line) is defined by equation (58), where
$$(I_{[A_k, A_{k+1}]}^{(M)}[x, y, z] \, F(x, y, z)) := sgn(x_{k+1} - x_k) \, I_{\omega[x_k, x_{k+1}]}^{(M_1)}[x] \, F_x(x, y_k(x), z_k(x)) +$$

$$sgn(y_{k+1} - y_k) \, I_{\omega[y_k, y_{k+1}]}^{(M_2)}[y] \, F_y(x_k(y), y, z_k(y)) +$$

$$sgn(z_{k+1} - z_k) \, I_{\omega[z_k, z_{k+1}]}^{(M_3)}[z] \, F_z(x_k(z), y_k(z), z),$$
where the functions $y = y_k(x)$, $x = x_k(y)$, $z = z_k(x)$, $x = x_k(z)$, $z = z(y)$, $y = y_k(x)$ are

defined by the equations

$$\frac{x-x_k}{x_{k+1}-x_k} = \frac{y-y_k}{y_{k+1}-y_k} = \frac{z-z_k}{z_{k+1}-z_k},$$

if the denominators are not zero. For example,

$$y = \frac{y_{k+1}-y_k}{x_{k+1}-x_k}(x-x_k) + y_k,$$

$$z = \frac{z_{k+1}-z_k}{x_{k+1}-x_k}(x-x_k) + z_k.$$

Here

$$I_{[a,b]}^{(M_1)}[x]\, F_x(x, y(x), z(x)) =$$

$$\int_0^b dx\, M_1(b-x)\, F_x(x, y(x), z(x)) - \int_0^a dx\, M_1(b-x)\, F_x(x, y(x), z(x)),$$

where $b > a \geq 0$,

$$I_{[c,d]}^{(M_2)}[x]\, F_y(x(y), y, z(y)) =$$

$$\int_0^d dy\, M_2(d-x)\, F_y(x(y), y, z(y)) - \int_0^c dy\, M_2(c-x)\, F_y(x(y), y, z(y)),$$

where $d > c \geq 0$,

$$I_{[e,f]}^{(M_3)}[x]\, F_y(x(z), y(z), z) =$$

$$\int_0^f dz\, M_3(f-x)\, F_z(x(z), y(z), z) - \int_0^e dz\, M_3(e-x)\, F_z(x(z), y(z), z),$$

where $f > e \geq 0$.

### 4.6 Line GFI for Piecewise Simple Lines

Similarly to the case of a broken line (the polygonal chains), we can define line GFI for line, which consists of simple lines and lines parallel to the axes.

Let us consider a line $L \subset \mathbb{R}_{0,+}^3$, which can be divided into several lines $L_k = L_k[A_k, A_{k+1}]$, $k = 1, \ldots, n$ that are simple lines or lines parallel to one of the axes:

$$L := \bigcup_{k=1}^{n-1} L_k[A_k, A_{k+1}], \tag{61}$$

where the line $L_k$ connects the points $A_k(x_k, y_k, z_k)$, and $A_{k+1}(x_{k+1}, y_{k+1}, z_{k+1})$ with $x_k, y_k, z_k \geq 0$ for all $k = 1, 2, \ldots, n-1$. Lines of this kind will be called the piecewise simple lines.

For piecewise simple line (61) in $\mathbb{R}_{0,+}^3$, and the vector field $F \in \mathbb{F}_{-1}(L)$, the line GFI is defined by the equation

$$(I_L^{(M)}[x,y,z]\, F(x,y,z)) := \sum_{k=1}^{n-1} (I_{[A_k, A_{k+1}]}^{(M)}[x,y,z]\, F(x,y,z)),$$

where

$$(I_{[A_k, A_{k+1}]}^{(M)}[x,y,z]\, F(x,y,z)) := sgn(x_{k+1} - x_k)\, I_{\omega[x_k, x_{k+1}]}^{(M_1)}[x]\, F_x(x, y_k(x), z_k(x)) +$$

$$sgn(y_{k+1} - y_k)\, I_{\omega[y_k, y_{k+1}]}^{(M_2)}[y]\, F_y(x_k(y), y, z_k(y)) +$$

$$sgn(z_{k+1} - z_k)\, I_{\omega[z_k, z_{k+1}]}^{(M_3)}[z]\, F_z(x_k(z), y_k(z), z),$$

where the functions $y = y_k(x)$, $x = x_k(y)$, $z = z_k(x)$, $x = x_k(z)$, $z = z_k(y)$, $y = y_k(x)$ define the lines $L_k$ that are simple or parallel to one of the axes.

## 4.7 General Fractional Circulation for Rectangle

Let us define concepts of general fractional circulation by using the GFC. Note that these concepts for the kernels

$$M_n(x_n) = h_{\alpha_n}(x_n) = \frac{x_n^{\alpha_n-1}}{\Gamma(\alpha_n)}, \quad \text{(for all n = 1,2,3)} \tag{62}$$

were proposed in [30] (see also [31]).

We will consider the piecewise simple line (61), where $A_n(a_n, c_n, e_n) = A_1(a_1, c_1, e_1)$, i.e. $L$ is closed line. Let $L$ be the piecewise simple line

$$L := \bigcup_{k=1}^{n-1} L_k, \tag{63}$$

where the lines $L_k = (A_k, A_{k+1})$ with points $A_k(a_k, c_k, d_k)$ ($a_{k+1} = b_k$, $c_{k+1} = d_k$, $e_{k+1} = f_k$) are simple that are described as

$$L_k := \{(x, y, z): x \in [a_k, b_k], \ y = y_k(x) \in C^1(0, \infty), \ z = z_k(x) \in C^1(0, \infty)\},$$

where the derivatives of the functions $y_k(x)$, $z_k(x)$ are nonzero and do not change sign on the interval $[a_k, b_k]$ foe all $k = 1, \ldots, n$.

The last requirement for derivatives of the functions $y_k(x)$, $z_k(x)$ allows us to represent the line in the form

$$L_k := \{(x, y, z): y \in [c_k, d_k], \ x = x_k(y) \in C^1(0, \infty), \ z = z_k(y) \in C^1(0, \infty)\},$$

$$L_k := \{(x, y, z): z \in [e_k, f_k], \ x = x_k(z) \in C^1(0, \infty), \ y = y_k(z) \in C^1(0, \infty)\},$$

where $y_k(a_k) = c_k$, $y_k(b_k) = d_k$, $z_k(a_k) = e_k$, $z_k(b_k) = f_k$.

Then, we can define the general fractional circulation in the following form

$$\mathbb{E}_L^{(M)}(F) = (I_L^{(M)}[x, y, z] F(x, y, z)) = \tag{64}$$

$$\sum_{k=1}^{n-1} (I_{L_k}^{(M)}[x, y, z] F(x, y, z)) =$$

$$\sum_{k=1}^{n-1} I_{L_k}^{(M_1)}[x] F_x(x, y_k(x), z_k(x)) + \sum_{k=1}^{n-1} I_{L_k}^{(M_2)}[y] F_y(x_k(y), y, z_k(y)) +$$

$$\sum_{k=1}^{n-1} I_{L_k}^{(M_3)}[z] F_z(x_k(z), y_k(z), z) = \sum_{k=1}^{n-1} I_{[a_k, b_k]}^{(M_1)}[x] F_x(x, y_k(x), z_k(x)) +$$

$$\sum_{k=1}^{n-1} I_{[c_k, d_k]}^{(M_2)}[y] F_y(x_k(y), y, z_k(y)) + \sum_{k=1}^{n-1} I_{[e_k, f_k]}^{(M_3)}[z] F_z(x_k(z), y_k(z), z).$$

If you include segments parallel to some axes, then for such segments the integrals are equal to zero.

As a result, we can formulate the following definition.

**Definition 4.6** *A general fractional circulation is a general fractional line integral of the vector field $F \in \mathbb{F}_{-1}(L)$ along a piecewise simple closed line $L$ that is defined by*

$$\mathbb{E}_L^{(M)}(F) = (I_L^{(M)}[x, y, z] F(x, y, z)) = \sum_{k=1}^{n-1} I_{[a_k, b_k]}^{(M_1)}[x] F_x(x, y_k(x), z_k(x)) + \tag{65}$$

$$\sum_{k=1}^{n-1} I_{[c_k, d_k]}^{(M_2)}[y] F_y(x_k(y), y, z_k(y)) + \sum_{k=1}^{n-1} I_{[e_k, f_k]}^{(M_3)}[z] F_z(x_k(z), y_k(z), z),$$

where $a_n = a_1$, $b_n = b_1$, $c_n = c_1$, $d_n = d_1$, $e_n = e_1$, $f_n = f_1$.

**Example 4.1** Let us consider the piecewise simple closed line $L \subset \mathbb{R}^2_{0,+}$ in XY-plane,

which consists of the simple lines $AB$, $BC$, $CD$, and $DA$, with the coordinates of the points
$$A(a, y_1(a)), \quad B(x_1(c), c), \quad C(b, y_1(b)), \quad D(x(d), d),$$
where
$$y_1(a) = y_2(a), \quad y_1(b) = y_2(b), \quad x_1(c) = x_2(c), \quad x_1(d) = x_2(d).$$
The lines $ABC$ and $ADC$ are $Y$-simple lines
$$ABC := \{(x, y): \; 0 \leq a \leq x \leq b, \; y = y_1(x) \geq 0\}, \tag{66}$$

$$ADC := \{(x, y): \; 0 \leq a \leq x \leq b, \; y = y_2(x) \geq 0\}, \tag{67}$$
where the functions $y_1(x)$ and $y_2(x)$ belong to the space $C^1(a, b)$, whose derivatives are nonzero and do not change sign on the interval $[a, b]$.

The lines $BAD$ and $BCD$ are $X$-simple lines
$$BAD := \{(x, y): \; 0 \leq c \leq y \leq d, \; x = x_1(y) \geq 0\}, \tag{68}$$

$$BCD := \{(x, y): \; 0 \leq c \leq y \leq c, \; x = x_2(y) \geq 0\}, \tag{69}$$
where the functions $x_1(x)$ and $x_2(x)$ belong to the space $C^1(c, d)$, whose derivatives do not change sign on the interval $[c, d]$.

Let us assume that the functions $F_x(x, y)$ and $F_y(x, y)$ satisfy the condition
$$F_x(x, y_1(x)), F_x(x, y_2(x)) \in C_{-1}(0, \infty), \tag{70}$$

$$F_y(x_1(y), y), F_y(x(y), y) \in C_{-1}(0, \infty), \tag{71}$$

Then general fractional circulation along the piecewise simple closed line
$$L = AB \cup BC \cup CD \cup DA$$
is defined by the equation
$$\mathbb{E}_L^{(M)}(F) = (I_L^{(M)}[x, y, z] \, F(x, y, z)) = I_L^{(M_1)}[x] \, F_x + I_L^{(M_2)}[y] \, F_y = \tag{72}$$

$$I_{ABC}^{(M_1)}[x] \, F_x(x, y_1(x)) - I_{ADC}^{(M_1)}[x] \, F_x(x, y_2(x)) +$$

$$I_{BCD}^{(M_2)}[y] \, F_y(x_2(y), y) - I_{BAD}^{(M_2)}[y] \, F_y(x_1(y), y).$$

**Example 4.2** Let us consider the piecewise simple line (61), where $A_n = A_1$, i.e. $L$ is closed line. For example, we consider the line GFI for the rectangle on $\mathbb{R}_{0,+}^2$ with vertices at the points
$$A(0,0), \quad B(a, 0), \quad C(a, b), \quad D(0, b). \tag{73}$$
The sides $AB$, $BC$, $CD$, $DA$ of the rectangle form the line $L$. For closed line $L = ABCD$, the line GFI operator $I_L^{(M)}$ is written as
$$I_L^{(M)} = e_x I_{AB}^{(M_1)}[x] + e_y I_{BC}^{(M_2)}[y] + e_x I_{CD}^{(M_1)}[x] + e_y I_{DA}^{(M_2)}[y] =$$

$$e_x I_{AB}^{(M_1)}[x] + e_y I_{BC}^{(M_2)}[y] - e_x I_{DC}^{(M_1)}[x] - e_y I_{AD}^{(M_2)}[y], \tag{74}$$
where we used that
$$I_{CD}^{(M_1)}[x] = -I_{DC}^{(M_1)}[x], \quad I_{DA}^{(M_2)}[y] = -I_{AD}^{(M_2)}[y].$$
For the vector field
$$F(x, y) = e_x F_x(x, y) + e_y F_y(x, y),$$

the line GFI has the form

$$(I_L^{(M)}, F) =$$

$$I_{AB}^{(M_1)}[x] \, F_x + I_{BC}^{(M_2)}[y] \, F_y - I_{DC}^{(M_1)}[x] \, F_x - I_{AD}^{(M_2)}[y] \, F_y,$$

$$I_{[0,a]}^{(M_1)}[x] \, F_x(x,0) + I_{[0,b]}^{(M_2)}[y] \, F_y(a,y) -$$

$$I_{[0,a]}^{(M_1)}[x] \, F_x(x,b) - I_{[0,b]}^{(M_2)}[y] \, F_y(0,y) =$$

$$I_{[0,a]}^{(M_1)}[x](F_x(x,0) - F_x(x,b)) + I_{[0,b]}^{(M_2)}[y](F_y(a,y) - F_y(0,y)). \tag{75}$$

**Example 4.3** The general fractional circulation for the line L that is rectangle with sides AB, BC, CD, DA, where the points have coordinates (73), is written as

$$\mathbb{E}_L^{(M)}(F) = (I_L^{(M)}, F) = I_{AB}^{(M_1)}[x]F_x + I_{BC}^{(M_2)}[y]F_y - I_{CD}^{(M_1)}[x]F_z - I_{DA}^{(M_2)}[y]F_z =$$

$$I_{(M_1)}^a[x]F_x(x,0) + I_{(M_2)}^b[y]F_y(a,y) - I_{(M_1)}^a[x]F_x(x,b) - I_{(M_2)}^b[y]F_y(0,y) =$$

$$I_{(M_1)}^a[x](F_x(x,0) - F_x(x,b)) + I_{(M_2)}^b[y](F_y(a,y) - F_y(0,y)). \tag{76}$$

As a resul, we get

$$\mathbb{E}_L^{(M)}(F) = \int_0^a dx \, M_1(a-x) \, (F_x(x,0) - F_x(x,b)) +$$

$$\int_0^b dy \, M_2(b-y) \, (F_y(a,y) - F_y(0,y)). \tag{77}$$

For kernels (62), the general fractional circulation (77) has the form

$$\mathbb{E}_L^{(M)}(F) = \int_0^a dx \, \frac{(a-x)^{\alpha_1 - 1}}{\Gamma(\alpha_1)} \, (F_x(x,0) - F_x(x,b)) +$$

$$\int_0^b dy \, \frac{(b-y)^{\alpha_2 - 1}}{\Gamma(\alpha_2)} \, (F_y(a,y) - F_y(0,y)), \tag{78}$$

which was proposed in [30] (see also [31]), where the Riemann-Liuoville fractional integrals are used. For $\alpha_1 = \alpha_2 = 1$, we get the standard circulation.

# 5 General Fractional Gradient Theorem

## 5.1 General Fractional Gradient

Let us give definitions of a set of scalar fields and a general fractional gradient for $\mathbb{R}^3_{0,+}$.

**Definition 5.1** Let $U(x,y,z)$ be a scalar field that satisfies the conditions
$$D^{x,*}_{(K_1)}[x']U(x',y,z) \in C_{-1}(\mathbb{R}^3_+),$$

$$D^{y,*}_{(K_2)}[y']U(x,y',z) \in C_{-1}(\mathbb{R}^3_+),$$

$$D^{z,*}_{(K_3)}[z']U(x,y,z') \in C_{-1}(\mathbb{R}^3_+).$$

Then the set of such scalar fields will be denoted as $C^1_{-1}(\mathbb{R}^3_+)$.

**Definition 5.2** Let $U(x,y,z)$ be a scalar field that belongs to the set $C^1_{-1}(\mathbb{R}^3_+)$.
Then the general fractional gradient $Grad^{(K)}_W$ for the region $W = \mathbb{R}^3_{0,+}$ is defined as
$$(Grad^{(K)}_W U)(x,y,z) =$$

$$e_x D^{x,*}_{(K_1)}[x']U(x',y,z) + e_y D^{y,*}_{(K_2)}[y']U(x,y',z) + e_z D^{z,*}_{(K_3)}[z']U(x,y,z').$$
This operator will be called the regional general fractional gradient (regional GF gradient).

**Remark 5.1** The formula defining the operator can be written in compact form. If the scalar field $U(x,y,z)$ belongs to the function space $C^1_{-1}(\mathbb{R}^3_{0,+})$, then the general fractional gradient for the region $W = \mathbb{R}^3_+$ is defined as
$$(Grad^{(K)}_W U)(x,y,z) = (D^{(K)}_W U)(x,y,z) = \sum_{k=1}^3 e_k D^{x_k,*}_{(K_k)}[x'_k] F_k,$$
where
$$D^{(K)}_W := e_x D^{x,*}_{(K_1)}[x'] + e_y D^{y,*}_{(K_2)}[y'] + e_z D^{z,*}_{(K_3)}[z'].$$

**Remark 5.2** The general fractional gradients can be defined not only for $W = \mathbb{R}^3_{0,+}$, but also for regions $W \subset \mathbb{R}^3_{0,+}$, surfaces $S \subset \mathbb{R}^3_{0,+}$ and line $L \subset \mathbb{R}^3_{0,+}$.

The gradient theorem is very important for the vector calculus and its generalizations, since this theorem is actually a fundamental theorem for standard vector calculus for line integral and gradient, and their generalizations.

In the following sections, we will analyze the differences between regional and line general fractional gradients to formulate general gradient theorems.

Note that the general fractional gradient for line $L \subset \mathbb{R}^3_{0,+}$ allows us to prove the general fractional gradient theorem for a wider class of lines $\mathbb{R}^3_{0,+}$ and $\mathbb{R}^2_{0,+}$.

## 5.2 Difficulties in Generalization of Gradient Theorem

Let $L$ be simple line in $\mathbb{R}^3_{0,+}$ that is described by the equations $y = y(x)$ and $z = z(x)$ for $x \in [a,b]$. Since the line is simple, the derivatives of the functions $y = y(x)$ and $z = z(x)$ are not equal to zero and do not change sign on the interval $[a,b]$. By definition, the linear integral of a vector field $F \in C^1(\mathbb{R}^3_+)$ for the simple line $L$ can be given by the equation

$$\int_L (F, dr) = \int_a^b F_x(x, y(x), z(x)) \, dx +$$

$$\int_c^d F_y(x(y), y, z(y)) \, dy + \int_e^f F_z(x(z), y(z), z) \, dz.$$

The standard gradient theorem is proved by using this definition of the line integral and the standard chain rule.

For the vector field $F = grad U$, we have

$$\int_L (grad U, dr) = \int_a^b U_x^{(1)}(x, y(x), z(x)) \, dx +$$

$$\int_c^d U_y^{(1)}(x(y), y, z(y)) \, dy + \int_e^f U_z^{(1)}(x(z), y(z), z) \, dz =$$

$$\int_c^d U_y^{(1)}(x(y), y, z(y)) \, dy + \int_e^f U_z^{(1)}(x(z), y(z), z) \, dz =$$

$$\int_a^b (U_x^{(1)}(x, y(x), z(x)) + U_y^{(1)}(x, y(x), z(x)) y_x^{(1)}(x) + U_z^{(1)}(x, y(x), z(x)) z_x^{(1)}(x)) \, dx.$$

Using the standard chain rule

$$\frac{dU(x, y(x), z(x))}{dx} = \left(\frac{\partial U(x, y, z)}{\partial x}\right)_{y=y(x), z=z(x)} +$$

$$\left(\frac{\partial U(x, y, z)}{\partial y}\right)_{y=y(x), z=z(x)} \frac{dy(x)}{dx} + \left(\frac{\partial U(x, y, z)}{\partial z}\right)_{y=y(x), z=z(x)} \frac{dz(x)}{dx},$$

we get

$$\int_L (grad U, dr) = \int_a^b \frac{dU(x, y(x), z(x))}{dx} =$$

$$U(b, y(b), z(b)) - U(a, y(a), z(a)) = U(b, d, f) - U(a, c, e).$$

In the fractional calculus and GFC, the standard chain rule is violated.

The line GFI for the simple line of the vector field $F$ is defined by the equation

$$(I_L^{(M)}, F) = I_{[a,b]}^{(M_1)}[x] F_x(x, y(x), z(x)) +$$

$$I_{[c,d]}^{(M_2)}[y] F_y(x(y), y, z(y)) + I_{[e,f]}^{(M_3)}[z] F_z(x(z), y(z), z).$$

For the vector field $F = Grad_W^{(K)} U$, where

$$F_x(x, y, z) = (Grad_W^{(K)} U)_x(x, y, z) = D_{(K_1)}^{x,*}[x'] U(x', y, z),$$

$$F_y(x, y, z) = (Grad_W^{(K)} U)_y(x, y, z) = D_{(K_2)}^{y,*}[y'] U(x, y', z),$$

$$F_y(x, y, z) = (Grad_W^{(K)} U)_y(x, y, z) = D_{(K_3)}^{z,*}[z'] U(x, y, z'),$$

we have

$$(I_L^{(M)}, F) = I_{[a,b]}^{(M_1)}[x](D_{(K_1)}^{x,*}[x'] \, U(x',y,z))_{y=y(x),z=z(x)} +$$

$$I_{[c,d]}^{(M_2)}[y](D_{(K_2)}^{y,*}[y'] \, U(x,y',z))_{x=x(y),z=z(y)} + I_{[e,f]}^{(M_3)}[z](D_{(K_3)}^{z,*}[z'] \, U(x,y,z'))_{x=x(z),y=y(z)}.$$

In this case, we should be emphasized that the fundamental theorem of GFC cannot be used. This fact is based on the following inequalities

$$I_{[a,b]}^{(M_1)}[x](D_{(K_1)}^{x,*}[x'] \, U(x',y,z))_{y=y(x),z=z(x)} \neq I_{[a,b]}^{(M_1)}[x]D_{(K_1)}^{x,*}[x'] \, (U(x',y,z))_{y=y(x),z=z(x)},$$

$$I_{[a,b]}^{(M_1)}[x](D_{(K_1)}^{x,*}[x'] \, U(x',y,z))_{y=y(x),z=z(x)} \neq (I_{[a,b]}^{(M_1)}[x]D_{(K_1)}^{x,*}[x'] \, U(x',y,z))_{y=y(x),z=z(x)}.$$

which are satisfied if the functions $y = y(x)$, $z = z(x)$ depend on the coordinate $x$. Similarly, the fundamental theorem does not hold for other variables.

As a result, for the regional GF gradient $Grad_W^{(K)}$ with $W = \mathbb{R}_{0,+}^3$, the general fractional theorem can be proved only if the line consists of sections parallel to the axes.

## 5.3 General Fractional Gradient Theorem for Regional GF Gradient

The gradient theorem can be considered as a fundamental theorem of standard calculus for line integrals and gradients.

Let us consider a line $L$ that consists of line segments that are parallel to the axes. To prove the general gradient theorem for such a broken line (polygonal chain), it is convenient to use the concept of an elementary broken line. We can state that any continuous line $L$ that consists only of lines parallel to the axes can be represented as a sequence of elementary lines $L_k$.

An elementary broken line is a line consisting of three (no more than three) segments parallel to different axes. There are 48 such elementary lines, which eight differ in the directions of the segments along the axis or against the axis, and six different different orders.

$$I_{L_k}^{(M)}[x,y,z] := sgn(x_{k+1} - x_k) \, e_x \, I_{\omega[x_k,x_{k+1}]}^{(M_1)}[x] + \qquad (79)$$

$$sgn(y_{k+1} - y_k) \, e_y \, I_{\omega[y_k,y_{k+1}]}^{(M_2)}[y] + sgn(z_{k+1} - z_k) \, e_z \, I_{\omega[z_k,z_{k+1}]}^{(M_3)}[z],$$

where $\omega[x_k, x_{k+1}]$ is defined by (60) and $sgn(x)$ is defined by equation (59).

Let us prove the following general gradient theorem for the regional GF gradient $Grad_W^{(K)}$ with $W = \mathbb{R}_{0,+}^3$.

**Theorem 5.1** *(General Gradient Theorem for Regional GF Gradient)*

*Let $L$ be continuous line in $\mathbb{R}_{0,+}^3$ that consists only of lines parallel to the axes can be represented as a sequence of elementary lines, for which $A_k(a,c,e)$ is the initial point and $D(b,d,f)$ is the final point.*

*Let $U(x,y,z)$ be a scalar field that belongs to the set $C_{-1}^1(\mathbb{R}_+^3)$.*

*Then the line GFI for the vector field $Grad_W^{(K)}U$ with $W = \mathbb{R}_{0,+}^3$ satisfies the equation*

$$(I_L^{(M)}, Grad_W^{(K)}U) = U(b,d,f) - U(a,c,e).$$

*Proof.* For simplicity, we will consider only an elementary line $L_k$, when moving along which from the initial point $A_k(a_k,c_k,e_k)$ to the final point $D_k(b_k,d_k,f_k)$ coordinate values do not decrease, and with order $XYZ$. Let the line $L_k$ be represented by the points
$$A_k(a_k, c_k, e_k), \quad B_k(b_k, c_k, e_k), \quad C_k(b_k, d_k, e_k), \quad D(b_k, d_k, f_k).$$
For the vector field $F \in \mathbb{F}^0_{-1,L}(\mathbb{R}^3_+)$, for which
$$F_x(x, c_k, e_k), F_y(b_k, y, e_k), F_z(b_k, d_k, z) \in C_{-1}(0, \infty),$$
the line GFI for the elementary line $L_k = A_k B_k C_k D_k$ is defined by the equation
$$(I_{L_k}^{(M)}, F) =$$

$$I_{A_k B_k}^{(M_1)}[x] F_x(x, c_k, e_k) + I_{B_k C_k}^{(M_2)}[y] F_y(b_k, y, e_k) + I_{C_k D_k}^{(M_3)}[z] F_z(b_k, d_k, z) =$$

$$I_{[a_k,b_k]}^{(M_1)}[x] F_x(x, c_k, e_k) + I_{[c_k,d_k]}^{(M_2)}[y] F_y(b_k, y, e_k) + I_{[e_k,f_k]}^{(M_3)}[z] F_z(b_k, d_k, z).$$

Let us condider the vector field
$$F = \text{Grad}_W^{(K)} U,$$
where $W = \mathbb{R}^3_{0,+}$, and
$$F_x(x, y, z) = D_{(K_1)}^{x,*}[x'] U(x', y, z),$$

$$F_y(x, y, z) = D_{(K_2)}^{y,*}[y'] U(x, y', z),$$

$$F_z(x, y, z) = D_{(K_3)}^{z,*}[z'] U(x, y, z').$$
Then line GFI for the elementary line $L_k = A_k B_k C_k D_k$ is given as
$$(I_{L_k}^{(M)}, \text{Grad}_W^{(K)} U) = I_{[a_k,b_k]}^{(M_1)}[x] D_{(K_1)}^{x,*}[x'] U(x', c_k, e_k) +$$

$$I_{[c_k,d_k]}^{(M_2)}[y] D_{(K_2)}^{y,*}[y'] U(b_k, y', e_k) + I_{[e_k,f_k]}^{(M_3)}[z] D_{(K_3)}^{z,*}[z'] U(b_k, d_k, z').$$
Using the fundamental theorem fo GFC, we obtain
$$(I_{L_k}^{(M)}, \text{Grad}_W^{(K)} U) = (U(b_k, c_k, e_k) - U(a_k, c_k, e_k)) +$$

$$U(b_k, d_k, e_k) - U(b_k, c_k, e_k)) + (U(b_k, d_k, f_k) - U(b_k, d_k, e_k)) =$$

$$U(b_k, d_k, f_k) - U(a_k, c_k, e_k).$$
Equations for other elementary lines are proved in a similar way.

Using that any continuous line $L$ that consists only of lines parallel to the axes can be represented as a sequence of elementary lines $L_k$, we have that general gradient theorem holds.

**Corollary 5.1** *The line GFI of the vector field $F$ is independent on the path, which is described by the lines $L$ that consist only of line segments parallel to the axes, if the vecror field $F$ can be represented as the regional GF gradient of a function $U(x, y, z) \in C^1_{-1}(\mathbb{R}^3_+)$.*

### 5.4 Line General Fractional Gradient in $\mathbb{R}^2_{0,+}$

The general fractional gradients can be defined not only for the region $W = \mathbb{R}^3_{0,+}$. Using the fact that GFD is integro-differential operator, we can define general fractional gradients for line $L \subset \mathbb{R}^2_{0,+}$ and $L \subset \mathbb{R}^3_{0,+}$.

Let us define a general fractional gradient for a simple line on the $XY$-plane. In this definition of line general fractional gradient (line GF Gradienr), we can use the fact that GFD can be represented as a sequential action of a first-order derivative and a general fractional integral:
$$D^{x,*}_{(K_1)}[x'] = I^{x,*}_{(K_1)}[x'] \frac{\partial}{\partial x'}.$$
The general fractional vector operators can also be defined as a sequential action of first-order derivatives and general fractional integrals.

**Definition 5.3** Let $L$ be a simple line in $\mathbb{R}^2_{0,+}$ of the XY-plane, and $U(x,y) \in \mathbb{F}^1_{-1,L}(\mathbb{R}^2_+)$ that means
$$U^{(1)}_x(x,y(x)) \in C_{-1}(0,\infty), \quad U^{(1)}_y(x(y),y) \in C_{-1}(0,\infty).$$
Then the line general fractional gradient for the line $L$ is defined by the equation
$$Grad^{(K)}_L U(x,y) = D^{(K)}_L U(x,y) :=$$

$$e_x I^{x,*}_{(K_1)}[x'] U^{(1)}_{x'}(x',y(x')) + e_y I^{y,*}_{(K_2)}[y'] U^{(1)}_{y'}(x(y'),y') =$$

$$e_x I^{x,*}_{(K_1)}[x'] (\frac{\partial U(x',y)}{\partial x'})_{y=y(x')} + e_y I^{y,*}_{(K_2)}[y'] (\frac{\partial U(x,y')}{\partial y'})_{x=x(y')},$$
where the kernels $(M_1(x), K_1(x))$ and $(M_2(y), K_2(y))$ belong to the Luchko set $\mathbb{L}$.

We can use the definitions of the line GFI and line GF Gradient to prove the following theorem.

**Theorem 5.2** *(General Fractional Gradient Theorem for line GF Gradient)*
Let $L$ be a simple line in $\mathbb{R}^2_{0,+}$ of the XY-plane, which connects the points $A(a,b)$ and $B(c,d)$, and the scalar field $U(x,y)$ belongs to the set $\mathbb{F}^1_{-1,L}(\mathbb{R}^2_+)$.
Then, the equality
$$(I^{(M)}_L, Grad^{(K)}_L U) = U(c,d) - U(a,b)$$
holds.

*Proof.* Using that the line GFI of a vector field $F$ for the simple line $L \subset \mathbb{R}^2_{0,+}$ of the XY-plane is defined by the equation
$$(I^{(M)}_L, F) = I^{(M_1)}_{[a,b]}[x] F_x(x,y(x)) + I^{(M_2)}_{[c,d]}[y] F_y(x(y),y),$$
we can consider the line GFI of the vector field that is defined by the line general fractional gradient
$$F_x(x,y(x)) := (Grad^{(K)}_L U(x,y))_x = I^{x,*}_{(K_1)}[x'] U^{(1)}_{x'}(x',y(x')),$$

$$F_y(x(y),x) := (Grad^{(K)}_L U(x,y))_y = I^{y,*}_{(K_2)}[y'] U^{(1)}_{y'}(x(y'),y').$$
Therefore, we can get the line GFI of the line FG Gradient
$$(I^{(M)}_L, Grad^{(K)}_L U) =$$

$$I^{(M_1)}_{[a,b]}[x] I^{x,*}_{(K_1)}[x'] U^{(1)}_{x'}(x',y(x')) + I^{(M_2)}_{[c,d]}[y] I^{y,*}_{(K_2)}[y'] U^{(1)}_{y'}(x(y'),y').$$
Then, we can use the fact that the kernels $(M_1(x), K_1(x))$ and $(M_2(y), K_2(y))$ belong to the Luchko set $\mathbb{L}$. In this case, we have the property
$$I^{(M_1)}_{[a,b]}[x] I^{x,*}_{(K_1)}[x'] f(x') = (M_1 * (K_1 * f))(b) - (M_1 * (K_1 * f))(a) =$$

$$((M_1 * K_1) * f)(b) - ((M_1 * K_1) * f)(a) = (\{1\} * f)(b) - (\{1\} * f)(a) =$$

$$\int_0^b f(x)\,dx - \int_0^a f(x)\,dx = \int_a^b f(x)\,dx.$$

Similarly, we obtain the equation

$$I_{[c,d]}^{(M_2)}[y]\,I_{(K_2)}^{y,*}[y']\,g(y') = (\{1\} * g)(d) - (\{1\} * g)(c) = \int_c^d g(y)\,dy.$$

Therefore, we get

$$I_{[a,b]}^{(M_1)}[x]\,I_{(K_1)}^{x,*}[x']\,U_{x'}^{(1)}(x',y(x')) + I_{[c,d]}^{(M_2)}[y]\,I_{(K_2)}^{y,*}[y']\,U_{y'}^{(1)}(x(y'),y') =$$

$$\int_a^b \left(\frac{\partial U(x',y)}{\partial x'}\right)_{y=y(x')} dx' + \int_c^d \left(\frac{\partial U(x,y')}{\partial y'}\right)_{x=x(y')} dy' =$$

$$\int_a^b \left(\frac{\partial U(x',y)}{\partial x'}\right)_{y=y(x')} dx' + \int_a^b \left(\frac{\partial U(x',y)}{\partial y}\right)_{y=y(x')} \frac{\partial y(x')}{\partial x'}\,dx' =$$

$$\int_a^b \frac{dU(x,y(x))}{dx}\,dx = U(c,d) - U(a,b),$$

where we use the standard gradient theorem.

As a result, we proved the general fractional gradient theorem for the line GF Gradient with the simple line in $\mathbb{R}_{0,+}^2$.

## 5.5 Line General Fractional Gradient in $\mathbb{R}_{0,+}^3$

The proposed approach to define line GF Gradien for lines in $\mathbb{R}_{0,+}^2$ can be used to define line GF Gradient for lines in $\mathbb{R}_{0,+}^3$. Let $L$ be a simple line in $\mathbb{R}_{0,+}^3$, which is defined in form (46), (47), (48), and the vector field $F(x,y,z)$ belongs to the sset $\mathbb{F}_{-1}^0(L)$. Then, the line GFI over a simple line $L$ in $\mathbb{R}_{0,+}^3$ can be defined by equation (56). The standard line integral of second kind over a simple line $L$ in $\mathbb{R}_{0,+}^3$ can be defined by the equation (55).

**Definition 5.4** Let $L$ be a simple line in $\mathbb{R}_{0,+}^3$, which is defined in form (46), (47), (48), and a scalar field $U(x,y,z)$ satisfies the conditions

$$U_1(x) := U_x^{(1)}(x,y(x),z(x)) = \left(\frac{\partial U(x,y,z)}{\partial x}\right)_{y=y(x),z=z(x)} \in C_{-1}^n(0,\infty),$$

$$U_2(y) := U_y^{(1)}(x(y),y,z(y)) = \left(\frac{\partial U(x,y,z)}{\partial y}\right)_{x=x(y),z=z(y)} \in C_{-1}^n(0,\infty),$$

$$U_3(z) := U_z^{(1)}(x(z),y(z),z) = \left(\frac{\partial U(x,y,z)}{\partial z}\right)_{x=x(z),y=y(z)} \in C_{-1}^n(0,\infty).$$

Then the set of such scalar fields $U(x,y,z)$ will be denoted as $U(x,y,z) \in \mathbb{F}_{-1}^n(L)$.

**Definition 5.5** Let $L$ be a simple line in $\mathbb{R}_{0,+}^3$, which is defined in form (46), (47), (48), and a scalar field $U(x,y,z)$ belongs to the set $\mathbb{F}_{-1}^1(L)$.

Then the line general fractional gradient for the line $L \in \mathbb{R}_{0,+}^3$ is defined by the equation

$$(Grad_L^{(K)} U)(x,y,z) = (D_L^{(M)} U)(x,y,z) = I_{[a,b]}^{(M_1)}[x]\,U_x^{(1)}(x,y(x),z(x)) +$$

$$I_{[c,d]}^{(M_2)}[y] \, U_y^{(1)}(x(y), y, z(y)) + I_{[e,f]}^{(M_3)}[z] \, U_z^{(1)}(x(z), y(z), z). \tag{80}$$

where the pairs of the kernels $(M_1(x), K_1(x))$, $(M_2(y), K_2(y))$ and $(M_3(z), K_3(z))$ belong to the Luchko set $\mathbb{L}$.

Let us prove the general fractional gradient theorem for line GF Gradient with simple lines.

**Theorem 5.3** *(General Fractional Gradient Theorem for Line GF Gradient)*

Let $L$ be a simple line in $\mathbb{R}_{0,+}^3$, which is defined in form (46), (47), (48), and connects the points $A(a, c, e)$ and $B(b, d, f)$.

Let $U(x, y, z)$ be a scalar filed that belongs to the set $\mathbb{F}_{-1,L}^1(\mathbb{R}_+^3)$.

Then, the equality

$$(I_L^{(M)}, Grad_L^{(K)} U) = U(a, c, e) - U(b, d, f)$$

holds.

*Proof.* Using that the line GFI of a vector field $F$ for the simple line $L \subset \mathbb{R}_{0,+}^3$ is defined by the equation

$$(I_L^{(M)}, F) = I_{[a,b]}^{(M_1)}[x] \, F_x(x, y(x), z(x)) +$$

$$I_{[c,d]}^{(M_2)}[y] \, F_y(x(y), y, z(y)) + I_{[e,f]}^{(M_3)}[z] \, F_z(x(z), y(z), z),$$

we can consider the line GFI of the vector field that is defined by the line general fractional gradient

$$F_x(x, y(x), z(x)) := (Grad_L^{(K)} U(x, y, z))_x = I_{(K_1)}^{x,*}[x'] \, U_{x'}^{(1)}(x', y(x'), z(x')),$$

$$F_y(x(y), y, z(y)) := (Grad_L^{(K)} U(x, y, z))_y = I_{(K_2)}^{y,*}[y'] \, U_{y'}^{(1)}(x(y'), y', z(y')).$$

$$F_z(x(z), y(z), z) := (Grad_L^{(K)} U(x, y, z))_z = I_{(K_3)}^{z,*}[z'] \, U_{z'}^{(1)}(x(z'), y(z'), z').$$

Therefore, we can get the line GFI of the line FG Gradient

$$(I_L^{(M)}, Grad_L^{(K)} U) = I_{[a,b]}^{(M_1)}[x] \, I_{(K_1)}^{x,*}[x'] \, U_{x'}^{(1)}(x', y(x'), z(x')) +$$

$$I_{[c,d]}^{(M_2)}[y] \, I_{(K_2)}^{y,*}[y'] \, U_{y'}^{(1)}(x(y'), y', z(y')) +$$

$$I_{[e,f]}^{(M_3)}[z] \, I_{(K_3)}^{z,*}[z'] \, U_{z'}^{(1)}(x(z'), y(z'), z').$$

Then, we can use the fact that the kernels $(M_1(x), K_1(x))$, $(M_2(y), K_2(y))$ and $(M_3(x), K_3(x))$ belong to the Luchko set $\mathbb{L}$. In this case, we have the property

$$I_{[a,b]}^{(M_1)}[x] \, I_{(K_1)}^{x,*}[x'] f(x') = (M_1 * (K_1 * f))(b) - (M_1 * (K_1 * f))(a) =$$

$$((M_1 * K_1) * f)(b) - ((M_1 * K_1) * f)(a) = (\{1\} * f)(b) - (\{1\} * f)(a) =$$

$$\int_0^b f(x) \, dx - \int_0^a f(x) \, dx = \int_a^b f(x) \, dx.$$

Similarly, we obtain the equations

$$I_{[c,d]}^{(M_2)}[y] \, I_{(K_2)}^{y,*}[y'] g(y') = (\{1\} * g)(d) - (\{1\} * g)(c) = \int_c^d g(y) \, dy,$$

$$I_{[e,f]}^{(M_3)}[z]\, I_{(K_3)}^{z,*}[z']\, h(z') = (\{1\} * h)(f) - (\{1\} * h)(e) = \int_e^f h(z)\, dz.$$

Therefore, we get

$$I_{[a,b]}^{(M_1)}[x]\, I_{(K_1)}^{x,*}[x']\, U_{x'}^{(1)}(x', y(x'), z(x')) +$$

$$I_{[c,d]}^{(M_2)}[y]\, I_{(K_2)}^{y,*}[y']\, U_{y'}^{(1)}(x(y'), y', z(y')) +$$

$$I_{[e,f]}^{(M_3)}[z]\, I_{(K_3)}^{z,*}[z']\, U_{z'}^{(1)}(x(z'), y(z'), z') =$$

$$\int_a^b \left(\frac{\partial U(x',y,z)}{\partial x'}\right)_{y=y(x'),z=z(x')} dx' + \int_c^d \left(\frac{\partial U(x,y',z)}{\partial y'}\right)_{x=x(y'),z=z(y')} dy' +$$

$$\int_e^f \left(\frac{\partial U(x,y,z')}{\partial z'}\right)_{x=x(z'),y=y(z')} dz' =$$

$$\int_a^b \left(\frac{\partial U(x',y,z)}{\partial x'}\right)_{y=y(x'),z=z(x')} dx' + \int_a^b \left(\frac{\partial U(x,y',z)}{\partial y'}\right)_{y=y(x'),z=z(x')} \frac{\partial y(x')}{\partial x'} dx' +$$

$$\int_a^b \left(\frac{\partial U(x,y,z')}{\partial z'}\right)_{y=y(x'),z=z(x')} \frac{\partial z(x')}{\partial x'} dx' =$$

$$\int_a^b \frac{dU(x,y(x),z(x))}{dx}\, dx = U(b,d,f) - U(a,c,e),$$

where we use the standard gradient theorem.

As a result, we proved the general fractional gradient theorem for the line GF Gradient with the simple line in $\mathbb{R}^3_{0,+}$.

**Corollary 5.2** *The line GFI of the vector field $F$ is independent on the path, which is described by the lines $L$ that consist of simple lines, if the vecror field $F$ can be represented as the line GF Gradient of a function $U(x,y,z) \in \mathbb{F}^1_{-1}(L)$.*

**Remark 5.3** The general fractional gradient theorem for line GF Gradient can be proves for the YZ-simple line, which can be considered as a union of XZ-simple lines and XY-simple lines.

# 6 General Fractional Green Theorem

It is known that the standard Green's theorem gives the relationship between a line integral around a simple closed curve $\partial S$ and a double integral over the plane region $S$ bounded by $\partial S$. Let us formulate a generalization of the Green's theorem for the general fractional vector calculus.

## 6.1 General Fractional Green Theorem for X and Y Simple Regions

Let us define a simple region on the plane.

**Definition 6.1** *Let $S$ be a region in $\mathbb{R}^2_{0,+}$ of the XY-plane that is bounded above and below by smooth lines $L_{1,x}$, $L_{2,x}$, and lines $L_{1,y}$, $L_{2,y}$, which are parallel to the Y-axis. We will assume that the lines $L_{1,x}$, $L_{2,x}$ are Y-simple lines, which are described by the equations*
$$L_{1,x} = \{(x,y): \ x \in [a,b], \ y = y_1(x) \geq 0\},$$

$$L_{2,x} = \{(x,y): \ x \in [a,b], \ y = y_2(x) \geq 0\},$$
*where the functions $y_1(x)$ and $y_3(x)$ belong to the function space $C^1(a,b)$ in the closed domain $[a,b]$ that is a projection of the region $S$ onto the X-axis, and $y_2(x) \geq y_1(x)$ for all $x \in [a,b]$.*

*The lines $L_{1,y}$ and $L_{2,y}$ parallel to the Y axis are described as*
$$L_{1,y} = \{(x,y): \ x = a \geq 0, \ 0 \leq y_1(a) \leq y \leq y_2(a)\},$$

$$L_{2,y} = \{(x,y): \ x = b \geq 0, \ 0 \leq y_1(b) \leq y \leq y_2(b)\}.$$
*If the conditions $y_1(a) = y_2(a)$, $y_1(b) = y_2(b)$ are satisfied, then the lines $L_{1,y}$ and $L_{2,y}$ are absent.*

*Then, such regions $S$ on the plane will be called the Y-simple region (simple area along the Y-axis).*

**Definition 6.2** *Let $S$ be a region in $\mathbb{R}^2_{0,+}$ of the XY-plane. The region $S \in \mathbb{R}^2_{0,+}$ is called simple, if $S$ is simple along $X$ and $Y$ axes. The region $S \in \mathbb{R}^2_{0,+}$ is called piecewise simple, if $S$ can be divided into a finite number of simle regions with respect to these two axes.*

Let us consider the vector field
$$F(x,y) = e_x F_x(x,y) + e_y F_y(x,y),$$
for the simple region $S$, where the functions $F_x = F_x(x,y)$ and $F_y = F_y(x,y)$ satisfy the conditions
$$F_x(x, y_1(x)), F_x(x, y_2(x)) \in C_{-1}(0, \infty), \tag{81}$$

$$F_y(x_1(y), y), F_y(x_2(y), y) \in C_{-1}(0, \infty), \tag{82}$$
and
$$F_x(\xi, y), F_y(x, \xi) \in C^1_{-1}(0, \infty) \quad \text{for each} \quad \xi \in [0, \infty). \tag{83}$$
In this case, we will use notation $F(x,y) \in \mathbb{F}^1_{-1}(S)$ for $S \subset \mathbb{R}^2_+$.

**Theorem 6.1** *(General Fractional Green's Theorem for X and Y Simple Region)*

Let $S$ be a $Y$-simple region in $\mathbb{R}^2_{0,+}$ of the $XY$-plane that can be described as

$$S := \{(x,y): \ 0 \le a \le x \le b, \ \ 0 \le y_1(x) \le y \le y_2(x)\}, \qquad (84)$$

where $y_1(x)$ and $y_2(x)$ belong to the space $C(a,b)$.

Let also region $S$ be represented by the union of the finite number of $X$-simple regions

$$S = \bigcup_{k=1}^{n} S_k \qquad (85)$$

of a finite number $X$-simple regions

$$S_k := \{(x,y): \ 0 \le x_{1,k}(y) \le x \le x_{2,k}(y), \ \ 0 \le c_k \le y \le d_k)\}, \qquad (86)$$

where $x_{1,k}(y)$ and $x_{2,k}(y)$ belong to the space $C(a,b)$ for all $k = 1, \dots n$.

Let the vector field

$$F(x,y) = e_x F_x(x,y) + e_y F_y(x,y),$$

satisfy the conditions

$$F_x(x, y_1(x)), F_x(x, y_2(x)) \in C_{-1}(0, \infty), \qquad (87)$$

$$F_y(x_{1,k}(y), y), F_y(x_{2,k}(y), y) \in C_{-1}(0, \infty), \qquad (88)$$

and

$$F_x(\xi, y), F_y(x, \xi) \in C^1_{-1}(0, \infty) \quad for each \quad \xi \in [0, \infty). \qquad (89)$$

Then, the equation

$$I_{\partial S}^{(M)}[x] F_x(x,y) + I_{\partial S}^{(M)}[y] F_y(x,y) =$$

$$I_S^{(M)}[x, y] \left( D_{(K_1)}^{x,*}[x'] F_y(x', y) - D_{(K_2)}^{y,*}[y'] F_x(x, y') \right) \qquad (90)$$

holds.

*Proof.* Let us consider the $Y$-simple region $S$ on the $XY$-plane that is described as

$$S := \{(x,y): \ 0 \le a \le x \le b, \ \ 0 \le y_1(x) \le y \le y_2(x)\}, \qquad (91)$$

where $y_1(x)$ and $y_2(x)$ are continuous functions on $[a,b]$. The boundary $L = \partial S$ of this area can be divided into the lines $AB$, $BC$, $CD$, $DA$, where the coordinates of the points are the following:

$$A(a, y_1(a)), \quad B(b, y_1(b)), \quad C(b, y_2(b)), \quad D(a, y_2(a)).$$

The lines $AB$ and $CD$ are $Y$-simple lines ($L_{1x} = AB$, $L_{2x} = CD$), and $BC$, $DA$ are lines parallel to the $Y$ axis.

Let us consider the vector field $F$ as a sum of two vectors $F_x$ and $F_y$ in the form

$$F(x,y) = F_x(x,y) + F_y(x,y): \ F_x = e_x F_x(x,y), \ F_y = e_y F_y(x,y),$$

where the functions $F_x(x,y)$ and $F_r(x,y)$ satisfy conditions (87), (88), (89).

1) The line GFI for the vector field $F_x = e_x F_x(x,y)$ and the line $L = \partial S$ is given by the equation

$$(I_L^{(M)}[x,y], F_x(x,y)) = I_{AB}^{(M_1)}[x] F_x(x, y_1(x)) +$$

$$I_{BC}^{(M_1)}[x] F_x(b, y) + I_{CD}^{(M_1)}[x] F_x(x, y_2(x)) + I_{DA}^{(M_1)}[x] F_x(a, y) =$$

$$I_{AB}^{(M_1)}[x] F_x(x, y_1(x)) - I_{DC}^{(M_1)}[x] F_x(x, y_2(x)) =$$

$$I^{(M_1)}_{[a,b]}[x]\,(F_x(x,y_1(x)) - F_x(x,y_2(x))),$$

where we use

$$I^{(M_1)}_{BC}[x]\,F_x(b,y) = 0, \quad I^{(M_1)}_{DA}[x]\,F_x(a,y) = 0,$$

since $dx = 0$ for $x = b$ and $x = a$.

Using the second fundamental theorem of GFC in the form

$$F_x(x,y_1(x)) - F_x(x,y_2(x)) = -I^{(M_2)}_{[y_1(x),y_2(x)]}[y]\,D^{y,*}_{(K_2)}[y]\,F_x(x,y),$$

where we assume that conditions (87) and (89) are satisfied. Then, we obtain

$$(I^{(M)}_L[x,y], F_x(x,y)) = -I^{(M_1)}_{[a,b]}[x]\,I^{(M_2)}_{[y_1(x),y_2(x)]}[y]\,D^{y,*}_{(K_2)}[y]\,F_x(x,y).$$

As a result, we obtain

$$(I^{(M)}_L[x,y], F_x(x,y)) = -I^{(M)}_S[x,y]\,D^{y,*}_{(K_2)}[y]\,F_x(x,y),$$

where we use the definition of double GFI by the iterated (repeated) GFI

$$I^{(M)}_S[x,y] = I^{(M_1)}_{[a,b]}[x]\,I^{(M_2)}_{[y_1(x),y_2(x)]}[y].$$

2) Let us consider the vector field $F_y = e_y\,F_y(x,y)$ and the $Y$-simple region $S$ on the $XY$-plane, which can be described by (85) and (86). Using the similar transformations and equations for the line GFI for $L = \partial S$, we obtain

$$(I^{(M)}_L[x,y], F_y(x,y)) = \sum_{k=1}^n \sum_{l=1}^m (I^{(M)}_{L_{kl}}[x,y], F_y(x,y)) =$$

$$\sum_{k=1}^n \sum_{l=1}^m I^{(M_2)}_{[c_{kl},d_{kl}]}[y]\,I^{(M_1)}_{[x_{1,kl}(y),x_{2,kl}(y)]}\,F_y(x,y) = I^{(M)}_S[x,y]\,D^{x,*}_{(K_1)}[x]\,F_y(x,y),$$

where we use the definition of double GFI by the iterated (repeated) GFI

$$I^{(M)}_S[x,y] = \sum_{k=1}^n \sum_{l=1}^m I^{(M)}_{S_{kl}}[x,y] = \sum_{k=1}^n \sum_{l=1}^m I^{(M_2)}_{[c_{kl},d_{kl}]}[y]\,I^{(M_1)}_{[x_{1,kl}(y),x_{2,kl}(y)]}[x].$$

3) Therefore, the Green formula for the vector field $F = F_x + F_y$ takes the form

$$(I^{(M)}_L[x,y], F(x,y)) = I^{(M)}_S[x,y]\,(D^{x,*}_{(K_1)}[x]\,F_y(x,y) - D^{y,*}_{(K_2)}[y]\,F_x(x,y)),$$

if conditions (87), (88) and (89) are satisfied.

**Remark 6.1** Using the proof of the theorem, we see that it is not required to use an approach similar to the approach that is proposed for the linearly general fractional gradient. WE should also note that this approach should be used for the general fractional Sokes theorem.

**Remark 6.2** Using the proof, we can see that the general fractional Green's theorem can be proved for the region $S$ that is union of funite mumber of $Y$-simple regions $S_j$ ($j = 1, \ldots, m$), each of which can be represented as a union of a finite number of simple regions $S_{jk}$, $k = 1, \ldots, n$.

## 6.2 General Fractional Green Theorem for Simple Region

Let us consider a simple region $S \subset \mathbb{R}^2_{0,+}$ of the $XY$-plane that can be described as

$$S := \{(x,y): \; 0 \le a \le x \le b, \; 0 \le y_1(x) \le y \le y_2(x)\}, \tag{92}$$

$$S := \{(x,y): \; 0 \le x_1(y) \le x \le x_2(y), \; 0 \le c \le y \le d)\}, \tag{93}$$

where $y_1(x)$ and $y_2(x)$ belong to the space $C^1(a,b)$, the derivatives of these functions are nonzero and do not change sign on the interval $[a,b]$. Then the functions $x_1(y)$ and $x_2(y)$ exit, belong to the space $C^1(c,d)$, where $c = x_1(a) = x_2(a)$ and $d = x_1(b) = x_2(b)$. In this case, we will use notation $S \in G^1(\mathbb{R}^2_{0,+})$.

The boundary of the simple region $S \in G^1(\mathbb{R}^2_{0,+})$ is a piecewise simple closed line $\partial S$ in the $XY$-plane that can be described as the union of two $Y$-simple lines
$$\partial S := \{(x,y): x \in [a,b], y = y_1(x)\} \cup \{(x,y): x \in [a,b], y = y_2(x)\}, \quad (94)$$
or as the union of two $X$-simple lines:
$$\partial S := \{(x,y): y \in [c,d], x = x_1(y)\} \cup \{(x,y): y \in [c,d], x = x_2(y)\}. \quad (95)$$

For $S \in G^1(\mathbb{R}^2_{0,+})$, we will consider the vector field
$$F(x,y) = e_x F_x(x,y) + e_y F_y(x,y),$$
where the functions $F_x = F_x(x,y)$ and $F_y = F_y(x,y)$ satisfy the conditions
$$F_x(x,y_1(x)), F_x(x,y_2(x)) \in C_{-1}(0,\infty), \quad (96)$$

$$F_y(x_1(y),y), F_y(x_2(y),y) \in C_{-1}(0,\infty), \quad (97)$$
and
$$F_x(\xi,y), F_y(x,\xi) \in C^1_{-1}(0,\infty) \quad \text{for each} \quad \xi \in [0,\infty). \quad (98)$$
In this case, we will use notation $F(x,y) \in \mathbb{F}_{-1}(S)$ for $S \subset \mathbb{R}^2_+$.

**Theorem 6.2** (General Fractional Green's Theorem for Simple Region)
Let $S$ be a simple region in $\mathbb{R}^2_{0,+}$ of the $XY$-plane such that $S \in G^1(\mathbb{R}^2_{0,+})$, and the vector field $F[x,y]$ belongs to the set $\mathbb{F}^1_{-1}(\mathbb{R}^2_+)$.

Then, the equation
$$I^{(M)}_{\partial S}[x]F_x(x,y) + I^{(M)}_{\partial S}[y]F_y(x,y) =$$

$$I^{(M)}_S[x,y]\left(D^{x,*}_{(K_1)}[x']F_y(x',y) - D^{y,*}_{(K_2)}[y']F_x(x,y')\right) \quad (99)$$

are satisfied.

*Proof.*
Let us consider the piecewise simple closed line $\partial S \subset \mathbb{R}^2_{0,+}$ in $XY$-plane, which consists of the simple lines $AB$, $BC$, $CD$, and $DA$,
$$\partial S = AB \cup BC \cup CD \cup DA$$
with the coordinates of the points
$$A(a,y_1(a)), \quad B(x_1(c),c), \quad C(b,y_1(b)), \quad D(x_1(d),d),$$
where
$$y_1(a) = y_2(a), \quad y_1(b) = y_2(b), \quad x_1(c) = x_2(c), \quad x_1(d) = x_2(d).$$
The lines $ABC$ and $ADC$ are $Y$-simple lines
$$ABC := \{(x,y): 0 \le a \le x \le b, y = y_1(x) \ge 0\}, \quad (100)$$

$$ADC := \{(x,y): 0 \le a \le x \le b, y = y_2(x) \ge 0\}, \quad (101)$$
where the functions $y_1(x)$ and $y_2(x)$ belong to the space $C^1(a,b)$, and derivatives of these functions are nonzero and do not change sign on the interval $[a,b]$.

The lines $BAD$ and $BCD$ are $X$-simple lines
$$BAD := \{(x,y): 0 \le c \le y \le d, x = x_1(y) \ge 0\}, \quad (102)$$

$$BCD := \{(x,y): 0 \le c \le y \le c, x = x_2(y) \ge 0\}, \quad (103)$$
where the functions $x_1(x)$ and $x_2(x)$ belong to the space $C^1(c,d)$, and derivatives of these functions do not change sign on the interval $[c,d]$.

Using that the functions $F_x(x,y)$ and $F_y(x,y)$ satisfy the condition
$$F_x(x, y_1(x)), F_x(x, y_2(x)) \in C_{-1}(0, \infty), \tag{104}$$

$$F_y(x_1(y), y), F_y(x(y), y) \in C_{-1}(0, \infty), \tag{105}$$

the general fractional circulation along the piecewise simple closed line $\partial S$ is described by the equation

$$(I_L^{(M)}[x,y] F(x,y)) = \tag{106}$$

$$I_{ABC}^{(M_1)}[x] F_x(x, y_1(x)) + I_{CDA}^{(M_1)}[x] F_x(x, y_2(x)) +$$

$$I_{BCD}^{(M_2)}[y] F_y(x_2(y), y) + I_{DAB}^{(M_2)}[y] F_y(x_1(y), y) =$$

$$I_{ABC}^{(M_1)}[x] F_x(x, y_1(x)) - I_{ADC}^{(M_1)}[x] F_x(x, y_2(x)) +$$

$$I_{BCD}^{(M_2)}[y] F_y(x_2(y), y) - I_{BAD}^{(M_2)}[y] F_y(x_1(y), y) =$$

$$I_{[a,b]}^{(M_1)}[x] F_x(x, y_1(x)) - I_{[a,b]}^{(M_1)}[x] F_x(x, y_2(x)) +$$

$$I_{[c,d]}^{(M_2)}[y] F_y(x_2(y), y) - I_{[c,d]}^{(M_2)}[y] F_y(x_1(y), y) =$$

$$I_{[a,b]}^{(M_1)}[x] (F_x(x, y_1(x)) - F_x(x, y_2(x))) +$$

$$I_{[c,d]}^{(M_2)}[y] (F_y(x_2(y), y) - F_y(x_1(y), y)) =$$

$$- I_{[a,b]}^{(M_1)}[x] (I_{[y_1(x),y_2(x)]}^{(M_2)}[y] D_{(K_2)}^{y,*}[y'] F_x(x, y')) -$$

$$I_{[c,d]}^{(M_2)}[y] (I_{[x_1(y),x_2(x)]}^{(M_1)}[x] D_{(K_1)}^{x,*}[x'] F_y(x', y)) =$$

$$I_{S]}^{(M)}[x,y] (D_{(K_1)}^{x,*}[x'] F_y(x', y) - D_{(K_2)}^{y,*}[y'] F_x(x, y')),$$

where we used the second fundamental theorem of general fractional calculus.

## 6.3 General Fractional Green Theorem for Rectangle

Let us formulate the general fractional Green's theorem for the region in the form of a rectangle.

**Theorem 6.3** (*General Fractional Green's Theorem for Rectangle*)
Let functions $F_x = F_x(x,y)$ and $F_y = F_y(x,y)$ belong to the function space $\mathbb{S}_{-1}^1(\mathbb{R}_+^2)$, and $S$ is the rectangle
$$S := \{(x,y): \ x \in [a,b], \ y \in [c,d]\}, \tag{107}$$

and the boundary of $S$ be the closed line $\partial S$. Then
$$I^{(M)}_{\partial S}[x] F_x(x,y) + I^{(M)}_{\partial S}[y] F_y(x,y) =$$

$$I^{(M)}_S[x,y] \, (D^{y,*}_{(K_2)}[y'] \, F_x(x,y') - D^{x,*}_{(K_1)}[x'] \, F_y(x',y)). \tag{108}$$

*Proof.* 1) To prove equation (108), the double fractional integral $I^{(M)}_S[x,y]$ is written as the iterated (repeated) fractional integrals and then the fundamental theorems of general fractional calculus is used. Let $S$ be the rectangular domain (107) with vertexes in the points
$$A(a,c), \quad B(b,c), \quad C(b,d), \quad D(a,d).$$
The sides $AB$, $BC$, $CD$, $DA$ of the rectangular domain (107) form the boundary $\partial S$ of $S$. For the rectangular region $S$ defined by $a \leq x \leq b$, and $c \leq y \leq d$, the iterated (repeated) integral is
$$I^{(M)}_S[x,y] = I^{(M_1)}_{[a,b]}[x] \, I^{(M_2)}_{[c,d]}[y], \tag{109}$$
where
$$I^{(M_1)}_{[a,b]}[x] = I^b_{(M_1)}[x] - I^a_{(M_1)}[x], \quad I^{(M_2)}_{[c,d]}[y] = I^d_{(M_2)}[y] - I^c_{(M_2)}[y].$$

2) To prove the general fractional Green's formula, we realize the following transformations
$$I^{(M)}_{\partial S}[x] F_x + I^{(M)}_{\partial S}[y] F_y =$$

$$I^{(M_1)}_{AB}[x] F_x + I^{(M_1)}_{CD}[x] F_x + I^{(M_2)}_{BC}[y] F_y + I^{(M_2)}_{DA}[y] F_y =$$

$$I^{(M_1)}_{[a,b]}[x] F_x(x,c) - I^{(M_1)}_{[a,b]}[x] F_x(x,d) + I^{(M_2)}_{[c,d]}[y] F_y(b,y) - I^{(M_2)}_{[c,d]}[y] F_y(a,y) =$$

$$I^{(M_1)}_{[a,b]}[x] (F_x(x,c) - F_x(x,d)) + I^{(M_2)}_{[c,d]}[y] (F_y(b,y)dy - F_y(a,y)). \tag{110}$$

Let us use the second fundamental theorem of GFC in the from
$$F_x(x,c) - F_x(x,d) = -I^{(M_2)}_{[c,d]}[y] \, D^{y,*}_{(K_2)}[y'] \, F_x(x,y'), \tag{111}$$

$$F_y(b,y) - F_y(a,y) = I^{(M_1)}_{[a,b]}[x] \, D^{x,*}_{(K_1)}[x'] \, F_y(x',y). \tag{112}$$

Then, expression (110) can be written as
$$I^{(M_1)}_{[a,b]}[x] \, (I^{(M_2)}_{[c,d]}[y] \, D^{y,*}_{(K_2)}[y'] \, F_x(x,y')) + I^{(M_2)}_{[c,d]}[y] \, (-I^{(M_1)}_{[a,b]}[x] \, D^{x,*}_{(K_1)}[x'] \, F_y(x',y)) =$$

$$I^{(M_1)}_{[a,b]}[x] \, I^{(M_2)}_{[c,d]}[y] \, (D^{y,*}_{(K_2)}[y'] \, F_x(x,y') - D^{x,*}_{(K_1)}[x'] \, F_y(x',y)) =$$

$$I^{(M)}_S[x,y] \, (D^{y,*}_{(K_2)}[y'] \, F_x(x,y') - D^{x,*}_{(K_1)}[x'] \, F_y(x',y)), \tag{113}$$

where we use (109) for the operator $I^{(M)}_S[x,y]$. As a result, we obtain the right-hand side of equation (108).

**Remark 6.3** *For the kernels*
$$M_1(x) = h_{\alpha_1}(x), \quad M_2(y) = h_{\alpha_2}(y), \quad K_1(x) = h_{1-\alpha_1}(x), \quad K_2(y) = h_{1-\alpha_2}(y),$$
*a general fractional Green's theorem was proved in [30], where the Riemann-Liuoville fractional integrals and the Caputo fractional derivatives are used. For $\alpha_1 = \alpha_2 = 1$, we get the standard Green's theorem.*

# 7 Double and Surface General Fractional Integral, and Flux

## 7.1 Definition of Double GFI by Iterated (Repeated) GFI

Let us consider a region $S$ in the first quarter of the coordinate plane $\mathbb{R}_+^2$, and assume that the area $S$ is bounded such that

$$S := \{(x,y): \ 0 \leq a \leq x \leq b, \ 0 \leq y_1(x) \leq y \leq y_2(x)\}, \tag{114}$$

where $y = y_1(x)$ and $y = y_2(x)$ are continuous function, and following two requirements are satisfied: (1) the projection of $S$ onto the $X$-axis is bounded by the two values $x = a$ and $x = b$; (2) any line parallel to the Y-axis that passes between these two values $x = a$ and $x = b$ intersects the line $y = y_1(x)$ (and $y = y_2(x)$ in the interval $(a,b)$ at no more than one point. Then the region $S$ will be called the $Y$-simple domain.

Let us define the double GFI through the iterated (repeated) GFI.

**Definition 7.1** If the function $f = f(x,y) \in C_{-1}(\mathbb{R}_{+1}^2)$ satisfies the condition

$$J(x) := I_{[y_1(x),y_2(x)]}^{(M_2)}[y]f(x,y) \in C_{-1}(0,\infty),$$

where

$$I_{[y_1(x),y_2(x)]}^{(M_2)}[y]f(x,y) := I_{(M_2)}^{y_2(x)}[y]f(x,y) - I_{(M_2)}^{y_1(x)}[y]f(x,y),$$

and $S$ is the $Y$-simple domain (114), then the iterated (repeated) GFI of the form

$$I_S^{(M)}[x,y] := I_{[a,b]}^{(M_1)}[x] \, I_{[y_1(x),y_2(x)]}^{(M_2)}[y], \tag{115}$$

is called the double GFI, where

$$I_{[a,b]}^{(M_1)}[x] = I_{(M_1)}^b[x] - I_{(M_1)}^a[x].$$

Using the function

$$F(x,y) = I_{(M_2)}^y[y']f(x,y') = \int_0^y M_2(y-y') \, f(x,y') \, dy',$$

equation (115) can be written in the form

$$I_S^{(M)}[x,y] \, f(x,y) := \int_0^b M_1(b-x) \, (F(x,y_2(x)) - F(x,y_1(x))) \, dx -$$

$$\int_0^a M_1(a-x) \, (F(x,y_2(x)) - F(x,y_1(x))) \, dx. \tag{116}$$

**Example 7.1** For the rectangular area

$$S := \{(x,y): \ 0 \leq a \leq x \leq b, \ 0 \leq c \leq y \leq d\}, \tag{117}$$

double GFI (115) has the form

$$I_S^{(M)}[x,y] := I_{[a,b]}^{(M_1)}[x] \, I_{[c,d]}^{(M_2)}[y]. \tag{118}$$

If the function $f = f(x,y) \in C_{-1}(\mathbb{R}_{+1}^2)$ and $a = c = 0$, then

$$I_S^{(M)}[x,y] \, f(x,y) = I_{(M_1)}^b[x] \, I_{(M_2)}^d[y] \, f(x,y) =$$

$$\int_0^b \int_0^d M_1(b-x) \, M_2(d-y) \, f(x,y) \, dx \, dy. \tag{119}$$

**Example 7.2** For the area
$$S := \{(x,y): \ 0 \le x \le b, \ 0 \le y \le g(x)\}, \tag{120}$$
and the function $f = f(x,y) \in C_{-1}(\mathbb{R}^2_{+1})$, the double GFI has the form
$$I_S^{(M)}[x,y]\, f(x,y) := I_{(M_1)}^b[x]\, I_{(M_2)}^{g(x)}[y]\, f(x,y) =$$

$$\int_0^b \int_0^{g(x)} M_1(b-x)\, M_2(g(x)-y)\, f(x,y)\, dx\, dy. \tag{121}$$

**Remark 7.1** Let us conseder double GFI (115) for the function
$$f(x,y) = D^{y,*}_{(K_2)}[y']F_x(x,y') = \int_0^y dy'\, K_2(y-y')\, F_x(x,y'),$$
where $D^{y,*}_{(K_2)}[y']$ is the GFD with respect to the variable $y$ with kernel $K_2(y)$ that is accociated to the kernel $M_2(y)$. Then
$$I_S^{(M)}[x,y]\, f(x,y) = I_{(M_1)}^a[x]\, I_{(M_2)}^{g(x)}[y]\, D^{y,*}_{(K_2)}[y']\, F_x(x,y').$$
For area (120), we obtain
$$I_S^{(M)}[x,y]\, D^{y,*}_{(K_2)}[y']\, F_x(x,y') = I_{(M_1)}^a[x]\, I_{(M_2)}^{g(x)}[y]\, D^{y,*}_{(K_2)}[y']\, F_x(x,y') =$$

$$I_{(M_1)}^a[x]\, (F_1(x,g(x)) - F_x(x,0)),$$
where the second fundamental theorem of GFC is used. In the general case, we obtain
$$I_S^{(M)}[x,y]\, D^{y)}_{(K_2)}[y]\, F_x(x,y) = I_{[a,b]}^{(M_1)}[x]\, I_{[y_1(x),y_2(x)]}^{(M_2)}[y] D^{y,*}_{(K_2)}[y']\, F_x(x,y') =$$

$$I_{[a,b]}^{(M_1)}[x]\, (F_x(x,y_2(x)) - F_x(x,y_1(x))),$$

$$\int_0^b dx\, M_1(b-x)\, (F_x(x,y_2(x)) - F_x(x,y_1(x))) -$$

$$\int_0^a dx\, M_1(a-x)\, (F_x(x,y_2(x)) - F_x(x,y_1(x))),$$
where $y_2(x) \ge y_1(x) > 0$ for all $0 \le x \le b$ (instead of $a \le x \le b$). Here we assume that
$$F_x(\xi,y) \in C^1_{-1}(0,\infty) \quad \text{forall} \quad \xi \in (0,\infty),$$

$$F_x(x,y_1(x)),\ F_x(x,y_2(x)) \in C_{-1}(0,\infty).$$

## 7.2 Surface General Fractional Integral

The standard surface integral of integer order is a generalization of double integrals to integration over surfaces. This integral can be considered as the double integral analogue of the line integral.

The surface general fractional integral over the surface $S$ in $\mathbb{R}^3_{0,+}$ can be defined through the double general fractional integrals over areas $S_{yz}$, $S_{xz}$, $S_{xy}$ in the $YZ$, $XZ$, $XY$ planes, where these areas are projections of the surface $S$ onto these planes.

Let us consider a two-sided smooth surface $S$, and fix one of its two sides, which is equivalent to choosing a certain orientation on the surface. We also assume that the surface is given by the equation
$$x = x(y,z) \ge 0,$$

where the point $(y,z)$ changes in area $S_{yz}$ in the $YZ$-plane, bounded by smooth contour $\partial S_{yz}$, and
$$S := \{(x,y,z): \quad x = x(y,z) \geq 0, \quad (y,z) \in S_{yz} \subset \mathbb{R}^2_{0,+}\}.$$
Let us first consider a vector field of the form
$$F_x := e_x\, F_x(x,y,z). \tag{122}$$
If the function $F_x(x,y,z)$ on the surface $S$ belongs to the function space $C_{-1}(\mathbb{R}^2_+)$ that is
$$F_x(x(y,z),y,z) \in C_{-1}(\mathbb{R}^2_+), \tag{123}$$
then the surface general fractional integral is defined by the equation
$$(I_S^{(M)}, F_x) = I_{S_{yz}}^{(M)}[y,z]\, F_x(x(y,z),y,z). \tag{124}$$

**Remark 7.2** If the area $S_{yz}$ in the $YZ$-plane is given as
$$S_{yz} := \{(y,z): \quad 0 \leq c \leq y \leq d, \quad 0 \leq z_1(y) \leq z \leq z_2(y)\},$$
where $z = z_1(y)$ and $z = z_2(y)$ are continuous functions, then the double GFI in equation (124) is defined as
$$I_{S_{yz}}^{(M)}[y,z]\, F_x(x(y,z),y,z) := I_{[c,d]}^{(M_2)}[z]\, I_{[z_1(y),z_2(y)]}^{(M_3)}[y]\, F_x(x(y,z),y,z).$$

Similar to vector field (122), we can obtain the surface GFIs for the other two projections of the surface $S$ on the $XZ$ and $YZ$ planes, and functions $F_y(x,y,z)$, $F_z(x,y,z)$ on the surface $S$. If, instead of the $YZ$-plane, we project the surface $S$ onto the $XZ$-plane or $YZ$-plane, then we get two other surface GFI of the second type:
$$(I_S^{(M)}, F_y) = I_{S_{xz}}^{(M)}[x,z]\, F_y(x,y(x,z),z),$$
$$(I_S^{(M)}, F_z) = I_{S_{xy}}^{(M)}[x,y]\, F_z(x,y,z(x,y)),$$
where
$$F_y(x,y(x,z),z) \in C_{-1}(\mathbb{R}^2_+), \quad F_z(x,y,z(x,y)) \in C_{-1}(\mathbb{R}^2_+). \tag{125}$$
In this case, the surface GFI is defined by projections onto all three planes ($YZ$, $XZ$, $XY$) in the form
$$(I_S^{(M)}, F) = I_{S_{yz}}^{(M)}[y,z]\, F_x(x(y,z),y,z) +$$
$$I_{S_{xz}}^{(M)}[x,z]\, F_y(x,y(x,z),z) + I_{S_{xy}}^{(M)}[x,y]\, F_z(x,y,z(x,y)),$$
where conditions (123) and (125) are satisfied.

Let us define a set of piecewise simple surfaces, a set of vector fields on these surfaces and the surface GFI.

**Definition 7.2** *Let $S$ be an oriented compact smooth surface in the region*
$$\mathbb{R}^3_{0,+} = \{(x,y,z): \quad x \geq 0, \quad y \geq 0, \quad z \geq 0\}.$$
*Let us choose a side of the surface $S$.*

*Let the surface $S$ be represented as the union of a finite number of $X$-simple surfaces $S_{X,i}$, as well as a finite number of $Y$-simple surfaces $S_{X,j}$, and $Z$-simple surfaces $S_{Z,k}$, where $i = 1,\ldots,n_x$, $j = 1,\ldots,n_y$, $k = 1,\ldots,n_z$:*
$$S := \cup_{k=1}^{n_x} S_{X,i} = \cup_{j=1}^{n_y} S_{Y,j} = \cup_{k=1}^{n_z} S_{Z,k}.$$
*Let $S_{i,yz}$, $S_{j,xz}$, $S_{k,xy}$ be projections of the surfaces $S_{X,i}$, $S_{Y,j}$, $S_{Z,k}$ onto the $YZ$, $XZ$, $XY$*

*planes, which can be described by continuous functions*

$$S_{X,i} = \{x = x_i(y,z) \geq 0 \quad if \quad (y,z) \in S_{i,yz} \subset \mathbb{R}^2_{0,+}\},$$
$$S_{Y,j} = \{y = y_j(x,z) \geq 0 \quad if \quad (x,z) \in S_{j,xz} \subset \mathbb{R}^2_{0,+}\},$$
$$S_{Z,k} = \{z = z_k(x,y) \geq 0 \quad if \quad (x,y) \in S_{k,xy} \subset \mathbb{R}^2_{0,+}\}.$$

*Such surfaces $S$ will be called piecewise simple surfaces. The set of such surfaces will be denoted by $\mathbb{P}(\mathbb{R}^3_{0,+})$.*

**Definition 7.3** *Let $S$ be a piecewise simple surface ($S \in \mathbb{P}(\mathbb{R}^3_{0,+})$). Let the vector field*

$$F := e_x F_x(x,y,z) + e_y F_y(x,y,z) + e_z F_z(x,y,z) \tag{126}$$

*on the surface $S$ satisfy the conditions*

$$F_x(x_i(y,z), y, z) \in C_{-1}(\mathbb{R}^2_+), \tag{127}$$
$$F_y(x, y_j(x,z), z) \in C_{-1}(\mathbb{R}^2_+), \tag{128}$$
$$F_z(x, y, z_k(x,y)) \in C_{-1}(\mathbb{R}^2_+) \tag{129}$$

*for all $i = 1, \ldots, n_x$, $j = 1, \ldots, n_y$, $k = 1, \ldots, n_z$.*
*The set of such vector fields $F$ on piecewise simple surface $S$ will be denoted by $\mathbb{F}_S(\mathbb{R}^3_{0,+})$.*

**Definition 7.4** *Let $S$ be a piecewise simple surface ($S \in \mathbb{P}(\mathbb{R}^3_{0,+})$) and a vector fields $F$ on this surface $S$ belongs to the set $\mathbb{F}_S(\mathbb{R}^3_{0,+})$.*
*Then the surface general fractional vector integral (surface GFI) of the second kind*

$$I_S^{(M)} = e_x I_{S_{yz}}^{(M)}[y,z] + e_y I_{S_{xz}}^{(M)}[z,x] + e_z I_{S_{xy}}^{(M)}[x,y] =$$

$$e_x \sum_{i=1}^{n_x} I_{S_{i,yz}}^{(M)}[y,z] + e_y \sum_{j=1}^{n_y} I_{S_{j,xz}}^{(M)}[z,x] + e_z \sum_{k=1}^{n_z} I_{S_{k,xy}}^{(M)}[x,y] \tag{130}$$

*for the vector fiels $F \in \mathbb{F}_S(\mathbb{R}^3_{0,+})$ is defined by the equation*

$$(I_S^{(M)}, F) := \sum_{i=1}^{n_x} I_{S_{i,yz}}^{(M)}[y,z] F_x(x_i(y,z), y, z) +$$

$$\sum_{j=1}^{n_y} I_{S_{j,xz}}^{(M)}[x,z] F_y(x, y_j(x,z), z) + \sum_{k=1}^{n_z} I_{S_{k,xy}}^{(M)}[x,y] F_z(x, y, z_k(x,y)).$$

### 7.3  General Fractional Flux

Let us consider the general fractional flux for simple surfaces.

**Definition 7.5** *Let $S$ be an oriented compact smooth simple surface in the region $\mathbb{R}^3_{0,+}$ such that $S_{yz}$, $S_{xz}$, $S_{xy}$ are the projections of the surface $S$ onto the YZ, XZ, XY planes, which can be described by continuous functions*

$$S = \{x = x(y,z) \geq 0 \quad if \quad (y,z) \in S_{yz} \subset \mathbb{R}^2_{0,+}\},$$

$$S = \{y = y(x,z) \geq 0 \quad if \quad (x,z) \in S_{xz} \subset \mathbb{R}^2_{0,+}\},$$

$$S = \{z = z(x,y) \geq 0 \quad if \quad (x,y) \in S_{xy} \subset \mathbb{R}^2_{0,+}\}.$$

*Then a general fractional flux of the vector field $F \in \mathbb{F}_S(\mathbb{R}^3_{0,+})$ across the simple surface $S$ is the surface general fractional integral of the form*

$$\Phi_S^{(M)}(F) = (I_S^{(M)}, F) = I_{S_{yz}}^{(M)}[y,z]F_x(x(y,z),y,z) +$$

$$I_{S_{xz}}^{(M)}[x,z]F_y(x,y(x,z),z) + I_{S_{xy}}^{(M)}[x,y]F_z(x,y,z(x,y)), \quad (131)$$

where

$$F_x(x(y,z),y,z), F_y(x,y(x,z),z), F_z(x,y,z(x,y)) \in C_{-1}(\mathbb{R}_+^2).$$

**Example 7.3** Let us consider a surface $S$ in the form of the rectangle
$$S = S_{yz} := \{(y,z): \ x = 0, \ 0 \le c \le y \le d, \ 0 \le e \le z \le f\}.$$
The general fractional flux of the vector field $F$ across this surface is written as
$$\Phi_S^{(M)}(F) = I_{S_{yz}}^{(M)}[y,z]\,F_x(0,y,z) = I_{[c,d]}^{(M_2)}[y]\,I_{[e,f]}^{(M_3)}[z]\,F_x(0,y,z).$$

**Example 7.4** Let us consider a surface $S$ in the form of boundary $\partial W$ of the parallelepiped area
$$W := \{\{(x,y,z): \ 0 \le a \le x \le b, \ 0 \le c \le y \le d, \ 0 \le e \le z \le f\}.$$
The general fractional flux of the vector field
$$F = e_x F_x(x,y,z) + e_y F_y(x,y,z) + e_z F_z(x,y,z)$$
across the surface $S = \partial W$ is written as
$$\Phi_S^{(M)}(F) = (I_S^{(M)}, F) =$$

$$I_{S_{yz}}^{(M)}[y,z]\,F_x(x,y,z) + I_{S_{xz}}^{(M)}[x,z]\,F_y(x,y,z) + I_{S_{xy}}^{(M)}[x,y]\,F_z(x,y,z) =$$

$$I_{[c,d]}^{(M_2)}[y]\,I_{[e,f]}^{(M_3)}[z]\,F_x(b,y,z) - I_{[c,d]}^{(M_2)}[y]\,I_{[e,f]}^{(M_3)}[z]\,F_x(a,y,z) +$$

$$I_{[a,b]}^{(M_1)}[x]\,I_{[e,f]}^{(M_3)}[z]\,F_y(x,c,z) - I_{[a,b]}^{(M_1)}[x]\,I_{[e,f]}^{(M_3)}[z]\,F_y(x,d,z) +$$

$$I_{[a,b]}^{(M_1)}[x]\,I_{[c,d]}^{(M_2)}[y]\,F_z(x,y,f) - I_{[a,b]}^{(M_1)}[x]\,I_{[c,d]}^{(M_2)}[y]\,F_z(x,y,e),$$

where

$$I_{[a_k,b_k]}^{(M_k)}[x_k] := I_{(M_k)}^{b_k}[x_k] - I_{(M_k)}^{a_k}[x_k].$$

**Example 7.5** For the rectangular area
$$S := \{0 \le x \le a, \ 0 \le y \le b, \ z = 0\}, \quad (132)$$
the general fractional flux of the vector field $F$ across surface (132) is written as
$$\Phi_S^{(M)}(F) = I_{S_{xy}}^{(M)}[x,y]\,F_z(x,y,0) = I_{S_{xy}}^{(M_1)}[x]\,I_{S_{xy}}^{(M_2)}[y]\,F_z(x,y,0) =$$
$$I_{(M_1)}^a[x]\,I_{(M_2)}^b[y]\,F_z(x,y,0) = \int_0^a dx \int_0^b dy\, M_1(a-x)\,M_2(b-y)\,F_z(x,y,0). \quad (133)$$
For kernels (62), the general fractional flux (133) for surface (132) has the form
$$\Phi_S^{(M)}(F) = \int_0^a dx \int_0^b dy\, \frac{(a-x)^{\alpha_1-1}}{\Gamma(\alpha_1)}\,\frac{(b-y)^{\alpha_2-1}}{\Gamma(\alpha_2)}\,F_z(x,y,0). \quad (134)$$
For $\alpha_1 = \alpha_2 = 1$, the general fractional flux (134) for surface (132) is the standard flux
$$\Phi_S(F) = \int_0^a dx \int_0^b dy\, F_z(x,y,0). \quad (135)$$

# 8 General Fractional Stokes Theorem

The Stokes theorem connects the surface integral with the line integral. The Stokes theorem generalizes Green theorem from $\mathbb{R}^2$ to $\mathbb{R}^3$. If the surface is a flat region lying in the plane, then the Stokes equation gives the Green equation.

## 8.1 Simple Domain on $XY$-Plane

We recall the definitions of simple domains.

**Definition 8.1** *The closed domain $D_{xy} \subset \mathbb{R}^2_{0,+}$ on the $XY$-plane will be called $Y$-simple domain, if $D_{xy}$ can be represented in the form*
$$D_{xy} := \{(x,y): \ 0 \leq a \leq x \leq b, \ 0 \leq y_1(x) \leq y \leq y_2(x)\},$$
*where $y = y_1(x)$, $y = y_2(x)$ are continuous functions for $x \in [a,b]$.*

Similarly, the closed domain $D_{xy}$ on the $XY$-plane is called $X$-simple domain, if $D_{xy}$ can be represented as
$$D_{xy} := \{(x,y): \ 0 \leq x_1(y) \leq x \leq x_2(y), \ 0 \leq c \leq y \leq d\},$$
where $x = x_1(y)$, $x = x_2(y)$ are continuous functions for $y \in [c,d]$.

**Definition 8.2** *If $D_{xy}$ is $X$-simple and $Y$-simple domain, then $D_{xy}$ is called the simple domain on the $XY$-plane.*
*The simple domain $D_{xy} \subset \mathbb{R}^2_{0,+}$ on the $XY$-plane can be represented in the form*
$$D_{xy} := \{(x,y): \ 0 \leq a \leq x \leq b, \ 0 \leq y_1(x) \leq y \leq y_2(x)\},$$

$$D_{xy} := \{(x,y): \ 0 \leq x_1(y) \leq x \leq x_2(y), \ 0 \leq c \leq y \leq d\},$$
where $y_1(x), y_2(x) \in C[a,b]$ and $x_1(y), x_2(y) \in C[c,d]$, and $y_1(a) = y_2(a) = c$, $y_1(b) = y_2(b) = d$, $x_1(c) = x_2(c) = a$, $x_1(d) = x_2(d) = b$.

The boundary of a $Y$-simple closed domain $D_{xy}$ in $\mathbb{R}^2_{0,+}$ can be represented as a closed line $L_{xy} = \partial D_{xy}$ consisting of the lines
$$L_{y1} := \{(x,y): \ 0 \leq a \leq x \leq b, \ y = y_1(x) \geq 0\},$$

$$L_{y2} := \{(x,y): \ 0 \leq a \leq x \leq b, \ y = y_2(x) \geq 0\},$$

$$L_{x1} := \{(x,y): \ x = a \geq 0 \ \ 0 \leq y_1(a) \leq y \leq y_2(a)\},$$

$$L_{x2} := \{(x,y): \ x = b \geq 0 \ \ 0 \leq y_1(b) \leq y \leq y_2(b)\}.$$
To simplify the proofs, we can sometimes use the case, when there are no straight lines $B_a$, $B_b$ that is, when $y_1(a) = y_2(a)$ and $y_1(b) = y_2(b)$.
Similarly, we obtain a representation of the boundary for a $X$-simple domian $D_{xy}$.

## 8.2 Simple Surface

Let us give a definition of $Z$-simple surface.

**Definition 8.3** Let $S$ be a smooth oriented surface in $\mathbb{R}^3_{0,+}$ that the surface $S$ is described by the equation
$$z = z(x,y) \geq 0,$$
where the function $z(x,y)$ is continuous in the closed domain $D_{xy}$ ($z(x,y) \in C(D_{xy})$), which is a projection of the surface $S$ onto the $XY$-plane.

We will assume that is bounded by an oriented closed smooth line $L = \partial S$. The boundary $\partial D_{xy}$ of the domain $D_{xy}$ is a projection of the line $L = \partial S$ onto the $XY$-plane.

Then the surface $S$ will be called the Z-simple surface, if the domain $D_{xy}$ is simple on the $XY$-plane.

**Definition 8.4** The surface $S$ will be called the simple surface in $\mathbb{R}^3_{0,+}$, if $S$ is simple with respect to the $X$, $Y$ and $Z$ axes.

The simple surface $S$ can be described by the continuous functions
$$S := \{(x,y,z): \ x = x(y,z) \in C(\mathbb{R}^2_+), \ (y,z) \in D_{yz} \subset \mathbb{R}^2_{0,+}\}, \qquad (136)$$

$$S := \{(x,y,z): \ y = y(x,z) \in C(\mathbb{R}^2_+), \ (x,z) \in D_{xz} \subset \mathbb{R}^2_{0,+}\}, \qquad (137)$$

$$S := \{(x,y,z): \ z = z(x,y) \in C(\mathbb{R}^2_+), \ (x,y) \in D_{x,y} \subset \mathbb{R}^2_{0,+}\}, \qquad (138)$$
where the domains $D_{yz}$, $D_{xz}$, $D_{xy}$ are simple domains in $\mathbb{R}^2_{0,+}$, which are projections of the surface $S$ onto the $YZ$, $XZ$, $XY$ planes.

**Remark 8.1** We can consider the Z-simple surfaces $S$ in $\mathbb{R}^3_{0,+}$ that can be described as a union of X-simple surfaces, as well as a union of Y-simple regions.

### 8.3   Vector Field on Surface

Let us define the properties of the scalar and vector fields in simple domains in $\mathbb{R}^2_{0,+}$.

**Definition 8.5** Let $D_{xy} \subset \mathbb{R}^2_{0,+}$ be a simple closed domain on the $XY$-plane that can be described in the forms
$$D_{xy} := \{(x,y): \ 0 \leq a \leq x \leq b, \ 0 \leq y_1(x) \leq y \leq y_2(x)\},$$

$$D_{xy} := \{(x,y): \ 0 \leq x_1(y) \leq x \leq x_2(y), \ 0 \leq c \leq y \leq d\},$$
where $y_1(x), y_2(x) \in C[a,b]$ and $x_1(y), x_2(y) \in C[c,d]$ .

Let $f(x,y)$ be a scalar field, which is defined in the simple domain $D_{xy} \subset \mathbb{R}^2_{0,+}$ on the $XY$-plane such that the following conditions are satisfied
$$f(x,y_1(x)), f(x,y_2(x)) \in C^n_{-1}(0,\infty),$$

$$f(x_1(y),x), f(x_2(y),y) \in C^n_{-1}(0,\infty).$$
Then this property of the field $f(x,y)$ will be denoted as $f(x,y) \in \mathbb{F}^n_{-1}(D_{xy})$.

Let us consider the vector field
$$F(x,y) = e_x F_x(x,y) + e_y F_y(x,y)$$

which is defined in the simple domain $D_{xy} \subset \mathbb{R}^2_{0,+}$. The fact that the field $F(x,y)$ satisfies the conditions
$$F_x(x, y_1(x)), \quad F_x(x, y_2(x)) \in C^n_{-1}(0, \infty),$$

$$F_y(x_1(y), y), \quad F_y(x(y), y) \in C^n_{-1}(0, \infty),$$
will be denoted as $F(x,y) \in \mathbb{F}^n_{-1}(D_{xy})$.

Let us define the properties of the scalar and vector fields on simple surfaces in $\mathbb{R}^3_{0,+}$.

**Definition 8.6** Let $S$ be a simple surface in $\mathbb{R}^3_{0,+}$ that is described by equations (136), (137), (138), and $D_{yz}, D_{xz}, D_{xy}$ are simple domains, which are projections of the surface $S$ onto the $YZ, XZ, XY$ planes.

Let us consider the field $f(x, y, z)$ that is defined on the simple surface $S \subset \mathbb{R}^3_{0,+}$. The fact that the field $f(x, y, z)$ satisfies the conditions
$$f(x, y, z(x, y)) \in \mathbb{F}^n_{-1}(D_{xy}),$$

$$f(x(y, z), y, z) \in \mathbb{F}^n_{-1}(D_{yz}),$$

$$f(x, y(x, z), z) \in \mathbb{F}^n_{-1}(D_{xz}),$$
will be denoted as $f(x, y, z) \in \mathbb{F}^n_{-1}(S)$.

**Definition 8.7** Let $S$ be a simple surface in $\mathbb{R}^3_{0,+}$ that is described by equations (136), (137), (138), and $D_{yz}, D_{xz}, D_{xy}$ are projections of the surface $S$ onto the $YZ, XZ, XY$ planes.

Let us consider the vector field
$$F(x, y, z) = e_x F_x(x, y, z) + e_y F_y(x, y, z) + e_z F_z(x, y, z)$$
that is defined on the simple surface $S \subset \mathbb{R}^3_{0,+}$.

The fact that the field $F(x, y, z)$ satisfies the conditions
$$F_x(x, y, z(x, y)), \quad F_y(x, y, z(x, y)) \in \mathbb{F}^n_{-1}(D_{xy}),$$

$$F_x(x, y(x, z), z), \quad F_z(x, y(x, z), z) \in \mathbb{F}^n_{-1}(D_{xz}),$$

$$F_y(x(y, z), y, z), \quad F_z(x(y, z), y, z) \in \mathbb{F}^n_{-1}(D_{yz}),$$
will be denoted as $F(x, y, z) \in \mathbb{F}^n_{-1}(S)$.

## 8.4   General Fractional Vector Integrals over Surface

Let us consider a smooth oriented surface $S$ in $\mathbb{R}^3_{0,+}$, which is bounded by an oriented closed smooth line $L = \partial S$. We will assume that $S$ is the simple surface in $\mathbb{R}^3_{0,+}$, and $D_{yz}, D_{xz}, D_{xy}$ are simple domains, which are projections of the surface $S$ onto the $YZ, XZ, XY$ planes.

**Definition 8.8** The surface GFI opereator over the simple surface $S \in \mathbb{R}^3_{0,+}$ is defined by the equation
$$I^{(M)}_S := I^{(M)}_{D_{xy}} + I^{(M)}_{D_{xz}} + I^{(M)}_{D_{yz}} =$$
$$n_x I^{(M)}_{D_{yz}} + n_y I^{(M)}_{D_{xz}} + n_z I^{(M)}_{D_{xy}},$$
where the vectors $n_x$, $n_y$, $n_z$ are the normal vector to surfaces $D_{yz}$, $D_{xz}$, $D_{xy}$ that are related with orientation of the closed lines $L_{yz} = \partial D_{yz}$, $L_{xz} = \partial D_{xz}$, $L_{xy} = \partial D_{xy}$ bounding its.

We can use the vectors $e_x$, $e_y$, $e_z$ instead of $n_x$, $n_y$, $n_z$, if all orientation of the closed lines $L_{yz} = \partial D_{yz}$, $L_{xz} = \partial D_{xz}$, $L_{xy} = \partial D_{xy}$ is positive.

For the simple domain $D_{xy}$, the vector operator $I^{(M)}_{D_{xy}}$ is defined by the equtaion
$$I^{(M)}_{D_{xy}} = e_z I^{(M)}_{D_{xy}} = e_z I^{(M_1)}_{[a,b]} [x] I^{(M_2)}_{[y_1(x),y_2(x)]} = - e_z I^{(M_2)}_{[c,d]} [x] I^{(M_1)}_{[x_1(y),x_2(y)]},$$
where $0 \leq y_1(x) \leq y_2(x)$ for $x \in [a,b]$, and $0 \leq x_1(y) \leq x_2(y)$ for $y \in [c,d]$. The minus sign in front of the last expression is due to the relationship between the normal vector of the surface $D_{xy}$ and the orientation of the closed contours $L_{xy} = \partial D_{xy}$.

Similarly, we obtain expressions for the simple domains $D_{xz}$, $D_{yz}$.

**Definition 8.9** The surface GFI of the vector field $F(x, y, z) \in \mathbb{F}^0_{-1}(S)$ over the simple surface $S \in \mathbb{R}^3_{0,+}$ is defined by the equation
$$(I^{(M)}_S, F) = I^{(M)}_{D_{yz}} [y, z] F_x(x(y, z), y, z) +$$

$$I^{(M)}_{D_{xz}} [x, z] F_y(x, y(x, z), z) + I^{(M)}_{D_{xy}} [x, y] F_z(x, y, z(x, y)),$$
if $n_x = e_x$, $n_y = e_y$, $n_z = e_z$.

**Definition 8.10** Let $S \in \mathbb{R}^3_{0,+}$ be a Z-simple surface that surface can be represented as the union of X-simple surfaces $S_{X,k}$ and Y-simple surfaces $S_{Y,k}$ such that
$$S := \bigcup_{k=1}^n S_{X,k} = \bigcup_{l=1}^m S_{Y,l},$$
and $D_{yz,k}$, $D_{xz,l}$ are projections of $S_{X,k}$ and $S_{Y,l}$ on the YZ and XZ planes.

The surface GFI of the vector field $F(x, y, z) \in \mathbb{F}^0_{-1}(S)$ over this surface $S \in \mathbb{R}^3_{0,+}$ is defined by the equation
$$(I^{(M)}_S, F) = \sum_{k=1}^n I^{(M)}_{D_{yz,k}} [y, z] F_x(x_k(y, z), y, z) +$$

$$\sum_{l=1}^l I^{(M)}_{D_{xz,l}} [x, z] F_y(x, y_l(x, z), z) + I^{(M)}_{D_{xy}} [x, y] F_z(x, y, z(x, y)),$$
if $n_x = e_x$, $n_y = e_y$, $n_z = e_z$.

## 8.5 General Fractional Curl Operators

### 8.5.1 Regional General Fractional Curl Operator

Let us give the definition of general fractional curl operator for $\mathbb{R}^3_{0,+}$.

**Definition 8.11** Let $F(x,y,z)$ be a vector field tha belongs to the set $C^1_{-1}(\mathbb{R}^3_+)$.
Then the regional general fractional curl (Regional GF Curl) for the region $W = \mathbb{R}^3_{0,+}$ is defined as
$$Curl^{(K)}_W F = e_x (Curl^{(K)}_W F)_x + e_y (Curl^{(K)}_W F)_y + e_z (Curl^{(K)}_W F)_z,$$
where
$$(Curl^{(K)}_W F)_x(x,y,z) = D^{y,*}_{(K_2)}[y'] F_z(x,y',z) - D^{z,*}_{(K_3)}[z'] F_y(x,y,z'),$$
$$(Curl^{(K)}_W F)_y(x,y,z) = D^{z,*}_{(K_3)}[z'] F_x(x,y,z') - D^{x,*}_{(K_1)}[x'] F_z(x',y,z),$$
$$(Curl^{(K)}_W F)_z(x,y,z) = D^{x,*}_{(K_1)}[x'] F_y(x',y,z) - D^{y,*}_{(K_2)}[y'] F_x(x,y',z)$$
with $x,y,z \geq 0$.

**Remark 8.2** If the vector field
$$F(x,y,z) = \sum_{k=1}^3 e_k F_k(x,y,z)$$
belongs to the function space $C^1_{-1}(\mathbb{R}^3_{0,+})$, then the regional general fractional curl operator for the region $W = \mathbb{R}^3_+$ can be written in the compact form
$$(Curl^{(K)}_W F)(x,y,z) = [D^{(K)}_W, F] = e_n \varepsilon_{nmk} D^{x_n,*}_{(K_n)}[x'_n] F_k,$$
where
$$D^{(K)}_W := e_x D^{x,*}_{(K_1)}[x'] + e_y D^{y,*}_{(K_2)}[y'] + e_z D^{z,*}_{(K_3)}[z'],$$
and $\varepsilon_{lmn}$ is Levi-Civita symbol, i.e. it is $1$, if $(i,j,k)$ is an even permutation of $(1,2,3)$, $0$ if any index is repeated, and $(-1)$ if it is an odd permutation, and .

**Remark 8.3** The general fractional curl operator can be defined not only for $W = \mathbb{R}^3_+$, but also for regions $W \subset \mathbb{R}^3_+$, surfaces $S \subset \mathbb{R}^3_+$ and line $L \subset \mathbb{R}^3_+$. Note that the surface general fractional curl operator is used in the general fractional Stokes theorem to be given and proved in the following sections.

**Definition 8.12** Let $F(x,y,z)$ belong to the space $C^1_{-1}(\mathbb{R}^3_+)$.
Then the regional general fractional curl (the regional GF Curl) for the region
$$W := \{(x,y,z) : x \geq a \geq 0, \ y \geq c \geq 0, \ z \geq e \geq 0\}$$
is defined as
$$Curl^{(K)}_W F = e_x (Curl^{(K)}_W F)_x + e_y (Curl^{(K)}_W F)_y + e_z (Curl^{(K)}_W F)_z,$$
where
$$(Curl^{(K)}_W F)_x(x,y,z) = D^{(K_2)}_{[c,y]}[y'] F_z(x,y',z) - D^{(K_3)}_{[e,z]}[z'] F_y(x,y,z'),$$
$$(Curl^{(K)}_W F)_y(x,y,z) = D^{(K_3)}_{[e,z]}[z'] F_x(x,y,z') - D^{(K_1)}_{[a,x]}[x'] F_z(x',y,z),$$
$$(Curl^{(K)}_W F)_y(x,y,z) = D^{(K_3)}_{[a,x]}[x'] F_y(x',y,z) - D^{(K_1)}_{[c,y]}[y'] F_x(x,y',z)$$
with $x,y,z \geq 0$.

### 8.5.2 Surface General Fractional Curl Operator

We have given the definition of the regional GF Curl operator. Let us now define the surface GF operator. In this definition we will use the notations

$$(D_z^1 F_x)(x, y(x,z), z) := \left(\frac{\partial F_x(x,y,z)}{\partial z}\right)_{y=y(x,z)}, \quad (D_y^1 F_x)(x, y, z(x,y)) := \left(\frac{\partial F_x(x,y,z)}{\partial y}\right)_{z=z(x,y)},$$

$$(D_x^1 F_y)(x, y, z(x,y)) := \left(\frac{\partial F_y(x,y,z)}{\partial x}\right)_{z=z(x,y)}, \quad (D_z^1 F_y)(x(y,z), y, z) := \left(\frac{\partial F_y(x,y,z)}{\partial z}\right)_{z=z(x,y)},$$

$$(D_x^1 F_z)(x, y(x,z), z) := \left(\frac{\partial F_z(x,y,z)}{\partial x}\right)_{y=y(x,z)}, \quad (D_y^1 F_z)(x(y,z), y, z) := \left(\frac{\partial F_z(x,y,z)}{\partial y}\right)_{x=x(y,z)}.$$

Let us give a definition of the set $\mathbb{F}_S^1(\mathbb{R}_{0,+}^3)$ of vector fields $F$ that are used to define the surface GF Curl operator (compare with Definition 7.3 of the set $\mathbb{F}_S^1(\mathbb{R}_{0,+}^3)$ that is used in definition of the surface GFI).

**Definition 8.13** *Let $S$ be a piecewise simple surface ($S \in \mathbb{P}(\mathbb{R}_{0,+}^3)$).*
*Let the vector field*
$$F := e_x F_x(x,y,z) + e_y F_y(x,y,z) + e_z F_z(x,y,z). \tag{139}$$
*on the surface $S$ satisfy the conditions*

$$(D_z^1 F_x)(x, y_j(x,z), z) \in C_{-1}(\mathbb{R}_+^2), \quad (D_y^1 F_x)(x, y, z_k(x,y)) \in C_{-1}(\mathbb{R}_+^2),$$

$$(D_x^1 F_y)(x, y, z(x,y)) \in C_{-1}(\mathbb{R}_+^2), \quad (D_z^1 F_y)(x(y,z), y, z) \in C_{-1}(\mathbb{R}_+^2),$$

$$(D_x^1 F_z)(x, y(x,z), z) \in C_{-1}(\mathbb{R}_+^2), \quad (D_y^1 F_z)(x(y,z), y, z) \in C_{-1}(\mathbb{R}_+^2).$$

*for all $i = 1, \dots, n_x$, $j = 1, \dots, n_y$, $k = 1, \dots, n_z$.*

*The set of such vector fields $F$ on piecewise simple surface $S$ will be denoted by $\mathbb{F}_S^1(\mathbb{R}_{0,+}^3)$.*

Then the surface general fractional curl is defined in the following form. The piecewise simple surface is defined in 7.2.

**Definition 8.14** *Let $S$ be a piecewise simple surface ($S \in \mathbb{P}(\mathbb{R}_{0,+}^3)$) and a vector fields $F$ on this surface $S$ belongs to the set $\mathbb{F}_S^1(\mathbb{R}_{0,+}^3)$.*

*Then the surface general fractional curl (surface GF Curl) of the vector field $F(x,y,z)$ for the piecewise simple surface $S$ is defined as*

$$\mathrm{Curl}_S^{(K)} F = e_x (\mathrm{Curl}_S^{(K)} F)_x + e_y (\mathrm{Curl}_S^{(K)} F)_y + e_z (\mathrm{Curl}_S^{(K)} F)_z,$$

*where*

$$(\mathrm{Curl}_S^{(K)} F)_x = \sum_{i=1}^{n_x} I_{(K_2)}^{y,*}[y'] (D_y^1 F_z)(x_i(y',z), y', z) - I_{(K_3)}^{z,*}[z'] (D_z^1 F_y)(x_i(y,z'), y, z'),$$

$$(\mathrm{Curl}_S^{(K)} F)_y = \sum_{j=1}^{n_y} I_{(K_3)}^{z,*}[z'] (D_z^1 F_x)(x, y_j(x,z'), z') - I_{(K_1)}^{x,*}[x'] (D_x^1 F_z)(x', y_j(x',z), z),$$

$$(\mathrm{Curl}_S^{(K)} F)_z = \sum_{k=1}^{n_z} I_{(K_1)}^{x,*}[x'] (D_x^1 F_y)(x', y, z_k(x',y)) - I_{(K_2)}^{y,*}[y'] (D_y^1 F_x)(x, y', z_k(x,y'))$$

*with $x, y, z \geq 0$.*

## 8.6 Difficulties in Generalization of Stokes Theorem

In the standard Stokes theorem, the following equality should be satisfied
$$\int_{\partial S} P(x,y,z)\, dx = \int_S \left( \frac{\partial P(x,y,z)}{\partial z} dx\, dz - \frac{\partial P(x,y,z)}{\partial y} dx\, dy \right). \tag{140}$$
Equation (140) can be regarded as the Stokes formula for the vector field $F = e_x P(x, y.z)$.

If $S \subset \mathbb{R}^3_{0,+}$ is $Z$-simple and $Y$-simple surface that is described by the equation $z = z(x,y)$ for $(x,y) \in D_{xy} \subset \mathbb{R}^2_{0,+}$, and the equation $y = y(x,z)$ for $(x,z) \in D_{xz} \subset \mathbb{R}^2_{0,+}$, then
$$\int_S Q_z(x,y,z)\, dx\, dy = \int_{D_{xy}} Q(x,y,z(x,y))\, dx\, dy,$$

$$\int_S Q_y(x,y,z)\, dx\, dz = \int_{D_{xz}} Q_y(x,y(x,z),z)\, dx\, dz,$$
where $D_{xy}$ and $D_{xz}$ are the projections of the surface $S$ on the $XY$ and $XZ$-planes.

The values of the function $P(x,y,z)$ on the line $L \subset \mathbb{R}^3_{0,+}$ are equal to the values of the function $P(x,y,z(x,y))$ on the line $L_{xy} \subset \mathbb{R}^2_{0,+}$, which is the projection of the line $L$ onto the $XY$-plane
$$\int_{\partial S} P(x,y,z)\, dx = \int_{\partial D_{xy}} P(x,y,z(x,y))\, dx,$$
where $L = \partial S$ and $L_{xy} = \partial D_{xy}$.

Therefore, equality (140) means that
$$\int_{\partial D_{xy}} P(x,y,z(x,y))\, dx = \int_{D_{xz}} \frac{\partial P(x,y,z(x,y))}{\partial y} dx\, dz =$$
$$\int_{D_{xz}} \left( \left(\frac{\partial P(x,y,z)}{\partial y}\right)_{z=z(x,y)} + \left(\frac{\partial P(x,y,z)}{\partial z}\right)_{z=z(x,y)} \frac{z(x,y)}{\partial y} \right) dx\, dz =$$

$$\int_{D_{xz}} \left(\frac{\partial P(x,y,z)}{\partial z}\right)_{y=y(x,z)} dx\, dz - \int_{D_{xy}} \left(\frac{\partial P(x,y,z)}{\partial y}\right)_{z=z(x,y)} dx\, dy. \tag{141}$$

Equality (140) is based on the standard chain rule
$$\frac{\partial Q(x,y,z(x,y))}{\partial y} = \left(\frac{\partial Q(x,y,z)}{\partial y}\right)_{z=z(x,y)} + \left(\frac{\partial Q(x,y,z)}{\partial z}\right)_{z=z(x,y)} \frac{\partial z(x,y)}{\partial y}.$$
For the fractional derivatives and GFD, the standard chain is violated [3, pp.97-98].

For the general fractional Stokes theorem, we should have the equality
$$I^{(M_1)}_{\partial D_{xy}}[x]\, P(x,y,z) = I^{(M)}_{D_{xz}}[x,z](D^{z,*}_{K_3}[z']P(x,y,z'))_{y=y(x,z)} -$$

$$I^{(M)}_{D_{xy}}[x,y](D^{y,*}_{K_2}[y']P(x,y',z))_{z=z(x,y)}, \tag{142}$$
which is a fractional analog of equation (141), where $S$ is the $Z$-simple and $Y$-simple surface.

Identity (142) does not hold in the general case, due to the inequalities
$$I^{(M_2)}_{[c,d]}[y](D^{y,*}[y']P(x,y',z))_{z=z(x,y)} \neq (I^{(M_2)}_{[c,d]}[y]D^{y,*}[y']P(x,y',z))_{z=z(x,y)}.$$
and
$$I^{(M_2)}_{[c,d]}[y](D^{y,*}[y']P(x,y',z))_{z=z(x,y)} \neq I^{(M_2)}_{[c,d]}[y]D^{y,*}[y'](P(x,y',z))_{z=z(x,y)}.$$
We can state that we have the equality only if $Z$-simple surface $S$ is described by equation $z = z(x)$ for all $(x,y) \in D_{x,y}$. Then the equalities have the form

$$I_{[c,d]}^{(M_2)}[y](D^{y,*}[y']P(x,y',z))_{z=z(x)} = (I_{[c,d]}^{(M_2)}[y]D^{y,*}[y']P(x,y',z))_{z=z(x)} =$$

$$P(x,d,z(x)) - P(x,c,z(x)).$$

and

$$I_{[c,d]}^{(M_2)}[y](D^{y,*}[y']P(x,y',z))_{z=z(x)} =$$

$$I_{[c,d]}^{(M_2)}[y]D^{y,*}[y'](P(x,y',z))_{z=z(x)} = P(x,d,z(x)) - P(x,c,z(x)).$$

A similar situation for $Y$-simple surface that is described by the equation $y = y(x,z)$ with $x,z \in D_{xz}$.

As a result, for the regional GF curl $Curl_W^{(K)}$ with $W = \mathbb{R}_{0,+}^3$, the general fractional theorem can be proved only if the surface consists of faces parallel to the $XY$, $XZ$, $YZ$ planes.

## 8.7 General Fractional Stokes Theorem for Box without Bottom

Let us consider a box-shaped surface $S$ without a bottom, i.e. a parallelepiped surface without a bottom face (base). The parallelepiped can be described as

$$W := \{(x,y,z): \ 0 \leq a \leq x \leq b, \ 0 \leq c \leq y \leq d, \ 0 \leq e \leq z \leq f\}. \tag{143}$$

The vertices of the parallelepiped have the coordinates

$$A(a,c,e), \ A'(a,c,f),$$
$$B(b,c,e), \ B'(b,c,f),$$
$$C(b,d,e), \ C'(b,d,f),$$
$$D(a,d,e), \ D'(a,d,f).$$

The boundary of this parallelepiped is

$$\partial W = S \cup ABCD,$$

where the surface $ABCD$ is a bottom face of the parallelepiped region, the surface $S$ consists of the following faces

$$S = D_{xy} \cup D_{xz} \cup D_{yz}, \tag{144}$$

where

$$D_{xy} = A'B'C'D', \tag{145}$$

$$D_{xz} = ABB'A' \cup CDD'C', \tag{146}$$

$$D_{yz} = ADD'A' \cup BCC'B' \tag{147}$$

without the face $ABCD$.

The closed line $L = \partial S$, which is the boundary of the surface $S$ in the form of rectangle $ABCD$, consists of the segments

$$L = ABCD = AB \cup BC \cup CD \cup DA. \tag{148}$$

The surface GFI operator is desribed as

$$I_S^{(M)} := e_x I_{D_{yz}}^{(M)}[y,z] + e_y I_{D_{xz}}^{(M)}[x,z] + e_z I_{D_{xy}}^{(M)}[x,y].$$

Then

$$(I_S^{(M)}, F) = I_{D_{yz}}^{(M)}[y,z] F_x(x,y,z) + I_{D_{xz}}^{(M)}[x,z] F_y(x,y,z) + I_{D_{xy}}^{(M)}[x,y] F_z(x,y,z).$$

To calculate the expession $(I_S^{(M)}, Curl_S^{(K)} F)$, we should consider the following surface GFIs

$$(I_S^{(M)}, Curl_S^{(K)} F) =$$

$$I_{D_{yz}}^{(M)}[y,z](D_{(K_2)}^{y,*}[y']F_z - D_{(K_3)}^{z,*}[z']F_y) +$$

$$I_{D_{xz}}^{(M)}[x,z](D_{(K_3)}^{z,*}[z']F_x - D_{(K_1)}^{x,*}[x']F_z) +$$

$$I_{D_{xy}}^{(M)}[x,y](D_{(K_1)}^{x,*}[x']F_y - D_{(K_2)}^{y,*}[y']F_x).$$

**Theorem 8.1** *(General fractional Stokes theorem for parallelepiped surface without bottom)*
*Let $S$ is a smooth oriented surface (144) in $\mathbb{R}_{0,+}^3$, which is bounded by an oriented closed smooth line $L = \partial S$ that is given by (148), and $D_{yz}$, $D_{xz}$, $D_{xy}$ are domains (145) - (147), which are projections of the surface $S$ onto the YZ, XZ, XY planes.*

*Then, for the vector field $F(x,y,z) \in \mathbb{F}_{-1}^1(S)$, we have the equation*
$$(I_{\partial S}^{(M)}, F) = (I_S^{(M)}, Curl_S^{(K)} F),$$
*which is the general fractional Stokes equation.*

*For the surface $S$ in the form of the parallelepiped without a bottom face, the general fractional Stokes equation has the form*

$$I_{D_{yz}}^{(M)}[y,z](D_{(K_2)}^{y,*}[y']F_z - D_{(K_3)}^{z,*}[z']F_y) +$$

$$I_{D_{xz}}^{(M)}[x,z](D_{(K_3)}^{z,*}[z']F_x - D_{(K_1)}^{x,*}[x']F_z) +$$

$$I_{D_{xy}}^{(M)}[x,y](D_{(K_1)}^{x,*}[x']F_y - D_{(K_2)}^{y,*}[y']F_x) =$$

$$I_{AB}^{(M_1)}[x] F_x - I_{DC}^{(M_1)}[x] F_x + I_{BC}^{(M_2)}[y] F_y - I_{AD}^{(M_2)}[y] F_y.$$

*Proof.* YZ) Let us consider the GFI operator $I_{D_{yz}}^{(M)}$. The surface $D_{yz}$ consists of two $X$-simple surfaces $BCC'B'$ and $ADD'A'$. Therefore, we get
$$I_{D_{yz}}^{(M)}[y,z] = I_{BCC'B'}^{(M)}[y,z] - I_{ADD'A'}^{(M)}[y,z],$$
where we took into account the direction of the normals to the outer surface with respect to the direction of the basis vector $e_x$.

YZ1) For the first $X$-simple surface $BCC'B'$.
$$I_{BCC'B'}^{(M)}[y,z](D_{(K_2)}^{y,*}[y']F_z - D_{(K_3)}^{z,*}[z']F_y) =$$

$$I_{BCC'B'}^{(M)}[y,z] D_{(K_2)}^{y,*}[y'] F_z(b,y',z) - I_{BCC'B'}^{(M)}[y,z] D_{(K_3)}^{z,*}[z'] F_y(b,y,z') =$$

$$I_{[e,f]}^{(M_3)}[z] I_{[c,d]}^{(M_2)}[y] D_{(K_2)}^{y,*}[y'] F_z(b,y',z) - I_{[c,d]}^{(M_2)}[y] I_{[e,f]}^{(M_3)}[z] D_{(K_3)}^{z,*}[z'] F_y(b,y,z') =$$

$$I_{[e,f]}^{(M_3)}[z] (F_z(b,d,z) - F_z(b,c,z)) - I_{[c,d]}^{(M_2)}[y] (F_y(b,y,f) - F_y(b,y,e)) =$$

$$I_{[e,f]}^{(M_3)}[z]\,F_z(b,d,z) - I_{[e,f]}^{(M_3)}[z]\,F_z(b,c,z) - I_{[c,d]}^{(M_2)}[y]\,F_y(b,y,f) + I_{[c,d]}^{(M_2)}[y]\,F_y(b,y,e) =$$

$$I_{CC'}^{(M_3)}[z]\,F_z(x,y,z) - I_{BB'}^{(M_3)}[z]\,F_z(x,y,z) - I_{B'C'}^{(M_2)}[y]\,F_y(x,y,z) + I_{BC}^{(M_2)}[y]\,F_y(x,y,z).$$

YZ2) For the second $X$-simple surface $ADD'A'$.

$$-I_{ADD'A'}^{(M)}[y,z](D_{(K_2)}^{y,*}[y']F_z - D_{(K_3)}^{z,*}[z']F_y) =$$

$$-I_{ADD'A'}^{(M)}[y,z]\,D_{(K_2)}^{y,*}[y']\,F_z(a,y',z) + I_{ADD'A'}^{(M)}[y,z]\,D_{(K_3)}^{z,*}[z']\,F_y(a,y,z') =$$

$$-I_{[e,f]}^{(M_3)}[z]\,I_{[c,d]}^{(M_2)}[y]\,D_{(K_2)}^{y,*}[y']\,F_z(a,y',z) + I_{[c,d]}^{(M_2)}[y]\,I_{[e,f]}^{(M_3)}[z]\,D_{(K_3)}^{z,*}[z']\,F_y(a,y,z') =$$

$$-I_{[e,f]}^{(M_3)}[z]\,(F_z(a,d,z) - F_z(a,c,z)) + I_{[c,d]}^{(M_2)}[y]\,(F_y(a,y,f) - F_y(a,y,e)) =$$

$$-I_{[e,f]}^{(M_3)}[z]\,F_z(a,d,z) + I_{[e,f]}^{(M_3)}[z]\,F_z(a,c,z) + I_{[c,d]}^{(M_2)}[y]\,F_y(a,y,f) - I_{[c,d]}^{(M_2)}[y]\,F_y(a,y,e) =$$

$$-I_{DD'}^{(M_3)}[z]\,F_z(x,y,z) + I_{AA'}^{(M_3)}[z]\,F_z(x,y,z) + I_{A'D'}^{(M_2)}[y]\,F_y(x,y,z) - I_{AD}^{(M_2)}[y]\,F_y(x,y,z).$$

XZ) Let us consider the GFI operator $I_{D_{xz}}^{(M)}[x,z]$. The surface $D_{xz}$ consists of two $Y$-simple surfaces $ABB'A'$ and $CDD'C'$. Therefore, we get

$$I_{D_{xz}}^{(M)}[x,z] = -I_{ABB'A'}^{(M)}[x,z] + I_{CDD'C'}^{(M)}[x,z].$$

Here we took into account the direction of the normals to the outer surface with respect to the direction of the basis vector $e_y$. The minus in front of the operator $I_{ABB'A'}^{(M)}$ is due to the fact that the normal to the surface is directed in the opposite direction with respect to the vector $e_y$.

XZ1) For the first $Y$-simple surface $ABB'A'$.

$$-I_{ABB'A'}^{(M)}[x,z](D_{(K_3)}^{z,*}[z']F_x - D_{(K_1)}^{x,*}[x']F_z) =$$

$$-I_{ABB'A'}^{(M)}[x,z]\,D_{(K_3)}^{z,*}[z']F_x + I_{ABB'A'}^{(M)}[x,z]\,D_{(K_1)}^{x,*}[x']F_z =$$

$$-I_{[a,b]}^{(M_1)}[x]\,I_{[e,f]}^{(M_3)}\,D_{(K_3)}^{z,*}[z']F_x(x,c,z') + I_{[e,f]}^{(M_3)}\,I_{[a,b]}^{(M_1)}[x]\,D_{(K_1)}^{x,*}[x']F_z(x',c,z) =$$

$$-I_{[a,b]}^{(M_1)}[x]\,(F_x(x,c,f) - F_x(x,c,e)) + I_{[e,f]}^{(M_3)}\,(F_z(b,c,z) - F_z(a,c,z)) =$$

$$-I_{A'B'}^{(M_1)}[x']\,F_x(x,y,z) + I_{AB}^{(M_1)}[x']\,F_x(x,y,z) + I_{BB'}^{(M_3)}[z']\,F_z(x,y,z) - I_{AA'}^{(M_3)}[z']\,F_z(x,y,z).$$

XZ2) For the second $Y$-simple surface $CDD'C'$.

$$I_{CDD'C'}^{(M)}[x,z](D_{(K_3)}^{z,*}[z']F_x - D_{(K_1)}^{x,*}[x']F_z) =$$

$$I_{CDD'C'}^{(M)}[x,z]\,D_{(K_3)}^{z,*}[z']F_x - I_{CDD'C'}^{(M)}[x,z]\,D_{(K_1)}^{x,*}[x']F_z) =$$

$$I_{[a,b]}^{(M_1)}[x]\,I_{[e,f]}^{(M_3)}[z]\,D_{(K_3)}^{z,*}[z']F_x(x,d,z') - I_{[e,f]}^{(M_3)}[z]\,I_{[a,b]}^{(M_1)}[x]\,D_{(K_1)}^{x,*}[x']F_z(x',d,z)) =$$

$$I_{[a,b]}^{(M_1)}[x]\,(F_x(x,d,f) - F_x(x,d,e)) - I_{[e,f]}^{(M_3)}[z]\,(F_z(b,d,z) - F_z(a,d,z)) =$$

$$I_{D'C'}^{(M_1)}[x]\,F_x(x,y,z) - I_{DC}^{(M_1)}[x]\,F_x(x,y,z) - I_{CC'}^{(M_3)}[z]\,F_z(x,y,z) + I_{DD'}^{(M_3)}[z]\,F_z(x,y,z).$$

XY) Let us consider the GFI operator $I_{D_{xy}}^{(M)}[x,y]$. The surface $D_{xy}$ consists of one $Z$-simple surface $A'B'C'D'$.

$$I_{D_{xy}}^{(M)}[x,y] = I_{A'B'C'D'}^{(M)}[x,y].$$

Therefore, we get

$$I_{A'B'C'D'}^{(M)}[x,y](D_{(K_1)}^{x,*}[x']F_y - D_{(K_2)}^{y,*}[y']F_x) =$$

$$I_{[a,b]}^{(M_1)}[x]\,I_{[c,d]}^{(M_2)}[y](D_{(K_1)}^{x,*}[x']F_y(x',y,f) - D_{(K_2)}^{y,*}[y']F_x(x,y',f)) =$$

$$I_{[c,d]}^{(M_2)}[y]\,I_{[a,b]}^{(M_1)}[x]\,D_{(K_1)}^{x,*}[x']F_y(x',y,f) - I_{[a,b]}^{(M_1)}[x]\,I_{[c,d]}^{(M_2)}[y]\,D_{(K_2)}^{y,*}[y']F_x(x,y',f)) =$$

$$I_{[c,d]}^{(M_2)}[y]\,(F_y(b,y,f) - F_y(a,y,f)) - I_{[a,b]}^{(M_1)}[x]\,(F_x(x,d,f) - F_x(x,c,f)) =$$

$$I_{B'C'}^{(M_2)}[y]\,F_y(x,y,z) - I_{A'D'}^{(M_2)}[y]\,F_y(x,y,z) - I_{D'C'}^{(M_1)}[x]\,F_x(x,y,z) + I_{A'B'}^{(M_1)}[x]\,F_x(x,y,z).$$

XYZ) As a result, we get

$$(I_S^{(M)}, Curl_S^{(K)}F) =$$

$$I_{CC'}^{(M_3)}[z]\,F_z - I_{BB'}^{(M_3)}[z]\,F_z - I_{B'C'}^{(M_2)}[y]\,F_y + I_{BC}^{(M_2)}[y]\,F_y +$$

$$- I_{DD'}^{(M_3)}[z]\,F_z + I_{AA'}^{(M_3)}[z]\,F_z + I_{A'D'}^{(M_2)}[y]\,F_y - I_{AD}^{(M_2)}[y]\,F_y +$$

$$- I_{A'B'}^{(M_1)}[x]\,F_x + I_{AB}^{(M_1)}[x]\,F_x + I_{BB'}^{(M_3)}[z]\,F_z - I_{AA'}^{(M_3)}[z]\,F_z +$$

$$I_{D'C'}^{(M_1)}[x]\,F_x - I_{DC}^{(M_1)}[x]\,F_x - I_{CC'}^{(M_3)}[z]\,F_z + I_{DD'}^{(M_3)}[z]\,F_z +$$

$$I_{B'C'}^{(M_2)}[y]\,F_y - I_{A'D'}^{(M_2)}[y]\,F_y - I_{D'C'}^{(M_1)}[x]\,F_x + I_{A'B'}^{(M_1)}[x]\,F_x.$$

Bringing down similar terms, we obtain

$$(I_S^{(M)}, Curl_S^{(K)}F) =$$

$$I_{AB}^{(M_1)}[x]\,F_x - I_{DC}^{(M_1)}[x]\,F_x + I_{BC}^{(M_2)}[y]\,F_y - I_{AD}^{(M_2)}[y]\,F_y =$$

$$I_{L_x}^{(M_1)}[x]\,F_x + I_{L_y}^{(M_2)}[y]\,F_y = (I_L^{(M)}, F),$$

where $L = ABCD$ is the closed line in the $XY$-plane, which is the boundary $L = \partial S$ of the surface $S$.

## 8.8 General Fractional Stokes Theorem for Surface GF Curl

Let us prove the general fractional Stokes theorem for surface general fractional gradient, where surface consist of simple surfaces or surfaces parallel to the coordinate planes.

**Theorem 8.2** *(General Fractional Stokes Theorem for Surface GF Curl)*
Let $S \subset \mathbb{R}^3_{0,+}$ be a simple surface (or a piecewise simple surface $S \in \mathbb{P}(\mathbb{R}^3_{0,+})$), and the vector field $F$ belongs to the set $\mathbb{F}^1_S(\mathbb{R}^3_{0,+})$.
Then, the equation
$$(I^{(M_1)}_{\partial S}, F) = (I^{(M)}_S, Curl^{(K)}_S F)$$
holds. where $Curl^{(K)}_S$ is the surface GF Curl.

*Proof.* Let us prove the general Stokes theorem for the vector field $F_x = e_x P(x, y, z)$. The general Stokes equarions for the vector fields $F_y = e_y F_y(x, y, z)$ and $F_z = e_z F_z(x, y, z)$ are proved similarly.

Let $S \subset \mathbb{R}^3_{0,+}$ be a Z-simple and Y-simple surface that is described by the equation $z = z(x, y)$ for $(x, y) \in D_{xy} \subset \mathbb{R}^2_{0,+}$, and the equation $y = y(x, z)$ for $(x, z) \in D_{xz} \subset \mathbb{R}^2_{0,+}$.

The values of the function $P(x, y, z)$ on the line $L \subset \mathbb{R}^3_{0,+}$ are equal to the values of the function $P(x, y, z(x, y))$ on the line $L_{xy} \subset \mathbb{R}^2_{0,+}$, which is a projection of the line $L$ onto the $XY$-plane,
$$I^{(M_1)}_{\partial S}[x]\, P(x, y, z) = I^{(M_1)}_{\partial D_{xy}}[x] P(x, y, z(x, y)),$$
where $L = \partial S$ and $L_{xy} = \partial D_{xy}$.

Let us assume that $L_{xy}$ consists of two Y-simple lines $L_{xy1}$ and $L_{xy2}$ that are described by equations $y = y_1(x) \geq 0$ and $y = y_2(x)$, where $y_2(x) \geq y_1(x) \geq 0$. Then
$$I^{(M_1)}_{L_{xy}}[x] P(x, y, z(x, y)) = I^{(M_1)}_{L_{xy1}}[x] P(x, y, z(x, y)) + I^{(M_1)}_{L_{xy2}}[x] P(x, y, z(x, y)) =$$

$$I^{(M_1)}_{[a,b]}[x] P(x, y_1(x), z(x, y_1(x))) - I^{(M_1)}_{[a,b]}[x] P(x, y_2(x), z(x, y_2(x))) =$$

$$I^{(M_1)}_{[a,b]}[x](P(x, y_1(x), z(x, y_1(x))) - P(x, y_2(x), z(x, y_2(x)))). \tag{149}$$

Using the fundamental theorem of general fractional calculus, expression (149) can be written as
$$I^{(M_1)}_{[a,b]}[x](P(x, y_1(x), z(x, y_1(x))) - P(x, y_2(x), z(x, y_2(x)))) =$$

$$- I^{(M_1)}_{[a,b]}[x]\, I^{(M_2)}_{[y_1(x), y_2(x)]}[y] D^{y,*}_{(K_2)}[y']\, P(x, y', z(x, y')). \tag{150}$$

Then using the definition of the general fractional derivative and the standard chain rule for the first-order derivative, we get
$$- I^{(M_1)}_{[a,b]}[x]\, I^{(M_2)}_{[y_1(x), y_2(x)]}[y] D^{y,*}_{(K_2)}[y']\, P(x, y', z(x, y')) =$$

$$- I^{(M_1)}_{[a,b]}[x]\, I^{(M_2)}_{[y_1(x), y_2(x)]}[y] I^{y,*}_{(K_2)}[y']\, \frac{\partial P(x, y', z(x, y'))}{\partial y'} =$$

$$- I^{(M_1)}_{[a,b]}[x]\, I^{(M_2)}_{[y_1(x), y_2(x)]}[y] I^{y,*}_{(K_2)}[y']\, \left((\frac{\partial P(x, y', z)}{\partial y'})_{z=z(x,y')} + (\frac{\partial P(x, y', z)}{\partial z})_{z=z(x,y')} \frac{\partial z(x, y')}{\partial y'}\right) =$$

$$- I_{[a,b]}^{(M_1)}[x]\, I_{[y_1(x),y_2(x)]}^{(M_2)}[y]\, I_{(K_2)}^{y,*}[y']\, \left(\frac{\partial P(x,y',z)}{\partial y'}\right)_{z=z(x,y')} - \qquad (151)$$

$$- I_{[a,b]}^{(M_1)}[x]\, I_{[y_1(x),y_2(x)]}^{(M_2)}[y]\, I_{(K_2)}^{y,*}[y']\, \left(\frac{\partial P(x,y',z)}{\partial z}\right)_{z=z(x,y')} \frac{\partial z(x,y')}{\partial y'}. \qquad (152)$$

Let us use the property
$$I_{[y_1,y_2]}^{(M_2)}[y]\, I_{(K_2)}^{y,*}[y']\, U(y') = (M_2 * (K_2 * U))(y_2) - (M_2 * (K_2 * U))(y_1) =$$

$$((M_2 * K_2) * U)(y_2) - ((M_2 * K_2) * U)(y_1) =$$

$$(\{1\} * U)(y_2) - (\{1\} * U)(y_1) = \int_{y_1}^{y_2} U(y)\, dy$$

for the expression of term (152) in the form
$$I_{[y_1(x),y_2(x)]}^{(M_2)}[y]\, I_{(K_2)}^{y,*}[y']\, \left(\frac{\partial P(x,y',z)}{\partial z}\right)_{z=z(x,y')} \frac{\partial z(x,y')}{\partial y'} =$$

$$\int_{y_1(x)}^{y_2(x)} \left(\frac{\partial P(x,y',z)}{\partial z}\right)_{z=z(x,y')} \frac{\partial z(x,y')}{\partial y'}\, dy'. \qquad (153)$$

Then we can use the equation, which is used in the standard Stokes theorem, $(\{1\} * U)(z) = (M_3 * K_3 * U)(z)$, and the fundamental theorem of GFC in the form
$$\int_{y_1(x)}^{y_2(x)} \left(\frac{\partial P(x,y',z)}{\partial z}\right)_{z=z(x,y')} \frac{\partial z(x,y')}{\partial y'}\, dy' =$$

$$- \int_{z_1(x)}^{z_2(x)} \left(\frac{\partial P(x,y,z')}{\partial z'}\right)_{y=y(x,z')}\, dz' =$$

$$- I_{[z_1(x),z_2(x)]}^{(M_3)}[z]\, I_{(K_3)}^{z,*}[z']\, \left(\frac{\partial P(x,y,z')}{\partial z'}\right)_{y=y(x,z')}. \qquad (154)$$

Using (154), expression (152) can be written in the form
$$- I_{[a,b]}^{(M_1)}[x]\, I_{[y_1(x),y_2(x)]}^{(M_2)}[y]\, I_{(K_2)}^{y,*}[y']\, \left(\frac{\partial P(x,y',z)}{\partial z}\right)_{z=z(x,y')} \frac{\partial z(x,y')}{\partial y'} =$$

$$I_{[a,b]}^{(M_1)}[x]\, I_{[z_1(x),z_2(x)]}^{(M_3)}[z]\, I_{(K_3)}^{z,*}[z']\, \left(\frac{\partial P(x,y,z')}{\partial z'}\right)_{y=y(x,z')}. \qquad (155)$$

Therefore, for two terms (151) and (152), we obtain the equations
$$- I_{[a,b]}^{(M_1)}[x]\, I_{[y_1(x),y_2(x)]}^{(M_2)}[y]\, I_{(K_2)}^{y,*}[y']\, \left(\frac{\partial P(x,y',z)}{\partial y'}\right)_{z=z(x,y')} = I_{D_{xy}}^{(M)}[x,y]\, (Curl_S^{(K)} F_x)_z, \qquad (156)$$

$$I_{[a,b]}^{(M_1)}[x]\, I_{[z_1(x),z_2(x)]}^{(M_3)}[z]\, I_{(K_3)}^{z,*}[z']\, \left(\frac{\partial P(x,y,z')}{\partial z'}\right)_{y=y(x,z')} = I_{D_{xz}}^{(M)}[x,z]\, (Curl_S^{(K)} F_x)_y, \qquad (157)$$

for the vector field $F_x = e_x P(x,y,z)$.

As a result, we proved the general Stokes theorem for the surface FG Curl and the vector field $F_x = e_x P(x,y,z)$.
$$(I_{\partial S}^{(M_1)}, F_x) = (I_S^{(M)}, Curl_S^{(K)} F_x).$$
Equations for the remaining components of the vector field $F_y = e_y F_y(x,y,z)$ and $F_z = e_z F_z(x,y,z)$ are proved similarly.

# 9 General Fractional Gauss Theorem

## 9.1 Definition of Triple GFI by Iterated GFI

Let us define concept of the $Z$-simple region $W$ in $\mathbb{R}^3_{0,+}$.

**Definition 9.1** Let $W$ be region in $\mathbb{R}^3_{0,+}$ that is bounded above and below by smooth surfaces $S_{2,xy}$, $S_{1,xy}$ and a lateral surface $S_z$, whose generatrices are parallel to the $Z$-axis. Let surfaces $S_{1,xy}$, $S_{2,xy}$ be described by the equations
$$z = z_1(x,y) \geq 0, \quad z = z_2(x,y) \geq 0, \tag{158}$$
where the functions are continuous in the closed domain $W_{xy}$ that is a projection of the region $W$ onto the $XY$-plane, and $z_2(x,y) \geq z_1(x,y)$ for all $(x,y) \in W_{xy}$.

Then, the region $W$ will be called the $Z$-simple region (simple area along the $Z$-axis). The region $W$ is called simple, if $W$ is simple along three axes ($X$, $Y$, $Z$). If $W$ can be divided into a finite number of such regions with respect to all three axes, then $W$ will be called piecewise simple region in $\mathbb{R}^3_{0,+}$.

**Definition 9.2** Let $W \subset \mathbb{R}^3_{0,+}$ be $Z$-simple domain that is that is bounded above and below by smooth surfaces $S_{2,xy}$, $S_{1,xy}$ descrined by equation (158).
Let scalar field $f(x,y,z)$ be satisfy the condition
$$J_2(x,y) = I^{(M_3)}_{[z_1(x,y),z_2(x,y)]}[z]\, f(x,y,z) \in C_{-1}(\mathbb{R}^2_+). \tag{159}$$
Then the triple general fractional integral (triple GFI) is defined by the equation
$$I^{(M)}_W[x,y,z]\, f(x,y,z) := I^{(M)}_{W_{xy}}[x,y]\, I^{(M_3)}_{[z_1(x,y),z_2(x,y)]}[z]\, f(x,y,z),$$
where $I^{(M)}_{W_{xy}}[x,y]$ is double GFI.

**Definition 9.3** Let $W \subset \mathbb{R}^3_{0,+}$ be $Z$-simple domain that is that is bounded above and below by smooth surfaces $S_{2,xy}$, $S_{1,xy}$ descrined by equation (158). Let $W_{xy} \subset \mathbb{R}^2_{0,+}$ be projection of $W$ on the $XY$-plane such that $W_{xy}$ is $Y$-simple region in $XY$-plane that is bounded by the lines $y = y_2(x)$ and $y = y_1(x)$, where $y = y_2(x)$ and $y = y_1(x)$ are continuous functions on the interval $[a,b]$, $b > a \geq 0$ and $y_2(x) \geq y_1(x) > 0$ for all $x \in [a,b]$.
Let scalar fiels $f(x,y,z)$ satisfy condition (159) and
$$J_1(x) = I^{(M_2)}_{[y_1(x),y_2(x)]}[y]\, I^{(M_3)}_{[z_1(x,y),z_2(x,y)]}[z]\, f(x,y,z) \in C_{-1}(0,\infty). \tag{160}$$
Then the triple general fractional integral (triple GFI) is defined in the form
$$I^{(M)}_W[x,y,z]\, f(x,y,z) := I^{(M_1)}_{[a,b]}[x]\, I^{(M_2)}_{[y_1(x),y_2(x)]}[y]\, I^{(M_3)}_{[z_1(x,y),z_2(x,y)]}[z]\, f(x,y,z).$$

A volume general fractional integral (volume GFI) of a scalar field is a triple general fractional integral for the $Z$-simple region $W \subset \mathbb{R}^3_{0,+}$.

**Example 9.1** Using the parallelepiped
$$W := \{0 \leq x \leq a,\ 0 \leq y \leq b,\ 0 \leq z \leq c\},$$
the volume general fractional integral can be written as
$$I^{(M)}_W[x,y,z]\, f(x,y,z) = \int_0^a dx \int_0^b dy \int_0^c dz\, M_1(a-x)\, M_2(b-y)\, M_3(c-z)\, f(x,y,z).$$

For the kernels (62), the general fractional flux (161) has the form

$$I_W^{(M)}[x,y,z]\,f(x,y,z) = \int_0^a dx \int_0^b dy \int_0^c dz\,\frac{(a-x)^{\alpha_1-1}}{\Gamma(\alpha_1)}\,\frac{(b-y)^{\alpha_2-1}}{\Gamma(\alpha_2)}\,\frac{(c-z)^{\alpha_3-1}}{\Gamma(\alpha_3)}\,f(x,y,z). \quad (162)$$

For $\alpha_1 = \alpha_2 = \alpha_3 = 1$, equation (162) is gives as

$$I_W^{(M)}[x,y,z]\,f(x,y,z) = \int\int\int_W dV\,f(x,y,z) = \int_0^a dx \int_0^b dy \int_0^c dz\,f(x,y,z), \quad (163)$$

which is the standard volume integral for the function $f(x,y,z)$.

## 9.2 General Fractional Divergence

In this section, we give the definition of general fractional divergence for $\mathbb{R}^3_{0,+}$.

Let us define sets of vector fields that will be used in the definition of theregional general fractional divergence.

**Definition 9.4** Let $F(x,y,z)$ be a vector field that satisfies the conditions
$$D^{x,*}_{(K_1)}[x']F_x(x',y,z) \in C_{-1}(\mathbb{R}^3_+),\quad D^{y,*}_{(K_2)}[y']F_y(x,y',z) \in C_{-1}(\mathbb{R}^3_+),\quad D^{z,*}_{(K_3)}[z']F_z(x,y,z') \in C_{-1}(\mathbb{R}^3_+).$$
Then the set of such vector fields will be denoted as $\mathbb{F}^1_{-1,\text{Div}}(\mathbb{R}^3_+)$.

We can also consider vector fields $F(x,y,z)$ that belong to the function space $C^1_{-1}(\mathbb{R}^3_+)$.

**Definition 9.5** Let $F(x,y,z)$ be a vector field that satisfies the conditions
$$F_x(x,y,z), F_y(x,y,z), F_z(x,y,z) \in C^1_{-1}(\mathbb{R}^3_+).$$
Then the set of such vector fields will be denoted as $C^1_{-1}(\mathbb{R}^3_+)$.

In other words, the condition $F(x,y,z) \in C^1_{-1}(\mathbb{R}^3_+)$ means that all general fractional derivatives of all components of the vector field $F(x,y,z)$ with respect to all coordinated belong to space $C_{-1}(\mathbb{R}^3_+)$.

Let us define the general fractional divergence.

**Definition 9.6** Let $F(x,y,z)$ be a vector field that belongs to the set $\mathbb{F}^1_{-1,\text{Div}}(\mathbb{R}^3_+)$ or $F(x,y,z) \in C^1_{-1}(\mathbb{R}^3_+)$.

Then the general fractional divergence $\text{Div}^{(K)}_W$ for the region $W = \mathbb{R}^3_{0,+}$ is defined as
$$\text{Div}^{(K)}_W F = (D^{(K)}_W, F) =$$
$$D^{x,*}_{(K_1)}[x']F_x(x',y,z) + D^{y,*}_{(K_2)}[y']F_y(x,y',z) + D^{z,*}_{(K_3)}[z']F_z(x,y,z').$$

**Remark 9.1** The formula defining the operator can be written in compact form. If the vector field
$$F(x,y,z) = \sum_{k=1}^3 e_k\,F_k(x,y,z)$$
belongs to the function space $C^1_{-1}(\mathbb{R}^3_+)$, then the general fractional divergence for the region $W = \mathbb{R}^3_+$ is defined as

$$(Div_W^{(K)} F)(x,y,z) = (D_W^{(K)}, F) = \sum_{k=1}^{3} D_{(K_k)}^{x_k,*}[x'_k]F_k,$$

where

$$D_W^{(K)} := e_x D_{(K_1)}^{x,*}[x'] + e_y D_{(K_2)}^{y,*}[y'] + e_z D_{(K_3)}^{z,*}[z'].$$

**Remark 9.2** The general fractional divergence can be defined not only for $W = \mathbb{R}^3_{0,+}$, but also for regions $W \subset \mathbb{R}^3_{0,+}$, surfaces $S \subset \mathbb{R}^3_{0,+}$ and line $L \subset \mathbb{R}^3_{0,+}$. Note that the general fractional divergence for regions $W \subset \mathbb{R}^3_{0,+}$ is used in the general fractional Gausss theorem to be given and proved in the following sections.

**Definition 9.7** Let $F(x,y,z)$ be a vector field that belongs to the set $C^1_{-1}(\mathbb{R}^3_+)$, and the region $W$ be defined in the form
$$W := \{(x,y,z): \ x \geq a \geq 0, \ y \geq c \geq 0, \ z \geq e \geq 0\}.$$
Then the general fractional divergence $Div_W^{(K)}$ for the region $W$ is defined as
$$Div_W^{(K)} F = (D_W^{(K)}, F) =$$

$$D_{[a,x]}^{(K_1)}[x']F_x(x',y,z) + D_{[c,y]}^{(K_2)}[y']F_y(x,y',z) + D_{[e,z]}^{(K_3)}[z']F_z(x,y,z').$$

## 9.3 General Fractional Gauss Theorem for Z-Simple Region

The standard Gauss theorem (the Gauss-Ostrogradsky theorem) relates the flux of a vector field through a closed surface to the divergence of the field in the region enclosed. The Gauss theorem states that the surface integral of a vector field over a closed surface, which is the flux through the surface, is equal to the volume integral of the divergence over the region inside the surface.

Let us define a set of vector field, for which general fractional Gauss theorem will be formulated.

**Definition 9.8** Let $W$ in $\mathbb{R}^3_{0,+}$ be Z-simple region such that $W$ is a piecewise Y-simple and X-simple region
$$W = \bigcup_{k=1}^{n} W_{X,k} = \bigcup_{j=1}^{m} W_{Y,j},$$
where $W_{X,k}$ are the X-simple regions that is described by $x = x_{k,1}(y,z)$, $x = x_{k,2}(y,z)$ for $y,z \in D_{yz}$, and $W_{Y,j}$ are the Y-simple regions that are described by the functions $y = y_{j,1}(x,z)$, $y = y_{j,2}(x,z)$ for $(x,z) \in D_{xz}$.

Let vector field $F(x,y,z)$ satisfy the conditions
$$F_x(x_{k,2}(y,z),y,z), F_x(x_{k,1}(y,z),y,z) \in C_{-1}(\mathbb{R}^2_+),$$

$$F_y(x, y_{j,2}(x,z), z), F_y(x, y_{j,1}(x,z), z) \in C_{-1}(\mathbb{R}^2_+),$$

$$F_z(x,y,z_1(x,y)), F_z(x,y,z_2(x,y)) \in C_{-1}(\mathbb{R}^2_+).$$
Then the set of such vector fields $F$ will be denoted as $\mathbb{F}(x,y,z) \in \mathbb{F}_{-1}(\partial W)$.

$$f_x(x,y,z) := D_{(K_1)}^{x,*}[x'] F_x(x',y,z) \in C_{-1}(\mathbb{R}^3_+),$$

$$f_y(x,y,z) := D^{y,*}_{(K_2)}[y'] F_y(x,y',z) \in C_{-1}(\mathbb{R}^3_+),$$

$$f_z(x,y,z) := D^{z,*}_{(K_3)}[z'] F_z(x,y,z') \in C_{-1}(\mathbb{R}^3_+).$$

Then the set of such vector fields $F$ will be denoted as $\mathbb{F}(x,y,z) \in \mathbb{F}^1_{-1}(\mathbb{R}^3_+)$.

**Theorem 9.1** *(General fractional Gauss theorem for Z-simple region)*

*Let $W$ in $\mathbb{R}^3_{0,+}$ be Z-simple region such that $W$ is a piecewise Y-simple and X-simple region. Let $W$ is bounded above and below by smooth surfaces $S_{2,xy}$, $S_{1,xy}$, which is described by equations (158), and a lateral surface $S_z$, whose generatrices are parallel to the Z-axis.*

*Let the vector field $F(x,y,z)$ belongs to the sets $\mathbb{F}_{-1}(W)$ and $\mathbb{F}^1_{-1}(\mathbb{R}^3_+)$.*

*Then the general fractional Gauss equation has the from*
$$I^{(M)}_W[x,y,z](Div^{(K)}_W F) = (I_{\partial W}, F)$$

*that can be written as*
$$I^{(M)}_W[x,y,z](D^{x,*}_{(K_1)}[x'] F_x(x',y,z) + D^{y,*}_{(K_2)}[y'] F_y(x,y',z) + D^{z,*}_{(K_3)}[z'] F_z(x,y,z')) =$$

$$I^{(M)}_S[y,z] F_x(x,y,z) + I^{(M)}_S[x,y] F_y(x,y,z) + I^{(M)}_S[x,y] F_z(x,y,z).$$

*Proof.* Let us consider the triple GFI
$$I^{(M)}_W[x,y,z] = I^{(M)}_{W_{xy}}[x,y] I^{(M_3)}_{[z_1(x,y),z_2(x,y)]}[z]$$
for the function
$$f_z(x,y,z) = D^{z,*}_{(K_3)}[z'] F_z(x,y,z') \in C_{-1}(0,\infty),$$
where $F_z(x,y,z') \in C^1_{-1}(0,\infty)$ for each point $(x,y) \in D_{x,y}$.

Using the second fundamental theorem of GFC, we get
$$I^{(M_3)}_{[z_1(x,y),z_2(x,y)]}[z] D^{z,*}_{(K_3)}[z'] F_z(x,y,z') =$$

$$I^{z_2(x,y)}_{(M_3)}[z] D^{z,*}_{(K_3)}[z'] F_z(x,y,z') - I^{z_1(x,y)}_{(M_3)}[z] D^{z,*}_{(K_3)}[z'] F_z(x,y,z') =$$

$$(F_z(x,y,z_2(x,y)) - F_z(x,y,0)) - (F_z(x,y,z_2(x,y)) - F_z(x,y,0))) =$$

$$F_z(x,y,z_2(x,y)) - F_z(x,y,z_1(x,y)).$$

If we assume that
$$F_z(x,y,z_1(x,y)), F_z(x,y,z_2(x,y)) \in C_{-1}(\mathbb{R}^2_+),$$
then we get
$$I^{(M)}_W[x,y,z] D^{z,*}_{(K_3)}[z'] F_z(x,y,z') = I^{(M)}_{W_{xy}}[x,y](F_z(x,y,z_2(x,y)) - F_z(x,y,z_1(x,y))).$$

Therefore, the triple GFI can be represented through the surface GFIs in the form
$$I^{(M)}_W[x,y,z] D^{z,*}_{(K_3)}[z'] F_z(x,y,z') = I^{(M)}_{S_{2,xy}}[x,y] F_z(x,y,z) + I^{(M)}_{S_{1,xy}}[x,y] F_z(x,y,z),$$

where the surface GFIs are represented by the definition in the form
$$I^{(M)}_{S_{2,xy}}[x,y] F_z(x,y,z) = I^{(M)}_{W_{xy}}[x,y] F_z(x,y,z_2(x,y)),$$

$$I^{(M)}_{S_{1,xy}}[x,y] F_z(x,y,z) = I^{(M)}_{W_{xy}}[x,y] F_z(x,y,z_1(x,y)).$$

Then, taking into account that the surface GFI oves the surface $S_z$ is equati to zero
$$I_{S_z}^{(M)}[x,y]\, F_z(x,y,z) = 0,$$
we obtain
$$I_W^{(M)}[x,y,z]\, D_{(K_3)}^{z,*}[z']\, F_z(x,y,z') = I_S^{(M)}[x,y]\, F_z(x,y,z),$$
where $\partial W = S = S_{1,xy} \cup S_{2,xy} \cup S_z$ is closed surface that contains the region $W$ inside the surface $\partial W$.

If $W$ is piecewise $Y$-simple and $X$-simple region, such that
$$W = \bigcup_{k=1}^n W_{X,k} = \bigcup_{j=1}^m W_{Y,j},$$
where $W_{X,k}$ are the $X$-simple regions, and $W_{Y,j}$ are the $Y$-simple regions. Then, we have
$$W_{yz} = \bigcup_{k=1}^n W_{yz,k}, \quad W_{xz} = \bigcup_{j=1}^m W_{xz,j},$$
where $W_{yz}$ and $W_{xz}$ are projections of the region $W$ into $YZ$ and $XZ$ planes, $W_{yz,k}$ are projections of the $X$-simple region $W_{X,k}$, and $W_{xz,j}$ are projections of the $Y$-simple region $W_{Y,j}$.

Then the following equations are proved similarly
$$I_W^{(M)}[x,y,z]\, D_{(K_1)}^{x,*}[x']\, F_x(x',y,z) =$$
$$\sum_{k=1}^n I_{W_{yz,k}}^{(M)}[y,z]\, (F_x(x_{k,2}(y,z),y,z) - F_x(x_{k,1}(y,z),y,z)),$$

$$I_W^{(M)}[x,y,z]\, D_{(K_2)}^{y,*}[y']\, F_y(x,y',z) =$$
$$\sum_{j=1}^m I_{W_{xz,j}}^{(M)}[x,z]\, (F_y(x,y_{j,2}(x,z),z) - F_y(x,y_{j,1}(x,z),z)).$$

Therefore, we get
$$I_W^{(M)}[x,y,z]\, D_{(K_1)}^{x,*}[x']\, F_x(x',y,z) = I_S^{(M)}[y,z]\, F_x(x,y,z),$$

$$I_W^{(M)}[x,y,z]\, D_{(K_2)}^{y,*}[y']\, F_y(x,y',z) = I_S^{(M)}[x,y]\, F_y(x,y,z).$$

Using that general fractional divercence is defined as
$$Div_W^{(K)} F = D_{(K_1)}^{x,*}[x']\, F_x(x',y,z) + D_{(K_2)}^{y,*}[y']\, F_y(x,y',z) + D_{(K_3)}^{z,*}[z']\, F_z(x,y,z'),$$
we can get the general fractional Gauss equation.

As a result, we derive the equation
$$I_W^{(M)}[x,y,z]\, (Div_W^{(K)} F) = (I_{\partial W}, F).$$
that can be written as
$$I_W^{(M)}[x,y,z]\, (D_{(K_1)}^{x,*}[x']\, F_x(x',y,z) + D_{(K_2)}^{y,*}[y']\, F_y(x,y',z) + D_{(K_3)}^{z,*}[z']\, F_z(x,y,z')) =$$

$$I_S^{(M)}[y,z]\, F_x(x,y,z) + I_S^{(M)}[x,y]\, F_y(x,y,z) + I_S^{(M)}[x,y]\, F_z(x,y,z).$$

**Remark 9.3** The general fractional Gauss theorem is proved similarly for region $W \subset \mathbb{R}^3_{0,+}$ that can be represented as unions of the $Z$-simple regions $W_k$, which are piecewise $Y$-simple and $X$-simple regions in $\mathbb{R}^3_{0,+}$.

## 9.4 General Fractional Gauss Theorem for Parallelepiped

The standard Gauss theorem (the Gauss-Ostrogradsky theorem) states the folllwing. Let $W$ be a region in $\mathbb{R}^3_{0,+}$ with boundary $\partial W$. Then the volume integral of the divergence of vector field $F$ over $W$ and the surface integral of $F$ over the boundary $\partial W$ are related by

$$(I^1_{\partial W}, F) = I^1_W \, divF. \tag{164}$$

For the parallelepiped region, a general fractional Gauss theorem can be formulated in the following form.

**Theorem 9.2** *(Fractional Gauss's Theorem for a Parallelepiped)*
Let $F_x(x,y,z)$, $F_y(x,y)$, $F_z(x,y,z)$ belong to the function space $\mathbb{F}^1_{-1}(\mathbb{R}^3_+)$, and the region $W \subset \mathbb{R}^3_{0,+}$ has the form of the parallelepiped

$$W := \{(x,y,z): \; 0 \leq a \leq x \leq b, \; 0 \leq c \leq y \leq d, \; 0 \leq e \leq z \leq f\}. \tag{165}$$

If the boundary of $W$ be a closed surface $S = \partial W$, then

$$(I^{(M)}_{\partial W}, F) = I^{(M)}_W \, Div^{(K)}_W F. \tag{166}$$

*Proof.* For Cartesian coordinates, we have the vector field $F = F_x e_x + F_y e_y + F_z e_z$, and the GFI operators

$$I^{(M)}_W = I^{(M)}_W[x,y,z], \quad I^{(M)}_{\partial W} = e_x I^{(M)}_{\partial W}[y,z] + e_y I^{(M)}_{\partial W}[x,z] + e_z I^{(M)}_{\partial W}[x,y]. \tag{167}$$

Then

$$(I^{(M)}_{\partial W}, F) = I^{(M)}_{S_{yz}}[y,z]F_x + I^{(M)}_{S_{xz}}[x,z]F_y + I^{(M)}_{S_{xy}}[x,y]F_z, \tag{168}$$

and

$$I^{(M)}_W \, Div^{(K)}_W F = I^{(M)}_W[x,y,z](D^{(K)}_S[x]F_x + D^{(K)}_S[y]F_y + D^{(K)}_S[z]F_z), \tag{169}$$

where $S_{xy}$, $S_{xz}$, $S_{yz}$ are projections of $S = \partial W$ into $XY$, $XZ$, $YZ$ planes.

If $W$ is parallelepiped (165), then GFI operators (167) are

$$I^{(M)}_W[x,y,z] = I^{(M_1)}_{[a,b]}[x] \, I^{(M_2)}_{[c,d]}[y] \, I^{(M_3)}_{[e,f]}[z], \tag{170}$$

and

$$I^{(M)}_{\partial W}[y,z] = I^{(M_2)}_{[c,d]}[y] \, I^{(M_3)}_{[e,f]}[z], \tag{171}$$

$$I^{(M)}_{\partial W}[x,z] = I^{(M_1)}_{[a,b]}[x] \, I^{(M_3)}_{[e,f]}[z], \tag{172}$$

$$I^{(M)}_{\partial W}[x,y] = I^{(M_1)}_{[a,b]}[x] \, I^{(M_2)}_{[c,d]}[y]. \tag{173}$$

Using the fundanental theorem of GFC, we can realize the following transformations

$$(I^{(M)}_{\partial W}, F) = I^{(M)}_{\partial W}[y,z]F_x + I^{(M)}_{\partial W}[z,x]F_y + I^{(M)}_{\partial W}[x,y]F_z =$$

$$I^{(M_2)}_{[c,d]}[y] \, I^{(M_3)}_{[e,f]}[z] \{F_x(b,y,z) - F_x(a,y,z)\} +$$

$$I^{(M_1)}_{[a,b]}[x] \, I^{(M_3)}_{[e,f]}[z] \{F_y(x,d,z) - F_y(x,c,z)\} +$$

$$I^{(M_1)}_{[a,b]}[x]\, I^{(M_2)}_{[c,d]}[y]\, \{F_z(x,y,f) - F_z(x,y,e)\} =$$

$$I^{(M_1)}_{[a,b]}[x]\, I^{(M_2)}_{[c,d]}[y]\, I^{(M_3)}_{[e,f]}[z]\, \{D^{x,*}_{(K_1)}[x']F_x(x',y,z) +$$

$$D^{y,*}_{(K_2)}[y']F_y(x,y',z) + D^{z,*}_{(K_3)}[z']F_z(x,y,z')\} =$$

$$I^{(M)}_W(D^{(M)}_W, F) = I^{(M)}_W Div^{(K)}_W F.$$

This ends the proof of the fractional Gauss's formula for parallelepiped region.

# 10 Equalities for General Fractional Vector Operators and Interpretations

## 10.1 Equalities for General Fractional Differential Vector Operators

Let us give the basic relations for the differential vector operators of general fractional calculus.

1. For the scalar field $U = U(x,y,x)$, we have
$$Div^{(K)}_W\, Grad^\alpha_W\, U = \sum_{j=1}^{3} D^{*,x_j}_{(K_j)}[x'_j]\, D^{*,x'_j}_{(K_j)}[x''_j]\, U. \quad (174)$$

In the general case,
$$D^{*,x_j}_{(K_j)}[x'_j]\, D^{*,x'_j}_{(K_j)}[x''_j]\, U \neq D^{*,x_j}_{(K_j * K_j)}[x'_j]\, U. \quad (175)$$

This is due to the fact that the general fractional integral and the first-order derivative (in definitions of GFDs) do not commute
$$D^{*,x}_{(K)}[x'] = I^x_{(K)}[x']\, f^{(1)}(x') \neq \frac{d}{dx} I^x_{(K)}[x']\, f(x') = D^x_{(K)}[x']. \quad (176)$$

For example, in the case of power-law kernels $K(x) = h_{1-\alpha}(x)$, we have [30] the inequality
$$\left(D^\alpha_{C,0+}\right)^2 \neq D^{2\alpha}_{C,0+}. \quad (177)$$

As a result, we have the definition of the "scalar" general fractional Laplacian
$$\Delta^{(K)}_W U = , Div^{(K)}_W\, Grad^\alpha_W\, U. \quad (178)$$

2. For the vector field $F = F(x,y,x)$, we have
$$Grad^\alpha_W\, Div^{(K)}_W\, F = \sum_{m=1}^{3} \sum_{j=1}^{3} e_m D^{*,x_m}_{(K_m)}[x'_m]\, D^{*,x'_j}_{(K_j)}[x''_j]\, F_j. \quad (179)$$

3. The second relation for the scalar field $U = U(x,y,z)$ is
$$Curl^{(K)}_W\, Grad^{(K)}_W U = e_l\, \varepsilon_{lmn} D^{*,x_m}_W[x_m]\, D^{*,x_n}_W[x_n] U = 0, \quad (180)$$
where $\varepsilon_{lmn}$ is Levi-Civita symbol, i.e. it is $1$ if $(i,j,k)$ is an even permutation of $(1,2,3)$, $(-1)$ if it is an odd permutation, and $0$ if any index is repeated. The fulfillment of equality (180) is due to the properties of the Levi-Civita symbol.

4. For the vector field $F = e_m F_m$, it is easy to prove the relation
$$Div_W^{(K)} Curl_W^{(K)} F(x,y,z) = \sum_{k=1}^{3} \sum_{m=1}^{3} D_{(K_k)}^{*,x_k}[x'_k] \, \varepsilon_{klm} \, D_{(K_m)}^{*,x_m}[x'_m] F_m(x,y,z) =$$

$$\sum_{k=1}^{3} \sum_{m=1}^{3} \varepsilon_{klm} D_{(K_k)}^{*,x_k}[x'_k] D_{(K_m)}^{*,x_m}[x'_m] F_m(x,y,z) = 0, \tag{181}$$

where we use antisymmetry of $\varepsilon_{klm}$ with respect to $l$ and $m$. Equality (181) is also satisfied by the properties of the Levi-Civita symbol.

5. There exists a relation for the double curl operator in the form
$$Curl_W^{(K)} Curl_W^{(K)} F(x,y,z) = e_l \varepsilon_{lmn} D^{*,x_m}(K_m)[x'_m] \, \varepsilon_{npq} D_{(K_p)}^{*,x_p}[x'_p] F_q(x,y,z) =$$

$$e_l \varepsilon_{lmn} \varepsilon_{npq} D^{*,x_m}(K_m)[x'_m] D_{(K_p)}^{*,x_p}[x'_p] F_q(x,y,z). \tag{182}$$

Using the relation
$$\varepsilon_{lmn}\varepsilon_{lpq} = \delta_{mp}\delta_{nq} - \delta_{mq}\delta_{np}, \tag{183}$$

we obtain
$$Curl_W^{(K)} Curl_W^{(K)} F = Grad_W^{(K)} Div_W^{(K)} F - \Delta_W^{(K)} F. \tag{184}$$

As a result, we have the definition of the "vector" general fractional Laplacian
$$\Delta_W^{(K)} F = , Grad_W^{(K)} Div_W^{(K)} F - Curl_W^{(K)} Curl_W^{(K)} F \tag{185}$$

that is a generalization of the standard vector Laplacian [67, 68].

6. In the general case,
$$D_{(K)}^{*,x}[x'](f(x')g(x')) \neq (D_{(K)}^{*,x}[x']f(x')) \, g(x) + f(x) \, (D_{(K)}^{*,x}[x']g(x')). \tag{186}$$

For example (see Theorem 15.1 of [1]), if $f(x)$ and $g(x)$ are analytic functions on $[a,b]$ and the kernel $K(x) = h_{1-\alpha}(x)$, then product rule for the Riemann-Liouville fractional deirvative has the form
$$D_{RL,0+}^{x}[x'](f(x')g(x')) = \sum_{j=0}^{\infty} a(\alpha,j)(D_{RL,0+}^{\alpha-j}[x']f(x'))(g^{(j)}(x)), \tag{187}$$

where
$$a(\alpha,j) = \frac{\Gamma(\alpha+1)}{\Gamma(j+1)\Gamma(\alpha-j+1)}.$$

As a result, we have
$$Grad_W^{(K)}(f \, g) \neq (Grad_W^{(K)} f) \, g + f \, (Grad_W^{(K)} g), \tag{188}$$

$$Div_W^{(K)}(U \, F) \neq (Grad_W^{(K)} U) \, F + U \, Div_W^{(K)} F. \tag{189}$$

$$Curl_W^{(K)}(U \, F) \neq (Grad_W^{(K)} U) \, F + U \, Curl_W^{(K)} F. \tag{190}$$

These relations state that we cannot use the standard Leibniz rule (the product rule) in the general fractional vector calculus.

## 10.2  Physical Interpretations of General Fractional Differential Vector Operators

The general fractional Grad, Curl, and Div operators (in the Regional, Surface and Line form) are nonlocal characteristics of scalar and vector fields in non-local media and continua. Therefore, these characteristics depend on the region, surface and line that are used in considerations in models of non-local media and systems.

Let us emphasize that the fractional vector analogs of fundamental theorems (such as the gradient, Stock's and Gauss theorems) are not fulfilled for all type (regional, surface and line) of the general fractional vector operators (the gradient, curl and divergence). In the general case, the general fractional (GF) gradient theorem should be considered for the line GF Gradient, the FG Stock's theorem should be considered for the surface GF Curl operator, and the GF Green theorem should be considered for the regional GF Divergence. This is due to violation of the chain rule for general fractional derivatives.

Let us give some basic interpretations of general fractional divergence, gradient and curl operators.

1) The general fractional gradient specifies the direction of the fastest increase in the scalar field $U(x,y,z)$ in a nonlocal medium. The length of the vector of the line general fractional gradient is equal to the rate of increase of this field in this direction, where this increase takes into account changes in this field at the previous points of the line. In other words, the general fractional gradient describes the increase of scalar field, that takes into account all values of all velocities of the highest integer orders (jet) with some weights [69].

2) The general fractional curl operators can be interpreted in the following representation. If the velocity field of a non-local medium is used as a vector field in a region, on a surface or a line, the Curl operators characterize the rotational component of this vector field. The regional, surface and line operators Curl characterize the degree (measure) of non-potentiality of a vector field in nonlocal media and systems. The equality to zero of the general fractional curl of a vector field in a certain region, surface or line means the potentiality of the vector field in a nonlocal medium.

3) The general fractional divergence characterizes the degree (measure) of non-solenoidality of a vector field in a nonlocal media and systems. The general fractional divergence in a certain area of space characterizes the amount of a nonlocal medium that arises or disappears within the considered area per unit time. We can say that the general fractional divergence of the vector field describes the power of sources and sinks in a nonlocal medium. The equality to zero of the general fractional divergence in a certain region of space means the absence of sources and sinks of the vector field in a nonlocal medium. In other words, the media or fields in the considered region of space do not disappear and does not arise.

# 11 General FVC for Orthogonal Curvilinear Coordinates

## 11.1 Orthogonal Curvilinear Coordinates (OCC)

Curvilinear system of coordinates, or curvilinear coordinates, is a coordinate system in the Euclidean space. In standard vector calculus. the curvilinear coordinates are usually used on a plane ($n = 2$) and in space ($n = 3$). For such systems, the coordinate lines may be curved.

In Euclidean space, the use of orthogonal curvilinear coordinates (OCC) is of particular importance, since the formulas related look simpler in orthogonal coordinates than in the general case. The orthogonality can simplify the calculations. The well-known examples of such curvilinear coordinate systems in three-dimensional Euclidean space are cylindrical and spherical coordinates.

The curvilinear coordinates $(q_1, q_2, q_3)$ may be derived from a set of Cartesian coordinates $(x, y, z)$ by using nonlinear coordinate transformations. It should be emphasized that the violation of the standard chain rule leads to the fact that general fractional vector operators defined in different coordinate systems (Cartesian, cylindrical and spherical) are not related to each other by coordinate transformations. Because of this, it is impossible to obtain general fractional vector integral and differential operators in spherical and cylindrical coordinates by using the coordinate transformation. Therefore, the definitions of general fractional vector operators should be formulated separately.

The specific form of standard vector differential operators may differ, but these forms will be equivalent due to the Leibniz (product) rule. For fractional and generalized fractional calculus, the standard Leibniz rule does not hold. Because of this, such forms of notation cannot be equivalent. In this case, the specific form of the cylindrical and spherical general fractional operators must be such that the theorems of Green, Stock and Gauss hold.

It is known that the expressions for the gradient, divergence, curl and line, surface and volume integrals can be directly expressed. For orthogonal curvilinear coordinates (OCC), these integral and differential operators of the vector calculus can be expressed through the functions:

$$H_i(q) = \sqrt{\left(\frac{\partial x}{\partial q_1}\right)^2 + \left(\frac{\partial y}{\partial q_2}\right)^2 + \left(\frac{\partial z}{\partial q_3}\right)^2}.$$

The positive values $H_i = H_i(q) = H_i(q_1, q_2, q_3)$, which depend on a point in $\mathbb{R}^3_{0,+}$, are called the Lame coefficients or scale factors.

For the cylindrical coordinates $q_1 = r$, $q_2 = \phi$, $q_3 = z$, the Lame coefficients are

$$H_1 = H_r = 1, \quad H_2 = H_\phi = r, \quad H_3 = H_z = 1,$$

where $r \in [0, \infty)$, $\phi \in [0, 2\pi]$, $z \in \mathbb{R}$.

For the spherical coordinates $q_1 = r$, $q_2 = \theta$, $q_3 = \phi$, the Lame coefficients are

$$H_1 = H_r = 1, \quad H_2 = H_\theta = r, \quad H_3 = H_\phi = r \sin\theta,$$

where $r \in [0, \infty)$, $\theta \in [0, \pi]$, $\phi \in [0, 2\pi]$.

## 11.2 General Fractional Vector Differential Operators in OCC

In this section, we give only equations that will be used in definitions of the general fractional gradient, divergence and curl operators in OCC, and their examples for spherical and cylindrical coordinates. Complete definitions of these operators with function sets for which these operators are defined will be given in the following sections. In this section, formulas will be given only for regional GF vector differential operators. Definitions of line and surface general fractional gradient, divergence and curl operators in OCC will be given in the following sections.

### 11.2.1 General Fractional Grad, Div, Curl in OCC

The regional general fractional gradient for orthogonal curvilinear coordinates (the GF Gradient in OCC) is expressed in the form

$$Grad_W^{(K)} U = e_1 (Grad_W^{(K)} U)_1 + e_2 (Grad_W^{(K)} U)_2 + e_3 (Grad_W^{(K)} U)_3$$

where

$$(Grad_W^{(K)} U)_1 := \frac{1}{H_1(q)} D_{(K_1)}^{q_1,*}[q'_1] U(q'_1, q_2, q_3),$$

$$(Grad_W^{(K)} U)_2 := \frac{1}{H_2(q)} D_{(K_2)}^{q_2,*}[q'_2] U(q_1, q'_2, q_3),$$

$$(Grad_W^{(K)} U)_3 := \frac{1}{H_3(q)} D_{(K_3)}^{q_3,*}[q'_3] U(q_1, q_2, q'_3),$$

where $H_k = H_k(q) = H_k(q_1, q_2, q_3)$, and $e_1 = e_{q_1}$, $e_2 = e_{q_2}$, $e_3 = e_{q_3}$, and the function $U(q)$ belongs to the set $C_{-1}^1(\mathbb{R}_+^3)$.

The regional general fractional divergence in orthogonal curvilinear coordinates (the GF Divergence in OCC) is expressed in the form

$$Div_W^{(K)} F = (Div_W^{(K)} F)_1 + (Div_W^{(K)} F)_2 + (Div_W^{(K)} F)_3,$$

where

$$(Div_W^{(K)} F)_1 := \frac{1}{H_1 H_2 H_3} D_{(K_1)}^{q_1,*}[q'_1] (F_1 H_2 H_3)(q'_1, q_2, q_3),$$

$$(Div_W^{(K)} F)_2 := \frac{1}{H_1 H_2 H_3} D_{(K_1)}^{q_1,*}[q'_1] (F_2 H_1 H_2)(q'_1, q_2, q_3),$$

$$(Div_W^{(K)} F)_3 := \frac{1}{H_1 H_2 H_3} D_{(K_1)}^{q_1,*}[q'_1] (F_3 H_1 H_2)(q'_1, q_2, q_3),$$

if $F(q) = F_1(q) e_1 + F_2(q) e_2 + F_3(q) e_3$ belong to the space $C_{-1}^1(\mathbb{R}_+^3)$.

The regional general fractional curl in orthogonal curvilinear coordinates (the GF Curl in OCC) is expressed by the eqations

$$Curl_W^{(K)} F = e_1 (Curl_W^{(K)} F)_1 + e_2 (Curl_W^{(K)} F)_2 + e_3 (Curl_W^{(K)} F)_3,$$

where

$$(Curl_W^{(K)} F)_1 = \frac{1}{H_2 H_3} \left( D_{(K_2)}^{q_2,*}[q'_2] (F_3 H_3) - D_{(K_3)}^{q_3,*}[q'_3] (F_2 H_2) \right),$$

$$(Curl_W^{(K)} F)_2 = \frac{1}{H_1 H_3} \left( D_{(K_3)}^{q_3,*}[q'_3] (F_1 H_1) - D_{(K_1)}^{q_1,*}[q'_1] (F_3 H_3) \right),$$

$$(Curl_W^{(K)} F)_3 = \frac{1}{H_1 H_2} \left( D_{(K_1)}^{q_1,*}[q'_1] (F_2 H_2) - D_{(K_2)}^{q_2,*}[q'_2] (F_1 H_1) \right).$$

### 11.2.2 General Fractional Grad, Div, Curl in Spherical Coorditates

The Cartesian coordinates $(x, y, z)$ and spherical coordinates $(r, \theta, \varphi)$ are connected by the equations

$$\begin{aligned} x &= r\sin\theta\,\cos\varphi, \\ y &= r\sin\theta\,\sin\varphi, \\ z &= r\cos\theta, \end{aligned} \qquad (191)$$

where $r \in [0, \infty)$ is the length of the radial vector connecting the origin to the point $P(x, y, z)$, $\theta \in [0, \pi]$ is the polar angle, $\varphi \in [0, 2\pi]$ is the azimuthal angle. The basic vectors of these coordinate systems

$$\begin{aligned} e_x &= e_r \sin\theta\,\cos\varphi + e_\theta \cos\theta\,\cos\varphi - e_\varphi \sin\varphi \\ e_y &= e_r \sin\theta\,\sin\varphi + e_\theta \cos\theta\,\sin\varphi + e_\varphi \cos\varphi \\ e_z &= e_r \cos\theta - e_\theta \sin\theta. \end{aligned}$$

Let us emphasize that the violation of the standard chain rule leads to the fact that general fractional vector operators defined in Cartesian, and spherical coordinate are not related to each other by coordinate transformations. Because of this, it is impossible to obtain the fractional vector integral and differential operators in spherical coordinated by using coordinate transformation (191). Therefore, the definitions of general fractional vector operators should be formulated separately.

Let us define the regions

$$Ball_+^3 = \{(r, \theta, \varphi): r \in (0, \infty), \theta \in (0, \pi], \varphi \in (0, 2\pi]\},$$

$$Ball_{0,+}^3 = \{(r, \theta, \varphi): r \in [0, \infty), \theta \in [0, \pi], \varphi \in [0, 2\pi]\}.$$

We can see that

$$Ball_+^3 \subset \mathbb{R}^3 +, \quad Ball_{0,+}^3 \subset \mathbb{R}_{0,+}^3.$$

We will consider the vector field

$$F = e_r\, F_r(r, \theta, \varphi) + e_\theta\, F_\theta(r, \theta, \varphi) + e_\varphi\, F_\varphi(r, \theta, \varphi)$$

that belongs to $\mathbb{F}_{-1}^1(\mathbb{R}_+^3)$.

The general fractional gradient, divergence and curl operators in spherical coordinates are defined in the following form. For simplicity, we only present the definitions of regional operators. Linear and surface general fractional operators in spherical coordinates are defined similarly.

**Definition 11.1** *If $F(r, \theta, \varphi) \in C_{-1}^1(\mathbb{R}_+^3)$, then the regional general fractional gradient in spherical coordinates for the region $W = Ball_{0,+}^3$ is defined as*

$$Grad_W^{(K)} F(r, \theta, \varphi) = e_r\, (Grad_W^{(K)} F)_r + e_\theta\, (Grad_W^{(K)} F)_\theta + e_\varphi\, (Grad_W^{(K)} F)_\varphi,$$

*where*

$$(Grad_W^{(K)} F)_r(r, \theta, \varphi) = D_{(K_r)}^{r,*}[r']\, F(r', \theta, \varphi)$$

$$(Grad_W^{(K)} F)_\theta(r, \theta, \varphi) = \frac{1}{r}\, D_{(K_\theta)}^{\theta,*}[\theta']\, F(r, \theta', \varphi),$$

$$(Grad_W^{(K)} F)_\varphi(r, \theta, \varphi) = \frac{1}{r\sin\theta}\, D_{(K_\varphi)}^{\varphi',*}[\varphi']\, F(r, \theta, \varphi').$$

**Definition 11.2** *If $F(r, \theta, \varphi) \in C_{-1}^1(\mathbb{R}_+^3)$, then the regional general fractional divergence in spherical coordinates for the region $W = Ball_{0,+}^3$ is defined as*

$$Div_W^{(K)} F = (Div_W^{(K)} F)_r + (Div_W^{(K)} F)_\theta + (Div_W^{(K)} F)_\varphi,$$

*where*

$$(Div_W^{(K)} F)_r(r,\theta,\varphi) = \frac{1}{r^2} D_{(K_r)}^{r,*}[r']\,((r')^2\,F_r(r',\theta,\varphi))$$

$$(Div_W^{(K)} F)_\theta(r,\theta,\varphi) = \frac{1}{r\sin\theta} D_{(K_\theta)}^{\theta,*}[\theta']\,(F_\theta(r,\theta',\varphi)\sin\theta'),$$

$$(Div_W^{(K)} F)_\varphi(r,\theta,\varphi) = \frac{1}{r\sin\theta} D_{(K_\varphi)}^{\varphi,*}[\varphi']\,F_\varphi(r,\theta,\varphi').$$

Let us note the violations of the standard product (Leibniz) rule. For example, we have the inequality

$$D_{(K_\theta)}^{\theta,*}[\theta']\,(F_\theta(r,\theta',\varphi)\sin\theta') \neq$$

$$(D_{(K_\theta)}^{\theta,*}[\theta']\,F_\theta(r,\theta',\varphi))\sin\theta + (F_\theta(r,\theta',\varphi)\,D_{(K_\theta)}^{\theta,*}[\theta']\sin\theta').$$

As a second example

$$\frac{1}{r^2} D_{(K_r)}^{r,*}[r']\,((r')^2\,F_r(r',\theta,\varphi)) \neq$$

$$D_{(K_r)}^{r,*}[r']\,F_r(r',\theta,\varphi) + \frac{2}{r} F_\theta(r,\theta,\varphi).$$

**Definition 11.3** *If $F(r,\theta,\varphi) \in C^1_{-1}(\mathbb{R}^3_+)$, then the regional general fractional curl operator in spherical coordinates for the region $W = Ball^3_{0,+}$ is defined as*

$$Curl_W^{(K)} F = e_r\,(Curl_W^{(K)} F)_r + e_\theta\,(Curl_W^{(K)} F)_\theta + e_\varphi\,(Curl_W^{(K)} F)_\varphi,$$

*where*

$$(Curl_W^{(K)} F)_r(r,\theta,\varphi) = \frac{1}{r\sin\theta} D_{(K_\theta)}^{\theta,*}[\theta']\,(F_\varphi(r,\theta',\varphi)\sin\theta') - \frac{1}{r\sin\theta} D_{(K_\varphi)}^{\varphi,*}[\varphi']\,F_\varphi(r,\theta,\varphi'),$$

$$(Curl_W^{(K)} F)_\theta(r,\theta,\varphi) = \frac{1}{r\sin\theta} D_{(K_\varphi)}^{\varphi,*}[\varphi']\,F_r(r,\theta,\varphi') - \frac{1}{r} D_{(K_r)}^{r,*}[r']\,(r\,F_\varphi(r',\theta,\varphi)),$$

$$(Curl_W^{(K)} F)_\varphi(r,\theta,\varphi) = \frac{1}{r} D_{(K_r)}^{r,*}[r']\,(r\,F_\theta(r',\theta,\varphi)) - \frac{1}{r} D_{(K_\theta)}^{\theta,*}[\theta']\,F_r(r,\theta',\varphi).$$

### 11.2.3 General Fractional Grad, Div, Curl in Cylindrical Cooriditates

The Cartesian and cylindrical coordinates
$$\begin{aligned} x &= r\cos\varphi \\ y &= r\sin\varphi \\ z &= z, \end{aligned} \qquad (192)$$
where $r \in [0,\infty)$, $\varphi \in [0,2\pi]$, and $z \in [0,\infty)$. The basic vectors of these coordinate systems
$$\begin{aligned} e_x &= e_r\cos\varphi - e_\varphi\sin\varphi \\ e_y &= e_r\sin\varphi + e_\varphi\cos\varphi \\ e_z &= e_z. \end{aligned}$$
Note that it is impossible to obtain the fractional vector integral and differential operators in cylindrical coordinated by using coordinate transformation (192). Therefore, the definitions of general fractional vector operators should be formulated separately.

Let us define the regions

$$Cylinder_+^3 = \{(r,\varphi,z): \ r \in (0,\infty), \ \varphi \in (0,2\pi], \ z \in (0,infty)\},$$

$$Cylinder_{0,+}^3 = \{(r,\varphi,z): \ r \in [0,\infty), \ \varphi \in [0,2\pi], \ z \in [0,infty)\}.$$

We can see that
$$Cylinder_+^3 \subset \mathbb{R}^3+, \quad Cylinder_{0,+}^3 \subset \mathbb{R}_{0,+}^3.$$

Let us consider the vector field
$$F = e_r \, F_r(r,\varphi,z) + e_\varphi \, F_\varphi(r,\varphi,z) + e_z \, F_y(r,\varphi,z)$$

that belongs to $\mathbb{F}_{-1}^1(\mathbb{R}_+^3)$

The general frationa gradient, divergence and curl operators in cylindrical coordinates are defined in the following form. For simplicity, we give definitions of regional operators. Linear and surface general fractional operators in cylindrical coordinates are defined similarly.

**Definition 11.4** *If $F(r,\varphi,z) \in C_{-1}^1(\mathbb{R}_+^3)$, then the regional general fractional gradient in cylindrical coordinates for the region $W = Cylinder_{0,+}^3$ is defined as*
$$Grad_W^{(K)} F(r',\varphi,z) = e_r \, (Grad_W^{(K)} F)_r + e_\varphi \, (Grad_W^{(K)} F)_\varphi + e_z \, (Grad_W^{(K)} F)_z,$$

*where*
$$(Grad_W^{(K)} F)_r(r,\varphi,z) = D_{(K_r)}^{r,*}[r'] \, F(r',\varphi,z)$$

$$(Grad_W^{(K)} F)_\varphi(r,\varphi,z) = \frac{1}{r} D_{(K_\varphi)}^{\varphi,*}[\varphi'] \, F(r,\varphi',z),$$

$$(Grad_W^{(K)} F)_z(r,\varphi,z) = D_{(K_z)}^{z',*}[z'] \, F(r,\varphi,z').$$

**Definition 11.5** *If $F(r,\varphi,z) \in C_{-1}^1(\mathbb{R}_+^3)$, then the regional general fractional divergence in cylindrical coordinates for the region $W = Cylinder_{0,+}^3$ is defined as*
$$Div_W^{(K)} F = (Div_W^{(K)} F)_r + (Div_W^{(K)} F)_\varphi + (Div_W^{(K)} F)_z,$$

*where*
$$(Div_W^{(K)} F)_r(r,\varphi,z) = \frac{1}{r} D_{(K_r)}^{r,*}[r'] \, (r' \, F_r(r',\varphi,z))$$

$$(Div_W^{(K)} F)_\varphi(r,\varphi,z) = \frac{1}{r} D_{(K_\varphi)}^{\varphi,*}[\varphi'] \, F_\varphi(r,\varphi',z),$$

$$(Div_W^{(K)} F)_z(r,\varphi,z) = D_{(K_z)}^{z',*}[z'] \, F_z(r,\varphi,z').$$

We should note that
$$\frac{1}{r} D_{(K_r)}^{r,*}[r'] \, (r' \, F_r(r',\varphi,z)) \neq D_{(K_r)}^{r,*}[r'] \, F_r(r',\varphi,z) + \frac{1}{r} F_r(r,\varphi,z)$$

since the standard Leibniz rule does not hold for fractional derivatives of non-integer order and for general fractional derivatives.

**Definition 11.6** *If $F(r,\varphi,z) \in C_{-1}^1(\mathbb{R}_+^3)$, then the regional general fractional curl operator in cylindrical coordinates for the region $W = Cylinder_{0,+}^3$ is defined as*
$$Curl_W^{(K)} F = e_r \, (Curl_W^{(K)} F)_r + e_\varphi \, (Curl_W^{(K)} F)_\varphi + e_z \, (Curl_W^{(K)} F)_z,$$

*where*

$$(Curl_W^{(K)} F)_r(r, \varphi, z) = \frac{1}{r} D_{(K_\varphi)}^{\varphi,*}[\varphi'] F_z(r, \varphi', z) - D_{(K_z)}^{z,*}[z'] F_z(r, \varphi, z'),$$

$$(Curl_W^{(K)} F)_\varphi(r, \varphi, z) = D_{(K_z)}^{z,*}[z'] F_r(r, \varphi, z') - D_{(K_r)}^{r,*}[r'] F_z(r', \varphi, z),$$

$$(Curl_W^{(K)} F)_z(r, \varphi, z) = \frac{1}{r} D_{(K_r)}^{r,*}[r'] (r' F_\varphi(r', \varphi, z)) - \frac{1}{r} D_{(K_\varphi)}^{\varphi,*}[\varphi'] F_r(r, \varphi', z).$$

## 11.3 General Fractional Integral Operators in OCC

### 11.3.1 GFI in OCC

The following GGI operators can be used to define the line, surface and volume general fractional integrals in orthogonal curvilinear coordinates (the GF integrals in OCC).

**Definition 11.7** Let $f(q) = f(q_1, q_2, q_3)$ be function in $\mathbb{R}_{0,+}^3$ that satisfies the conditions

$$H_1(q_1, a, b) f(q_1, a, b) \in C_{-1}(0, \infty),$$

$$H_2(a, q_2, b) f(a, q_2, b) \in C_{-1}(0, \infty),$$

$$H_3(a, b, q_3) f(a, b, q_3) \in C_{-1}(0, \infty),$$

for all $a, b \geq 0$.
Then the set of such functions is denoted by $C_{H,-1}(\mathbb{R}_+^3)$.

Let function $f(q) = f(q_1, q_2, q_3)$ belong to the set $C_{H,-1}(\mathbb{R}_+^3)$. For the orthogonal curvilinear coordinate $q_k$ ($k = 1,2,3$), the GFI operators are defined as

$$\widehat{I_{(M_k)}^{q_k}}[q'_k] f(q) := I_{(M_k)}^{q_k}[q'_k]( H_k(q) f(q)) = \int_0^{q_k} M_k(q_k - q'_k) f(q) H_k(q) dq'_k, \quad (193)$$

and

$$\widehat{I^{(M_k)}}_{[a_k, b_k]}[q_k] f(q) := \widehat{I_{(M_k)}^{b_k}}[q_k] f(q) - \widehat{I_{(M_k)}^{a_k}}[q_k] f(q),$$

where notations (193) means the following

$$\widehat{I_{(M_1)}^{q_1}}[q'_1] f(q'_1, q_2, q_3) = I_{(M_1)}^{q_1}[q'_1]( H_1(q'_1, q_2, q_3) f(q'_1, q_2, q_3)) =$$

$$\int_0^{q_1} M_1(q_1 - q'_1) f(q'_1, q_2, q_3) H_1(q'_1, q_2, q_3) dq'_1,$$

$$\widehat{I_{(M_2)}^{q_2}}[q'_2] f(q_1, q'_2, q_3) := I_{(M_2)}^{q_2}[q'_2]( H_2(q_1, q'_2, q_3) f(q_1, q'_2, q_3)) =$$

$$\int_0^{q_2} M_2(q_2 - q'_2) f(q_1, q'_2, q_3) H_2(q_1, q'_2, q_3) dq'_2,$$

$$\widehat{I_{(M_3)}^{q_3}}[q'_3] f(q_1, q_2, q'_3) := I_{(M_3)}^{q_3}[q'_3]( H_3(q_1, q_2, q'_3) f(q_1, q_2, q'_3)) =$$

$$\int_0^{q_3} M_3(q_3 - q'_3) f(q_1, q_2, q'_3) H_3(q_1, q_2, q'_3) dq'_3.$$

### 11.3.2 GFI in Spherical Coordinates

To define the line spherical GFI, we can use the following GFI operators.
1) The radial GFI operator
$$\widehat{I^r_{(M_r)}}[r'] f(r, \theta', \varphi) := I^r_{(M_r)}[r'](r' f(r, \theta', \varphi)) = \int_0^r M_r(r - r') f(r, \theta', \varphi) r' \, dr'.$$
2) The polar GFI operator
$$\widehat{I^\theta_{M_\theta}}[\theta'] f(,\theta', \varphi) = I^\theta_{M_\theta}[\theta'] (r f(r, \theta', \varphi)) =$$

$$\int_0^\theta M_\theta(\theta - \theta') f(r, \theta', \varphi) r \, d\theta' = r \int_0^\theta M_\theta(\theta - \theta') f(r, \theta', \varphi) \, d\theta'.$$
3) The azimuthal GFI operator
$$\widehat{I^\varphi_{(M_\varphi)}}[\varphi'] f(r, \theta, \varphi') = I^\varphi_{(M_\varphi)}[\varphi'] (r \sin\theta \, f(r, \theta, \varphi')) =$$

$$\int_0^\varphi M_\varphi(\varphi - \varphi') f(r, \theta, \varphi') r \sin\theta \, d\varphi' = r \sin\theta \int_0^\varphi M_\varphi(\varphi - \varphi') f(r, \theta, \varphi') \, d\varphi'.$$

We can define the operators on the positive intervals
$$\widehat{I^{(M_r)}_{[r_0,r]}}[r'] f(r', \theta, \varphi) = \widehat{I^r_{(M_r)}}[r'] f(r', \theta, \varphi) - \widehat{I^{r_0}_{(M_r)}}[r'] f(r', \theta, \varphi)$$
for $r_0 > 0$. For $r_0 = r_1 > 0$, $r = r_1$, we have
$$\widehat{I^{(M_r)}_{[r_1,r_2]}}[r'] f(r', \theta, \varphi) = \widehat{I^{r_2}_{(M_r)}}[r'] f(r', \theta, \varphi) - \widehat{I^{r_1}_{(M_r)}}[r'] f(r', \theta, \varphi).$$
Integrals with respect to other variables are defined in a similar way.
$$\widehat{I^{M_\theta}_{[\theta_0,\theta]}}[\theta'] f(r, \theta', \varphi) = \widehat{I^\theta_{M_\theta}}[\theta'] f(r, \theta', \varphi) - \widehat{I^{\theta_0}_{M_\theta}}[\theta'] f(r, \theta', \varphi).$$

$$\widehat{I^{M_\theta}_{[\theta_1,\theta_2]}}[\theta'] f(r, \theta', \varphi) = \widehat{I^{\theta_2}_{M_\theta}}[\theta'] f(r, \theta', \varphi) - \widehat{I^{\theta_1}_{M_\theta}}[\theta'] f(r, \theta', \varphi),$$

and

$$\widehat{I^{(M_\varphi)}_{[\varphi_0,\varphi]}}[\varphi'] f(r, \theta, \varphi') = \widehat{I^\varphi_{(M_\varphi)}}[\varphi'] f(r, \theta, \varphi') - \widehat{I^{\varphi_0}_{(M_\varphi)}}[\varphi'] f(r, \theta, \varphi').$$

$$\widehat{I^{(M_\varphi)}_{[\varphi_1,\varphi_1]}}[\varphi'] f(r, \theta, \varphi') = \widehat{I^{\varphi_2}_{(M_\varphi)}}[\varphi'] f(r, \theta, \varphi') - \widehat{I^{\varphi_1}_{(M_\varphi)}}[\varphi'] f(r, \theta, \varphi').$$
For intervals $r_0$, $\theta_0 = 0$, and $\varphi = 0$, we define
$$\widehat{I^{(M_r)}_{[0,r]}}[r'] := \widehat{I^r_{(M_r)}}[r'], \quad \widehat{I^{(M_\theta)}_{[0,\theta]}}[\theta'] := \widehat{I^\theta_{(M_\theta)}}[\theta'], \quad \widehat{I^{(M_\varphi)}_{[0,\varphi]}}[\varphi'] := \widehat{I^\varphi_{(M_\varphi)}}[\varphi'].$$

### 11.3.3 GFI in Cylindrical Coordinates

Let us define the cylindrical GFI, we can use the following GFI operators.
1) If $f(r) \in C_{-1}(0, \infty)$, then
$$\widehat{I^r_{(M_r)}}[r'] f(r') := I^r_{(M_r)}[r'] f(r') := \int_0^r M_r(r - r') f(r') \, dr'$$
and
$$\widehat{I^{(M_r)}_{[r_1,r_2]}}[r] f(r') := \widehat{I^{r_2}_{(M_r)}}[r] f(r) - \widehat{I^{r_1}_{(M_r)}}[r] f(r) =$$

$$\int_0^{r_2} M_r(r_2 - r) f(r) \, dr - \int_0^{r_1} M_r(r_1 - r) f(r) \, dr,$$

where $r_2 > r_1 > 0$, and $\widehat{I_{[0,r_2]}^{(M_r)}}[r] := \widehat{I_{(M_r)}^{r_2}}[r]\,f(r)$.

2) If $f(\varphi) \in C_{-1}(0,\infty)$, then
$$\widehat{I_{(M_\varphi)}^{\varphi}}[\varphi']\,f(\varphi') := I_{(M_\varphi)}^{\varphi}[\varphi']\,(r\,f(\varphi')) :=$$

$$\int_0^\varphi M_\varphi(\varphi - \varphi')\,f(\varphi')\,r\,d\varphi' = r\int_0^\varphi M_\varphi(\varphi - \varphi')\,f(\varphi')\,d\varphi',$$

and
$$\widehat{I_{[\varphi_1,\varphi_2]}^{(M_\varphi)}}[\varphi]\,f(\varphi) := \widehat{I_{(M_\varphi)}^{\varphi_2}}[\varphi]\,f(\varphi) - \widehat{I_{(M_\varphi)}^{\varphi_1}}[\varphi]\,f(\varphi),$$

if $0 < \varphi_1 < \varphi_2 < 2\pi$ and
$$\widehat{I_{[0,\varphi_2]}^{(M_\varphi)}}[\varphi]\,f(\varphi) := \widehat{I_{(M_\varphi)}^{\varphi_2}}[\varphi]\,f(\varphi).$$

3) If $f(z) \in C_{-1}(0,\infty)$, then
$$\widehat{I_{(M_z)}^{z}}[z']\,f(z') = \int_0^z M_z(z - z')\,f(z')\,dz'.$$

and
$$\widehat{I_{[z_1,z_2]}^{M_z}}[z]\,f(z) = \widehat{I_{(M_z)}^{z_2}}[z]\,f(z) - \widehat{I_{(M_z)}^{z_1}}[z]\,f(z),$$

if $z_1 > 0$ and
$$\widehat{I_{[0,z_2]}^{(M_z)}}[z]\,f(z) = \widehat{I_{(M_z)}^{z_2}}[z]\,f(z).$$

## 11.4 General Fractional Operators in Curvilinear Coordinates

### 11.4.1 Definition of Line GFI for Vector Field in OCC

Let us define the line GFI in orthogonal curvilinear coordinates (OCC) in $\mathbb{R}^2_{0,+}$ of the $Q_1Q_2$-plane.

**Definition 11.8** Let $L$ be a simple line in $\mathbb{R}^2_{0,+}$ of the $Q_1Q_2$-plane. Let the functions
$$f_1(q_1) := H_1(q_1, y(q_1)) F_1(q_1, y(q_1)), \quad f_2(q_2) := H_2(x(q_2), q_2) F_2(x(q_2), q_2)$$
belong to the set $C_{-1}(0, \infty)$.
Then line GFI for the line $L$ is defined by the equation
$$(\hat{I}_L^{(M)}, F) = I_{[a,b]}^{(M_1)}[q_1] f_1(q_1) + I_{[c,d]}^{(M_2)}[q_2] f_2(q_2) =$$

$$I_{[a,b]}^{(M_1)}[q_1] H_1(q_1, y(q_1)) F_1(q_1, y(q_1)) + I_{[c,d]}^{(M_2)}[q_2] H_2(x(q_2), q_2) F_y(x(q_2), q_2). \quad (194)$$

Line GFI (194) exists, if the kernel pairs $(M_1(q_1), K_1(q_1))$ and $(M_2(q_2), K_2(q_2))$ belong to the Luchko set $\mathbb{L}$.

Let us define the line GFI in orthogonal curvilinear coordinates (OCC) in $\mathbb{R}^3_{0,+}$ of the $Q_1Q_2Q_3$-space. To give this definition, we will describe conditions on the vector field $F(q_1, q_2, q_3)$. We will assume that the vector field
$$F = F(q) = e_1 F_1(q_1, q_2, q_3) + e_2 F_2(q_1, q_2, q_3) + e_3 F_3(q_1, q_2, q_3)$$
on the simple line $L \subset \mathbb{R}^3_{0,+}$ is described by the following functions that belong to the space $C_{-1}(0, \infty)$:
$$f_1(q_1) := H_1(q_1, y(q_1), z(q_1)) F_1(q_1, y(q_1), z(q_1)) \in C_{-1}(0, \infty),$$

$$f_2(q_2) := H_2(x(q_2), q_2, z(q_2)) F_2(x(q_2), q_2, z(q_2)) \in C_{-1}(0, \infty),$$

$$f_3(q_3) := H_3(x(q_3), y(q_3), q_3) F_3(x(q_3), y(q_3), q_3) \in C_{-1}(0, \infty).$$
If these conditions are satisfied, then we will write $F(q) \in \mathbb{F}_{H,-1}(L)$.
In the case $F \in \mathbb{F}_{H,-1}(\mathbb{R}^3_+)$, the line general fractional integral for the vector field $F$ and the line $L = AB$ with endpoints $A(a_1, a_2, a_3)$ and $B(b_1, b_1, b_3)$ with all $b_k > a_k \geq 0$ is defined by the equation
$$(\hat{I}_L^{(M)}, F) = \sum_{k=1}^{3} I_{[a_k, b_k]}^{(M_k)}[q_k] f_k(q_k) =$$

$$I_{[a_1, b_1]}^{(M_1)}[q_1] f_1(q_1) + I_{[a_2, b_2]}^{(M_2)}[q_2] f_2(q_2) + I_{[a_3, b_3]}^{(M_3)}[q_3] f_3(q_3) =$$

$$I_{[a_1, b_1]}^{(M_1)}[q_1] H_1(q_1, y(q_1), z(q_1)) F_1(q_1, y(q_1), z(q_1)) +$$

$$I_{[a_2, b_2]}^{(M_2)}[q_2] H_2(x(q_3), q_2, z(q_2)) F_2(x(q_3), q_2, z(q_2)) +$$

$$I_{[a_3, b_3]}^{(M_3)}[q_3] H_3(x(q_3), y(q_3), q_3) F_3(x(q_3), y(q_3), q_3). \quad (195)$$

Using the Laplace convolution, the line GFI can be written as
$$(\hat{I}_L^{(M)}, F) = \sum_{k=1}^{3} ((M_k * f_k)(b_k) - (M_k * f_k)(a_k)),$$
where $b_k > a_k \geq 0$, and
$$(M_k * f_k)(c_k) = \int_0^{c_k} dq_k\, M(c_k - q_k)\, f_k(q_k),$$
where $c_k = a_k$ or $c_k = b_k$ with $M_k(q_k) \in C_{-1,0}(0,\infty)$ and $f_k(q_k) \in C_{-1}(0,\infty)$.

### 11.4.2 Line GFI for Piecewise Simple Lines in OCC

Let us consider a line $L \subset \mathbb{R}_{0,+}^3$, which can be divided into several lines $L_k = L_k[A_k, A_{k+1}]$, $k = 1, \dots, n$ that are simple lines or lines parallel to one of the axes:
$$L := \bigcup_{k=1}^{n-1} L_k[A_k, A_{k+1}], \tag{196}$$
where the line $L_k$ connects the points $A_k(x_k, y_k, z_k)$, and $A_{k+1}(x_{k+1}, y_{k+1}, z_{k+1})$ with $x_k, y_k, z_k \geq 0$ for all $k = 1, 2, \dots, n-1$. Lines of this kind will be called the piecewise simple lines.

For piecewise simple line (196) in OCC of $\mathbb{R}_{0,+}^3$, and the vector field $F \in \mathbb{F}_{-1}(L)$, the line GFI is defined by the equation
$$(\hat{I}_L^{(M)}[q_1, q_2, q_3], F(q_1, q_2, q_3)) := \sum_{k=1}^{n-1} (I_{[A_k, A_{k+1}]}^{(M)}[q_1, q_2, q_3], F(q_1, q_2, q_3)),$$
where
$$(\hat{I}_{[A_k, A_{k+1}]}^{(M)}[q_1, q_2, q_3], F(q_1, q_2, q_3)) :=$$

$$sgn(q_{1,k+1} - q_{1,k})\, I_{\omega[q_{1,k}, q_{1,k+1}]}^{(M_1)}[q_1]\, H_1(q_1, y_k(q_1), z_k(q_1))\, F_1(q_1, y_k(q_1), z_k(q_1)) +$$

$$sgn(q_{2,k+1} - q_{2,k})\, I_{\omega[q_{2,k}, q_{2,k+1}]}^{(M_2)}[q_2]\, H_2(x_k(q_2), q_2, z_k(q_2))\, F_2(x_k(q_2), q_2, z_k(q_2)) +$$

$$sgn(q_{3,k+1} - q_{3,k})\, I_{\omega[q_{3,k}, q_{3,k+1}]}^{(M_3)}[q_3]\, H_3(x_k(q_3), y_k(q_3), q_3)\, F_3(x_k(q_3), y_k(q_3), q_3),$$
where the functions $y = y_k(q_1)$, $x = x_k(q_2)$, $z = z_k(q_1)$, $x = x_k(q_3)$, $z = z_k(q_2)$, $y = y_k(q_1)$ define the lines $L_k$ that are simple or parallel to one of the axes of OCC.

### 11.4.3 Regional GF Gradient in OCC

Let us give definitions of a set of scalar fields and a general fractional gradient in OCC for $\mathbb{R}_{0,+}^3$.

**Definition 11.9** *Let $U(q) = U(q_1, q_2, q_3)$ be a scalar field that satisfies the conditions*
$$D_{(K_1)}^{q_1,*}[q'_1]U(q'_1, q_2, q_3) \in C_{-1}(\mathbb{R}_+^3),$$

$$D_{(K_2)}^{q_2,*}[q'_2]U(q_1, q'_2, q_3) \in C_{-1}(\mathbb{R}_+^3),$$

$$D_{(K_3)}^{q_3,*}[q'_3]U(q_1, q_2, q'_3) \in C_{-1}(\mathbb{R}_+^3).$$
*Then the set of such scalar fields $U(q) = U(q_1, q_2, q_3)$ will be denoted as $C_{-1}^1(\mathbb{R}_+^3)$.*

**Definition 11.10** *Let* $U(q) = U(q_1, q_2, q_3)$ *be a scalar field that belongs to the set* $C_{-1}^1(\mathbb{R}_+^3)$.

*Then the general fractional gradient* $Grad_W^{(K)}$ *in OCC for the region* $W = \mathbb{R}_{0,+}^3$ *is defined as*

$$(Grad_W^{(K)} U)(x, y, z) = e_1 \frac{1}{H_1(q)} D_{(K_1)}^{q_1,*}[x'] U(q'_1, q_2, q_3) +$$

$$e_2 \frac{1}{H_2(q)} D_{(K_2)}^{q_2,*}[y'] U(q_1, q'_2, q_3) + e_3 \frac{1}{H_3(q)} D_{(K_3)}^{q_3,*}[z'] U(q_1, q_2, q'_3).$$

*This operator will be called the regional general fractional gradient (regional GF gradient) in OCC.*

### 11.4.4 Line GF Gradient in OCC for $\mathbb{R}_{0,+}^2$

The general fractional gradients (Line GF Gradient) in OCC can be defined not only for the region $W = \mathbb{R}_{0,+}^3$. Using the fact that GFD is integro-differential operator, we can define GF Gradients for line $L \subset \mathbb{R}_{0,+}^2$ and $L \subset \mathbb{R}_{0,+}^3$. The general fractional vector operators can be defined as a sequential action of first-order derivatives and general fractional integrals.

**Definition 11.11** *Let* $L$ *be a simple line in* $\mathbb{R}_{0,+}^2$ *of the* $Q_1 Q_2$-*plane, and* $U(q_1, q_2) \in \mathbb{F}_{-1,L}^1(\mathbb{R}_+^2)$ *that means*

$$U_{q_1}^{(1)}(q_1, y(q_1)) \in C_{-1}(0, \infty), \quad U_{q_2}^{(1)}(x(q_2), q_1) \in C_{-1}(0, \infty).$$

*Then, the line general fractional gradient in OCC (Line GF Gradient in OCC) for the line* $L$ *is defined by the equation*

$$Grad_L^{(K)} U(q_1, q_2) :=$$

$$e_1 \frac{1}{H_1(q)} I_{(K_1)}^{q_1,*}[q'_1] U_{q'_1}^{(1)}(q'_1, y(q'_1)) + e_2 \frac{1}{H_2(q)} I_{(K_2)}^{q_2,*}[q'_2] U_{q'_2}^{(1)}(x(q'_2), q'_1) =$$

$$e_1 \frac{1}{H_1(q_1, q_2)} I_{(K_1)}^{q_1,*}[q'_1] (\frac{\partial U(q'_1, q_2)}{\partial q'_1})_{q_2 = y(q'_1)} + e_2 \frac{1}{H_2(q_1, q_2)} I_{(K_2)}^{q_1,*}[q'_2] (\frac{\partial U(q_1, q'_2)}{\partial q'_2})_{q_1 = x(q'_2)},$$

*where the kernels* $(M_1(q_1), K_1(q_1))$ *and* $(M_2(q_2), K_2(q_1))$ *belong to the Luchko set* $\mathbb{L}$.

We can use the definitions of the line GFI and line GF Gradient in OCC to prove the following theorem.

**Theorem 11.1** *(General Fractional Gradient Theorem for line GF Gradient for* $\mathbb{R}_{0,+}^2$*)*
*Let* $L$ *be a simple line in* $\mathbb{R}_{0,+}^2$ *of the* $Q_1 Q_2$-*plane, which connects the points* $A(a, b)$ *and* $B(c, d)$, *and the scalar field* $U(q_1, q_2)$ *belongs to the set* $\mathbb{F}_{-1,L}^1(\mathbb{R}_+^2)$.
*Then, the equality*

$$(I_L^{(M)}, Grad_L^{(K)} U) = U(c, d) - U(a, b)$$

*holds, where* $b \geq a \geq 0$ *and* $d \geq c \geq 0$.

### 11.4.5 Theorem for Line GF Gradient in OCC of $\mathbb{R}^3_{0,+}$

In the theorem we will use the following definitions.

**Definition 11.12** Let $L \subset \mathbb{R}^3_{0,+}$ be a line in OCC that is described by the functions
$$q_2 = y(q_1) \geq 0, \quad q_3 = z(q_1) \geq 0, \quad q_1 \in [a,b] \subset \mathbb{R}_{0,+}, \tag{197}$$
which are continuously differentiable functions for $q_1 \in (a,b) \subset \mathbb{R}_{0,+}$, i.e. $y(q_1), z(q_1) \in C^1(a,b)$. Then this line will be called $Q_2Q_3$-simple line. The line $L$ is called simple line in OCC of $\mathbb{R}^3_{0,+}$, if $L$ is $Q_1Q_2$-, $Q_1Q_3$- and $Q_2Q_3$-simple line.

Let us consider a simple line $L \subset \mathbb{R}^3_{0,+}$ in OCC, which connects the points $A(a,c,e)$ and $B(b,d,f)$, and can be described by following equivalent forms
$$L = \{(q_1,q_2,q_3): \ q_1 \in [a,b], \ q_2 = y(q_1) \in C^1[a,b], \ q_3 = z(q_1) \in C^1[a,b]\}, \tag{198}$$

$$L = \{(q_1,q_2,q_3): \ q_2 \in [c,d], \ q_1 = x(q_2) \in C^1[c,d], \ q_3 = z(q_2) \in C^1[c,d]\}, \tag{199}$$

$$L = \{(q_1,q_2,q_3): \ q_3 \in [e,f], \ q_1 = x(Q_3) \in C^1[e,f], \ q_2 = y(q_3) \in C^1[e,f]\}, \tag{200}$$
where $y(q_1 = a) = c$, $y(1_1 = b) = d$, $z(q_1 = a) = e$, $z(q_1 = b) = f$.

**Definition 11.13** Let $L$ be a simple line in $\mathbb{R}^3_{0,+}$, which is defined in form (198), (199), (200), and a scalar field $U(q_1,q_2,q_3)$ satisfies the conditions
$$U_1(x) := U^{(1)}_x(x,y(x),z(x)) = \left(\frac{\partial U(q_1,q_2,q_3)}{\partial x}\right)_{y=y(x),z=z(x)} \in C^n_{-1}(0,\infty),$$

$$U_2(y) := U^{(1)}_y(x(y),y,z(y)) = \left(\frac{\partial U(q_1,q_2,q_3)}{\partial y}\right)_{x=x(y),z=z(y)} \in C^n_{-1}(0,\infty),$$

$$U_3(z) := U^{(1)}_z(x(z),y(z),z) = \left(\frac{\partial U(q_1,q_2,q_3)}{\partial z}\right)_{x=x(z),y=y(z)} \in C^n_{-1}(0,\infty).$$
Then the set of such scalar fields $U(q_1,q_2,q_3)$ will be denoted as $U(q_1,q_2,q_3) \in \mathbb{F}^n_{-1}(L)$.

**Definition 11.14** Let $L$ be a simple line in $\mathbb{R}^3_{0,+}$, which is defined in form (198), (199), (200), and a scalar field $U(q_1,q_2,q_3)$ belongs to the set $\mathbb{F}^1_{-1}(L)$.
Then the line general fractional gradient for the line $L \in \mathbb{R}^3_{0,+}$ is defined by the equation
$$(Grad^{(K)}_L U)(q_1,q_2,q_3) = (D^{(M)}_L U)(q_1,q_2,q_3) =$$
$$\frac{1}{H_1(q)} I^{(M_1)}_{[a,b]}[q_1] U^{(1)}_1(q_1,y(q_1),z(q_1)) +$$

$$\frac{1}{H_2(q)} I^{(M_2)}_{[c,d]}[q_2] U^{(1)}_2(x(q_2),q_2,z(q_2)) + \frac{1}{H_3(q)} I^{(M_3)}_{[e,f]}[q_3] U^{(1)}_3(x(q_3),y(q_3),q_3). \tag{201}$$
where the pairs of the kernels $(M_k(q_k), K_{(q_k)})$, $(k=1,2,3)$ belong to the Luchko set $\mathbb{L}$.

For OCC, the general fractional gradient theorem for the line GF Gradient is fomulated for simple lines in the form.

**Theorem 11.2** *(General Fractional Gradient Theorem for OCC and Line GF Gradient)*

Let $L$ be a simple line in OCC for $\mathbb{R}^3_{0,+}$, which is described in form (198), (199), (200), and connects the points $A(a,c,e)$ and $B(b,d,f)$.

Let $U(q_1, q_2, q_3)$ be a scalar filed that belongs to the set $\mathbb{F}^1_{-1,L}(\mathbb{R}^3_+)$.

Then, the equality
$$(\hat{I}^{(M)}_L, Grad^{(K)}_L U) = U(a,c,e) - U(b,d,f)$$
holds.

*Proof.* Using that the line GFI of a vector field $F$ for the simple line $L \subset \mathbb{R}^3_{0,+}$ is defined by the equation
$$(\hat{I}^{(M)}_L, F) = I^{(M_1)}_{[a,b]}[q_1]\, H_1(q_1, y(q_1), z(q_1))\, F_1(q_1, y(q_1), z(q_1)) +$$
$$I^{(M_2)}_{[c,d]}[q_2]\, H_2(x(q_2), q_2, z(q_2))\, F_2(x(q_2), q_2, z(q_2)) +$$
$$I^{(M_3)}_{[e,f]}[q_3]\, H_3(x(q_3), y(q_3), q_3)\, F_3(x(q_3), y(q_3), q_3),$$
we can consider the line GFI of the vector field that is defined by the line general fractional gradient
$$F_1(q_1, y(q_1), z(q_1)) := (Grad^{(K)}_L U(q_1, q_2, q_3))_1 =$$
$$\frac{1}{H_1(q_1, y(q_1), z(q_1))} I^{q_1,*}_{(K_1)}[q'_1]\, U^{(1)}_{q'_1}(q'_1, y(q'_1), z(q'_1)),$$

$$F_2(x(q_2), q_2, z(q_2)) := (Grad^{(K)}_L U(q_1, q_2, q_3))_2 =$$
$$\frac{1}{H_2(x(q_2), q_2, z(q_2))} I^{q_2,*}_{(K_2)}[q'_2]\, U^{(1)}_{q'_2}(x(q'_2), q'_2, z(q'_2)),$$

$$F_3(x(q_3), y(q_3), q_3) := (Grad^{(K)}_L U(q_1, q_2, q_3))_3 =$$
$$\frac{1}{H_3(x(q_3), y(q_3), q_3)} I^{q_3,*}_{(K_3)}[q'_3]\, U^{(1)}_{q'_3}(x(q'_3), y(q'_3), q'_3).$$

Therefore, we can get the line GFI of the line FG Gradient
$$(\hat{I}^{(M)}_L, Grad^{(K)}_L U) =$$

$$\widehat{I^{(M_1)}_{[a,b]}}[q_1] \frac{1}{H_1(q_1, y(q_1), z(q_1))} I^{q_1,*}_{(K_1)}[q'_1]\, U^{(1)}_{q'_1}(q'_1, y(q'_1), z(q'_1)) +$$

$$\widehat{I^{(M_2)}_{[c,d]}}[q_2] \frac{1}{H_2(x(q_2), q_2, z(q_2))} I^{q_2,*}_{(K_2)}[q'_2]\, U^{(1)}_{q'_2}(x(q'_2), q'_2, z(q'_2)) +$$

$$\widehat{I^{(M_3)}_{[e,f]}}[q_3] \frac{1}{H_3(x(q_3), y(q_3), q_3)} I^{q_3,*}_{(K_3)}[q'_3]\, U^{(1)}_{q'_3}(x(q'_3), y(q'_3), q'_3) =$$

$$I^{(M_1)}_{[a,b]}[q_1]\, I^{q_1,*}_{(K_1)}[q'_1]\, U^{(1)}_{q'_1}(q'_1, y(q'_1), z(q'_1)) +$$

$$I_{[c,d]}^{(M_2)}[q_2]\, I_{(K_2)}^{q_2,*}[q'_2]\, U_{q'_2}^{(1)}(x(q'_2), q'_2, z(q'_2)) +$$

$$I_{[e,f]}^{(M_3)}[q_3]\, I_{(K_3)}^{q_3,*}[q'_3]\, U_{q'_3}^{(1)}(x(q'_3), y(q'_3), q'_3),$$

where we use

$$\widehat{I_{[a,b]}^{(M_1)}}[q_1] = I_{[a,b]}^{(M_1)}[q_1]\, H_1(q_1, y(q_1), z(q_1)),$$

$$\widehat{I_{[c,d]}^{(M_2)}}[q_2] = I_{[c,d]}^{(M_2)}[q_2]\, H_2(x(q_2), q_2, z(q_2)),$$

$$\widehat{I_{[e,f]}^{(M_3)}}[q_3] = I_{[e,f]}^{(M_3)}[q_3]\, H_3(x(q_3), y(q_3), q_3).$$

Further, transformations are made similar to the transformations performed in the proof of Theorem 11.2. Using the fact that the kernels $(M_k(q_k), K_k(q_k))$, $k=1,2,3$, belong to the Luchko set $\mathbb{L}$, we obtain

$$I_{[a,b]}^{(M_1)}[q_1]\, I_{(K_1)}^{q_1,*}[q'_1]\, f_1(q'_1) = (M_1 * (K_1 * f_1))(b) - (M_1 * (K_1 * f_1))(a) =$$

$$((M_1 * K_1) * f_1)(b) - ((M_1 * K_1) * f_1)(a) = (\{1\} * f_1)(b) - (\{1\} * f_1)(a) =$$

$$\int_0^b f_1(q_1)\, dq_1 - \int_0^a f_1(q_1)\, dq_1 = \int_a^b f_1(q_1)\, dq_1.$$

Similarly, we can obtain the equations for $q_2$ and $q_3$.

Therefore, we get

$$(\hat{I}_L^{(M)}, Grad_L^{(K)} U) =$$

$$\int_a^b \left(\frac{\partial U(q'_1, q_2, q_3)}{\partial q'_1}\right)_{q_2=y(q'_1), q_3=z(q'_1)} dq'_1 + \int_c^d \left(\frac{\partial U(q_1, q'_2, q_3)}{\partial q'_2}\right)_{q_1=x(q'_2), q_3=z(q'_2)} dq'_2 +$$

$$\int_e^f \left(\frac{\partial U(q_1, q_2, q'_3)}{\partial q'_3}\right)_{q_1=x(q'_3), q_2=y(q'_3)} dq'_3 =$$

$$\int_a^b \left(\frac{\partial U(q'_1, q_2, q_3)}{\partial q'_1}\right)_{q_2=y(q'_1), q_3=z(q'_1)} dq'_1 + \int_a^b \left(\frac{\partial U(q_1, q'_2, q_3)}{\partial q'_2}\right)_{q_2=y(q'_1), q_3=z(q'_1)} \frac{\partial y(q'_1)}{\partial q'_1} dq'_1 +$$

$$\int_z^b \left(\frac{\partial U(q_1, q_2, q'_3)}{\partial q'_3}\right)_{q_2=y(q'_1), q_3=z(q'_1)} \frac{\partial z(q'_1)}{\partial q'_1} dq'_1 =$$

$$\int_a^b \frac{dU(q_1, y(q_1), z(q_1))}{dq_1} dq_1 = U(b, d, f) - U(a, c, e),$$

where the standard gradient theorem is used.

As a result, we proved the general fractional gradient theorem for OCC with the line GF Gradient and the simple line in $\mathbb{R}_{0,+}^3$.

## 11.5 Regional and Surface GF Curl in OCC

In this section, we proposed definitions of two type GF Curl operators in orthogonal curvilinear coordinates (OCC) for $\mathbb{R}_{0,+}^3$.

### 11.5.1 Regional GF Curl in OCC

Let us give the definition of Regional GF Curl operator in OCC for $\mathbb{R}^3_{0,+}$.

**Definition 11.15** Let $F(q_1, q_2, q_3)$ be a vector field which satisfies conditions
$$\hat{F}_k(q_1, q_2, q_3) := H_k(q_1, q_2, q_3) F_k(q_1, q_2, q_3) \in C^1_{-1}(\mathbb{R}^3_+) \qquad (202)$$
for all $k = 1, 2, 3$ in OCC.
Then the regional general fractional curl in OCC for the region $W = \mathbb{R}^3_{0,+}$ is defined as
$$Curl_W^{(K)} F = e_1 (Curl_W^{(K)} F)_1 + e_2 (Curl_W^{(K)} F)_2 + e_3 (Curl_W^{(K)} F)_3,$$
where
$$(Curl_W^{(K)} F)_1(q_1, q_2, q_3) = \frac{1}{H_2 H_3} \left( D^{q_2,*}_{(K_2)}[q'_2] \hat{F}_3(q_1, q'_2, q_3) - D^{q_3,*}_{(K_3)}[q'_3] \hat{F}_2(q_1, q_2, q'_3) \right),$$

$$(Curl_W^{(K)} F)_2(q_1, q_2, q_3) = \frac{1}{H_1 H_3} \left( D^{q_3,*}_{(K_3)}[q'_3] \hat{F}_1(q_1, q_2, q'_3) - D^{q_1,*}_{(K_1)}[q'_1] \hat{F}_3(q'_1, q_2, q_3) \right),$$

$$(Curl_W^{(K)} F)_3(q_1, q_2, q_3) = \frac{1}{H_1 H_2} \left( D^{q_1,*}_{(K_1)}[q'_1] \hat{F}_2(q'_1, q_2, q_3) - D^{q_2,*}_{(K_2)}[q'_2] \hat{F}_1(q_1, q'_2, q_3) \right)$$
with $q_1, q_2, q_3 \geq 0$.

### 11.5.2 Surface GF Integral in OCC

Using Definition 7.2 piecewise simple surface ($S \in \mathbb{P}(\mathbb{R}^3_{0,+})$) in OCC, we propose definition of surface GF integral in OCC.
Let us define a set $\mathbb{F}_{S,H}(\mathbb{R}^3_{0,+})$. of vector fields that is used in definition of surface GF integral in OCC.

**Definition 11.16** Let $S$ be a piecewise simple surface ($S \in \mathbb{P}(\mathbb{R}^3_{0,+})$). Let the vector field
$$F := e_1 F_1(q_1, q_2, q_3) + e_2 F_2(q_1, q_2, q_3) + e_3 F_3(q_1, q_2, q_3) \qquad (203)$$
on the surface $S$ satisfy the conditions
$$\tilde{F}_1(x_i(q_2, q_3), q_2, q_3) \in C_{-1}(\mathbb{R}^2_+), \qquad (204)$$

$$\tilde{F}_2(q_1, y_j(q_1, q_2), q_3) \in C_{-1}(\mathbb{R}^2_+), \qquad (205)$$

$$\tilde{F}_3(q_1, q_2, z_k(q_2, q_3)) \in C_{-1}(\mathbb{R}^2_+), \qquad (206)$$
for all $i = 1, \ldots, n_x$, $j = 1, \ldots, n_y$, and $k = 1, \ldots, n_z$, where we use the notations
$$\tilde{F}_1(q_1, q_2, q_3) := H_2(q_1, q_2, q_3) H_3(q_1, q_2, q_3) F_1(q_1, q_2, q_3),$$

$$\tilde{F}_2(q_1, q_2, q_3) := H_1(q_1, q_2, q_3) H_3(q_1, q_2, q_3) F_2(q_1, q_2, q_3),$$

$$\tilde{F}_3(q_1, q_2, q_3) := H_1(q_1, q_2, q_3) H_2(q_1, q_2, q_3) F_3(q_1, q_2, q_3).$$
The set of such vector fields $F$ on piecewise simple surface $S$ will be denoted by $\mathbb{F}_{S,H}(\mathbb{R}^3_{0,+})$.

We can use the representationof the tilde $F_k$ in the form

$$\tilde{F}_k(q_1,q_2,q_3) := H_1(q_1,q_2,q_3)\, H_2(q_1,q_2,q_3)\, H_3(q_1,q_2,q_3)\, \frac{F_k(q_1,q_2,q_3)}{H_k(q_1,q_2,q_3)}. \tag{207}$$

Let us give a definition of surface GF integral in OCC.

**Definition 11.17** Let $S$ be a piecewise simple surface ($S \in \mathbb{P}(\mathbb{R}^3_{0,+})$) and a vector fields $F$ on this surface $S$ belongs to the set $\mathbb{F}_{S,H}(\mathbb{R}^3_{0,+})$.

Then, the surface general fractional vector integral in OCC (Surface GFI in OCC) of the second kind

$$\hat{I}^{(M)}_S = e_1\, \hat{I}^{(M)}_{S_{23}}[q_2,q_3] + e_2\, \hat{I}^{(M)}_{S_{13}}[q_1,q_3] + e_3\, \hat{I}^{(M)}_{S_{12}}[q_1,q_2] =$$

$$e_1 \sum_{i=1}^{n_1} \hat{I}^{(M)}_{S_{i,23}}[q_2,q_3] + e_1 \sum_{j=1}^{n_2} \hat{I}^{(M)}_{S_{j,13}}[q_1,q_3] + e_2 \sum_{k=1}^{n_3} \hat{I}^{(M)}_{S_{k,12}}[q_1,q_2] \tag{208}$$

for the vector fiels $F \in \mathbb{F}_{S,H}(\mathbb{R}^3_{0,+})$ is defined by the equation

$$(\hat{I}^{(M)}_S, F) := \sum_{i=1}^{n_1} \hat{I}^{(M)}_{S_{i,23}}[q_2,q_3]\, F_1(x_i(y,z),y,z) +$$

$$\sum_{j=1}^{n_2} \hat{I}^{(M)}_{S_{j,13}}[q_1,q_3]\, F_2(q_1, y_j(q_1,q_3), q_3) + \sum_{k=1}^{n_3} \hat{I}^{(M)}_{S_{k,12}}[q_1,q_2]\, F_3(q_1,q_2, z_k(q_1,q_2)) =$$

$$\sum_{i=1}^{n_1} I^{(M)}_{S_{i,23}}[q_2,q_3]\, \tilde{F}_1(x_i(q_2,q_3), q_2, q_3) + \sum_{j=1}^{n_2} I^{(M)}_{S_{j,13}}[q_1,q_3]\, \tilde{F}_2(q_1, y_j(q_1,q_3), q_3) +$$

$$\sum_{k=1}^{n_3} I^{(M)}_{S_{k,12}}[q_1,q_2]\, \tilde{F}_3(q_1,q_2, z_k(q_1,q_2)).$$

Here the areas

$$S_{23} = \bigcup_{i=1}^{n_1} S_{i,23}, \quad S_{13} = \bigcup_{j=1}^{n_2} S_{j,13}, \quad S_{12} = \bigcup_{k=1}^{n_3} S_{k,12}$$

are the projections of the surface $S$ onto the $Q_2 Q_3$, $Q_1 Q_3$, $Q_1 Q_2$ planes in OCC, where $S_{i,23}$, $S_{j,13}$ and $S_{k,12}$ are simple areas in these planes (or simple along some axes).

### 11.5.3 Surface GF Curl in OCC

Let us now give a definition of the Surface GF Curl in OCC for $\mathbb{R}^3_{0,+}$ for piecewise simple surface. The piecewise simple surface is given in Definition in 7.2.

Let us give a definition of the set $\mathbb{F}^1_{S,H}(\mathbb{R}^3_{0,+})$ of vector fields $F$ that are used to define the surface GF Curl operator in OCC (compare with Definition 7.3 of the set $\mathbb{F}^1_S(\mathbb{R}^3_{0,+})$ that is used in definition of the surface GFI).

**Definition 11.18** Let $S$ be a piecewise simple surface ($S \in \mathbb{P}(\mathbb{R}^3_{0,+})$).

Let the vector field

$$F := e_1\, F_1(q_1,q_2,q_3) + e_2\, F_2(q_1,q_2,q_3) + e_3\, F_3(q_1,q_2,q_3) \tag{209}$$

on the surface $S$ satisfy the conditions

$$\left(D^1_{q_3}\hat{F}_1\right)(q_1, y_j(q_1,q_3), q_3) := \left(\frac{\partial \hat{F}_1(q_1,q_2,q_3)}{\partial q_3}\right)_{q_2 = y_j(q_1,q_3)} \in C_{-1}(\mathbb{R}^2_+),$$

$$\left(D^1_{q_2}\hat{F}_1\right)(q_1, q_2, z_k(q_1,q_2)) := \left(\frac{\partial \hat{F}_1(q_1,q_2,q_3)}{\partial q_2}\right)_{q_3 = z_k(q_1,q_2)} \in C_{-1}(\mathbb{R}^2_+),$$

$$(D_{q_1}^1 \hat{F}_2)(q_1, q_2, z_k(q_1, q_2)) := (\frac{\partial \hat{F}_2(q_1, q_2, q_3)}{\partial q_1})_{q_3 = z_k(q_1, q_2)} \in C_{-1}(\mathbb{R}_+^2),$$

$$(D_{q_3}^1 \hat{F}_2)(x_i(q_2, q_3), q_2, q_3) := (\frac{\partial \hat{F}_2(q_1, q_2, q_3)}{\partial q_3})_{q_1 = x_i(q_2, q_3)} \in C_{-1}(\mathbb{R}_+^2),$$

$$(D_{q_1}^1 \hat{F}_3)(q_1, y_j(q_1, q_3), q_3) := (\frac{\partial \hat{F}_3(q_1, q_2, q_3)}{\partial q_1})_{q_2 = y_j(q_1, q_3)} \in C_{-1}(\mathbb{R}_+^2),$$

$$(D_{q_2}^1 \hat{F}_3)(x_i(q_2, q_3), q_2, q_3) := (\frac{\partial \hat{F}_3(q_1, q_2, q_3)}{\partial q_2})_{q_1 = x_i(q_2, q_3)} \in C_{-1}(\mathbb{R}_+^2)$$

for all $i = 1, \ldots, n_1$, $j = 1, \ldots, n_2$, and $k = 1, \ldots, n_3$.

The set of such vector fields $F$ on piecewise simple surface $S$ will be denoted by $\mathbb{F}_{S,H}^1(\mathbb{R}_{0,+}^3)$.

Then the surface general fractional curl is defined in the following form.

**Definition 11.19** Let $S$ be a piecewise simple surface ($S \in \mathbb{P}(\mathbb{R}_{0,+}^3)$) and a vector fields $F$ on this surface $S$ belongs to the set $\mathbb{F}_S^1(\mathbb{R}_{0,+}^3)$.

Then the surface general fractional curl in OCC (Surface GF Curl in OCC) of the vector field $F(q_1, q_2, q_3)$ for the piecewise simple surface $S$ is defined as

$$Curl_S^{(K)} F = e_1 (Curl_S^{(K)} F)_1 + e_2 (Curl_S^{(K)} F)_2 + e_3 (Curl_S^{(K)} F)_3,$$

where

$$(Curl_S^{(K)} F)_1 = \sum_{i=1}^{n_1} \frac{1}{H_2(x_i(q_2,q_3),q_2,q_3) \, H_3(x_i(q_2,q_3),q_2,q_3)}$$

$$(I_{(K_2)}^{q_2,*}[q'_2] (D_{q'_2}^1 \hat{F}_3)(x_i(q'_2, q_3), q'_2, q_3) - I_{(K_3)}^{q_3,*}[q'_3] (D_{q'_3}^1 \hat{F}_2)(x_i(q_2, q'_3), q_2, q'_3)),$$

$$(Curl_S^{(K)} F)_2 = \sum_{j=1}^{n_2} \frac{1}{H_1(q_1,y_j(q_1,q_3),q_3) \, H_3(q_1,y_j(q_1,q_3),q_3)}$$

$$(I_{(K_3)}^{q_3,*}[q'_3] (D_{q'_3}^1 \hat{F}_1)(q_1, y_j(q_1, q'_3), q'_3) - I_{(K_1)}^{q_1,*}[q'_1] (D_{q'_1}^1 \hat{F}_3)(q'_1, y_j(q'_1, q_3), q_3)),$$

$$(Curl_S^{(K)} F)_3 = \sum_{k=1}^{n_2} \frac{1}{H_1(q_1,q_2,z_k(q_1,q_2)) \, H_2(q_1,q_2,z_k(q_1,q_2))}$$

$$(I_{(K_1)}^{q_1,*}[q'_1] (D_{q'_1}^1 \hat{F}_2)(q'_1, q_2, z_k(q'_1, q_2)) - I_{(K_2)}^{q_2,*}[q'_2] (D_{q'_2}^1 \hat{F}_1)(q_1, q'_2, z_k(q_1, q'_2)))$$

with $q_1, q_2, q_3 \geq 0$.

Let us give the expression of the surface GFI of the vector field in the form of the Surface GF Curl in OCC.

**Remark 11.1** The surface GFI of the surface GF Curl in OCC has the form

$$(\hat{I}_S^{(M)}, Curl_S^{(K)} F) := \sum_{i=1}^{n_1} \hat{I}_{S_{i,23}}^{(M)}[q_2, q_3] (Curl_S^{(K)} F)_1 (x_i(y, z), y, z) +$$

$$\sum_{j=1}^{n_2} \hat{I}_{S_{j,13}}^{(M)}[q_1,q_3]\, (Curl_S^{(K)}F)_2(q_1,y_j(q_1,q_3),q_3) +$$

$$\sum_{k=1}^{n_3} \hat{I}_{S_{k,12}}^{(M)}[q_1,q_2]\, (Curl_S^{(K)}F)_3(q_1,q_2,z_k(q_1,q_2)) =$$

$$\sum_{i=1}^{n_1} I_{S_{i,23}}^{(M)}[q_2,q_3]\, H_2(x_i(q_2,q_3),q_2,q_3)\, H_3(x_i(q_2,q_3),q_2,q_3)\, (Curl_S^{(K)}F)_1(x_i(y,z),y,z) +$$

$$\sum_{j=1}^{n_2} I_{S_{j,13}}^{(M)}[q_1,q_3]\, H_1(q_1,y_j(q_1,q_3),q_3)\, H_3(q_1,y_j(q_1,q_3),q_3)\, (Curl_S^{(K)}F)_2(q_1,y_j(q_1,q_3),q_3) +$$

$$\sum_{k=1}^{n_3} I_{S_{k,12}}^{(M)}[q_1,q_2]\, H_1(q_1,q_2,z_k(q_1,q_2))\, H_2(q_1,q_2,z_k(q_1,q_2))(Curl_S^{(K)}F)_3(q_1,q_2,z_k(q_1,q_2)) =$$

$$\sum_{i=1}^{n_1} I_{S_{i,23}}^{(M)}[q_2,q_3]\, (I_{(K_2)}^{q_2,*}[q'_2]\, (D_{q'_2}^1 \hat{F}_3)(x_i(q'_2,q_3),q'_2,q_3) - I_{(K_3)}^{q_3,*}[q'_3]\, (D_{q'_3}^1 \hat{F}_2)(x_i(q_2,q'_3),q_2,q'_3)) +$$

$$\sum_{j=1}^{n_2} I_{S_{j,13}}^{(M)}[q_1,q_3]\, (I_{(K_3)}^{q_3,*}[q'_3]\, (D_{q'_3}^1 \hat{F}_1)(q_1,y_j(q_1,q'_3),q'_3) - I_{(K_1)}^{q_1,*}[q'_1]\, (D_{q'_1}^1 \hat{F}_3)(q'_1,y_j(q'_1,q_3),q_3)) +$$

$$\sum_{k=1}^{n_3} I_{S_{k,12}}^{(M)}[q_1,q_2]\, (I_{(K_1)}^{q_1,*}[q'_1]\, (D_{q'_1}^1 \hat{F}_2)(q'_1,q_2,z_k(q'_1,q_2)) - I_{(K_2)}^{q_2,*}[q'_2]\, (D_{q'_2}^1 \hat{F}_1)(q_1,q'_2,z_k(q_1,q'_2))).$$

For the case, $F_2 = F_3 = 0$, we have
$$F = F_1 = e_1 P,$$
and
$$(\hat{I}_S^{(M)}, Curl_S^{(K)}F) := \sum_{j=1}^{n_2} I_{S_{j,13}}^{(M)}[q_1,q_3]\, I_{(K_3)}^{q_3,*}[q'_3]\, (D_{q'_3}^1 \hat{P})(q_1,y_j(q_1,q'_3),q'_3) -$$

$$\sum_{k=1}^{n_3} I_{S_{k,12}}^{(M)}[q_1,q_2]\, I_{(K_2)}^{q_2,*}[q'_2]\, (D_{q'_2}^1 \hat{P})(q_1,q'_2,z_k(q_1,q'_2)).$$

This case will be considered in the proof of the general fractional Stokes theorem for Surface GF Curl with simple surface, which is given in the next section.

### 11.6 General Fractional Stokes Theorem for Surface GF Curl in OCC

Let us prove the general fractional Stokes theorem for surface general fractional curl operator, where surface consists of simple surfaces or surfaces parallel to the coordinate planes.

**Theorem 11.3** *(General Fractional Stokes Theorem for Surface GF Curl in OCC)*
Let $S \subset \mathbb{R}_{0,+}^3$ be a simple surface (or a piecewise simple surface $S \in \mathbb{P}(\mathbb{R}_{0,+}^3)$), and the vector field $F$ belongs to the set $\mathbb{F}_{S,H}^1(\mathbb{R}_{0,+}^3)$.
Then, the equation
$$(\hat{I}_{\partial S}^{(M)}, F) = (\hat{I}_S^{(M)}, Curl_S F)$$
holds, where $Curl_S^{(K)}$ is the surface GF Curl that is defined by Definition 11.19.

*Proof.* Let us prove the theorem for the vector field $F_1 = e_1 P(q_1, q_2, q_3)$. The general Stokes equarions for the vector fields $F_2 = e_2 F_2(q_1, q_2, q_3)$ and $F_3 = e_3 F_3(q_1, q_2, q_3)$ are proved similarly. The proof for the vector field $F = \sum_{k=0}^{3} e_k F_k$ is realized by the sum of the vector fields $F_k$.

Let $S \subset \mathbb{R}_{0,+}^3$ be a $Q_3$-simple and $Q_2$-simple surface that is described by the equation $q_3 = z(q_1, q_2)$ for $(q_1, q_2) \in S_{12} \subset \mathbb{R}_{0,+}^2$, and the equation $q_2 = y(q_1, q_3)$ for $(q_1, q_3) \in S_{13} \subset \mathbb{R}_{0,+}^2$.

The values of the function $P(q_1, q_2, q_3)$ on the line $L \subset \mathbb{R}_{0,+}^3$ are equal to the values of the function $P(q_1, q_2, z(q_1, q_2))$ on the line $L_{12} \subset \mathbb{R}_{0,+}^2$, which is a projection of the line $L$ onto the $Q_1 Q_2$-plane,

$$(I_{\partial S}^{(M)}, F) = \hat{I}_{\partial S}^{(M_1)}[q_1] P(q_1, q_2, q_3) = \hat{I}_{\partial S_{12}}^{(M_1)}[q_1] P(q_1, q_2, z(q_1, q_2)),$$

where $L = \partial S$ and $L_{12} = \partial S_{12}$.

Let us assume that $L_{12}$ consists of two $Q_2$-simple lines $L_{12a}$ and $L_{12b}$ that are described by equations $q_2 = y_a(q_1) \geq 0$ and $q_2 = y_b(q_1)$, where $y_b(q_1) \geq y_a(q_1) \geq 0$. Then

$$\hat{I}_{L_{12}}^{(M_1)}[q_1] P(q_1, q_2, z(q_1, q_2)) =$$

$$\hat{I}_{L_{12a}}^{(M_1)}[q_1] P(q_1, q_2, z(q_1, q_2)) + \hat{I}_{L_{12b}}^{(M_1)}[q_1] P(q_1, q_2, z(q_1, q_2)) =$$

$$I_{[a,b]}^{(M_1)}[q_1] \hat{P}(q_1, y_1(q_1), z(q_1, y_1(q_1))) - I_{[a,b]}^{(M_1)}[q_1] \hat{P}(q_1, y_2(q_1), z(q_1, y_2(q_1))) =$$

$$I_{[a,b]}^{(M_1)}[q_1](\hat{P}(q_1, y_1(q_1), z(q_1, y_1(q_1))) - \hat{P}(q_1, y_2(q_1), z(q_1, y_2(q_1)))), \qquad (210)$$

where

$$\hat{P}(q_1, q_2, q_3) := H_1(q_1, q_2, q_3) P(q_1, q_2, q_3).$$

Using the fundamental theorem of general fractional calculus, expression (210) can be written as

$$I_{[a,b]}^{(M_1)}[q_1](\hat{P}(q_1, y_1(q_1), z(q_1, y_1(q_1))) - \hat{P}(q_1, y_2(q_1), z(q_1, y_2(q_1)))) =$$

$$- I_{[a,b]}^{(M_1)}[q_1] I_{[y_1(q_1), y_2(q_1)]}^{(M_2)}[q_2] D_{(K_2)}^{q_2,*}[q'_2] \hat{P}(q_1, q'_2, z(q_1, q'_2)). \qquad (211)$$

Then using the definition of the general fractional derivative and the standard chain rule for the first-order derivative, we get

$$- I_{[a,b]}^{(M_1)}[q_1] I_{[y_1(q_1), y_2(q_1)]}^{(M_2)}[q_2] D_{(K_2)}^{q_2,*}[q'_2] \hat{P}(q_1, q'_2, z(q_1, q'_2)) =$$

$$- I_{[a,b]}^{(M_1)}[q_1] I_{[y_1(q_1), y_2(q_1)]}^{(M_2)}[q_2] I_{(K_2)}^{q_2,*}[q'_2] \frac{\partial \hat{P}(q_1, q'_2, z(q_1, q'_2))}{\partial q'_2} =$$

$$- I_{[a,b]}^{(M_1)}[q_1] I_{[y_1(q_1), y_2(q_1)]}^{(M_2)}[q_2] I_{(K_2)}^{q_2,*}[q'_2] \left(\left(\frac{\partial \hat{P}(q_1, q'_2, q_3)}{\partial q'_2}\right)_{q_3 = z(q_1, q'_2)} + \right.$$

$$\left.\left(\frac{\partial \hat{P}(q_1, q'_2, q_3)}{\partial q_3}\right)_{q_3 = z(q_1, q'_2)} \frac{\partial z(q_1, q'_2)}{\partial q'_2}\right) =$$

$$- I^{(M_1)}_{[a,b]}[q_1] \, I^{(M_2)}_{[y_1(q_1),y_2(q_1)]}[q_2] \, I^{q_2,*}_{(K_2)}[q'_2] \left( \frac{\partial \hat{P}(q_1,q'_2,q_3)}{\partial q'_2} \right)_{q_3=z(q_1,q'_2)} - \quad (212)$$

$$- I^{(M_1)}_{[a,b]}[q_1] \, I^{(M_2)}_{[y_1(q_1),y_2(q_1)]}[q_2] \, I^{q_2,*}_{(K_2)}[q'_2] \left( \frac{\partial \hat{P}(q_1,q'_2,q_3)}{\partial q_3} \right)_{q_3=z(q_1,q'_2)} \frac{\partial z(q_1,q'_2)}{\partial q'_2}. \quad (213)$$

Let us use the property

$$I^{(M_2)}_{[y_1,y_2]}[q_2] \, I^{q_2,*}_{(K_2)}[q'_2] U(q'_2) = (M_2 * (K_2 * U))(y_2) - (M_2 * (K_2 * U))(y_1) =$$

$$((M_2 * K_2) * U)(y_2) - ((M_2 * K_2) * U)(y_1) =$$

$$(\{1\} * U)(y_2) - (\{1\} * U)(y_1) = \int_{y_1}^{y_2} U(q_2) \, dq_2$$

for the expression of term (213) in the form

$$I^{(M_2)}_{[y_1(q_1),y_2(q_1)]}[q_2] \, I^{q_2,*}_{(K_2)}[q'_2] \left( \frac{\partial \hat{P}(q_1,q'_2,q_3)}{\partial q_3} \right)_{q_3=z(q_1,q'_2)} \frac{\partial z(q_1,q'_2)}{\partial q'_2} =$$

$$\int_{y_1(q_1)}^{y_2(q_1)} \left( \frac{\partial \hat{P}(q_1,q'_2,q_3)}{\partial q_3} \right)_{q_3=z(q_1,q'_2)} \frac{\partial z(q_1,q'_2)}{\partial q'_2} \, dq'_2. \quad (214)$$

Then we can use the equation, which is used in the standard Stokes theorem, $(\{1\} * U)(q_3) = (M_3 * K_3 * U)(q_3)$, and the fundamental theorem of GFC in the form

$$\int_{y_1(q_1)}^{y_2(q_1)} \left( \frac{\partial \hat{P}(q_1,q'_2,q_3)}{\partial q_3} \right)_{q_3=z(q_1,q'_2)} \frac{\partial z(q_1,q'_2)}{\partial q'_2} \, dq'_2 =$$

$$- \int_{z_1(q_1)}^{z_2(q_1)} \left( \frac{\partial \hat{P}(q_1,q_2,q'_3)}{\partial q'_3} \right)_{q_2=y(q_1,q'_3)} \, dq'_3 =$$

$$- I^{(M_3)}_{[z_1(q_1),z_2(q_1)]}[q_3] \, I^{q_3,*}_{(K_3)}[q'_3] \left( \frac{\partial \hat{P}(q_1,q_2,q'_3)}{\partial q'_3} \right)_{q_2=y(q_1,q'_3)}. \quad (215)$$

Using (215), expression (213) can be written in the form

$$- I^{(M_1)}_{[a,b]}[q_1] \, I^{(M_2)}_{[y_1(q_1),y_2(q_1)]}[q'_2] \, I^{q_2,*}_{(K_2)}[q'_2] \left( \frac{\partial \hat{P}(q_1,q'_2,q_3)}{\partial q_3} \right)_{q_3=z(q_1,q'_2)} \frac{\partial z(q_1,q'_2)}{\partial q'_2} =$$

$$I^{(M_1)}_{[a,b]}[q_1] \, I^{(M_3)}_{[z_1(q_1),z_2(q_1)]}[q_3] \, I^{q_3,*}_{(K_3)}[q'_3] \left( \frac{\partial \hat{P}(q_1,q_2,q'_3)}{\partial q'_3} \right)_{q_2=y(q_1,q'_3)}. \quad (216)$$

Therefore, for two terms (212) and (213), we obtain the equations

$$- I^{(M_1)}_{[a,b]}[q_1] \, I^{(M_2)}_{[y_1(q_1),y_2(q_1)]}[q_2] \, I^{q_2,*}_{(K_2)}[q'_2] \left( \frac{\partial \hat{P}(q_1,q'_2,q_3)}{\partial q'_2} \right)_{q_3=z(q_1,q'_2)} =$$

$$- I^{(M_1)}_{[a,b]}[q_1] \, I^{(M_2)}_{[y_1(q_1),y_2(q_1)]}[q_2] \frac{1}{H_1(q_1,q_2,z(q_1,q_2)) \, H_2(q_1,q_2,z(q_1,q_2))} \cdot$$

$$\cdot I^{q_2,*}_{(K_2)}[q'_2] \left( \frac{\partial \hat{P}(q_1,q'_2,q_3)}{\partial q'_2} \right)_{q_3=z(q_1,q'_2)} = \hat{I}^{(M)}_{S_{12}}[q_1,q_2] \, (Curl^{(K)}_S F_1)_3, \quad (217)$$

and

$$I^{(M_1)}_{[a,b]}[q_1] \, I^{(M_3)}_{[z_1(q_1),z_2(q_1)]}[q_3] \, I^{q_3,*}_{(K_3)}[q'_3] \left( \frac{\partial \hat{P}(q_1,q_2,q'_3)}{\partial q'_3} \right)_{q_2=y(q_1,q'_3)} =$$

$$I^{(M_1)}_{[a,b]}[q_1] \, I^{(M_3)}_{[z_1(q_1),z_2(q_1)]}[q_3] \frac{1}{H_1(q_1,q_2,z(q_1,q_2)) \, H_2(q_1,q_2,z(q_1,q_2))} \cdot$$

$$\cdot \ I_{(K_3)}^{q_3,*}[q'_3] \left(\frac{\partial \hat{P}(q_1,q_2,q'_3)}{\partial q'_3}\right)_{q_2=y(q_1,q'_3)} = \hat{I}_{S_{13}}^{(M)}[q_1,q_3] \, (Curl_S^{(K)} F_1)_2, \tag{218}$$

for the vector field $F_1 = e_1 P(q_1, q_2, q_3)$.

As a result, we proved the general Stokes theorem for the surface FG Curl and the vector field $F_1 = e_1 P(q_1, q_2, q_3)$.

$$(\hat{I}_{\partial S}^{(M_1)}, F_1) = (\hat{I}_S^{(M)}, Curl_S^{(K)} F_1).$$

Equations for the remaining components of the vector field $F_2 = e_2 F_2(q_1, q_2, q_3)$ and $F_3 = e_3 F_3(q_1, q_2, q_3)$ are proved similarly.

## 11.7 General Fractional Gauss Theorem in OCC

### 11.7.1 Definition of Triple GFI by Iterated GFI in OCC

Let us use the concept 9.1 of $Z$-simple region for OCC in $\mathbb{R}_{0,+}^3$.

**Definition 11.20** Let $W \subset \mathbb{R}_{0,+}^3$ be $Q_3$-simple domain in OCC that is bounded above and below by smooth surfaces $S_{12b}$, $S_{12a}$ descrined by the equations

$$q_3 = z_1(q_1, q_2) \geq 0, \quad q_3 = z_2(q_1, q_2) \geq 0, \tag{219}$$

where the functions are continuous in the closed domain $W_{12}$ that is a projection of the region $W$ onto the $Q_1 Q_2$-plane, and $z_3(q_1, q_2) \geq z_1(q_1, q_2)$ for all $(q_1, q_2) \in W_{12}$.

Let scalar field $U(q_1, q_2, q_3)$ be satisfy the condition

$$I_{[z_1(q_1,q_2),z_2(q_1,q_2)]}^{(M_3)}[q_3] \, \tilde{U}(q_1, q_2, q_3) \in C_{-1}(\mathbb{R}_+^2), \tag{220}$$

where

$$\tilde{U}(q_1, q_2, q_3) := H_1(q_1, q_2, q_3) H_2(q_1, q_2, q_3) H_3(q_1, q_2, q_3) U(q_1, q_2, q_3). \tag{221}$$

Then the triple general fractional integral in OCC is defined by the equation

$$\hat{I}_W^{(M)}[q_1, q_2, q_3] \, U(q_1, q_2, q_3) := I_{W_{12}}^{(M)}[q_1, q_2] \, I_{[z_1(q_1,q_2),z_2(q_1,q_2)]}^{(M_3)}[q_3] \, \tilde{U}(q_1, q_2, q_3),$$

where $I_{W_{12}}^{(M)}[q_1, q_2]$ is double GFI.

**Definition 11.21** Let $W \subset \mathbb{R}_{0,+}^3$ be $Q_3$-simple domain that is that is bounded above and below by smooth surfaces $S_{12b}$, $S_{12a}$ descrined by equation (219). Let $W_{12} \subset \mathbb{R}_{0,+}^2$ be projection of $W$ on the $Q_1 Q_2$-plane such that $W_{12}$ is $Q_2$-simple region in $Q_1 Q_2$-plane that is bounded by the lines $q_2 = y_1(q_1)$ and $q_2 = y_2(q_1)$, where $q_2 = y_1(q_1)$ and $q_2 = y_2(q_1)$ are continuous functions on the interval $[a, b]$, $b > a \geq 0$ and $y_2(q_1) \geq y_1(q_1) > 0$ for all $q_1 \in [a, b]$.

Let scalar field $U(q_1, q_2, q_3)$ satisfy condition (220) and

$$I_{[y_1(q_1),y_2(q_1)]}^{(M_2)}[q_2] \, I_{[z_1(q_1,q_2),z_2(q_1,q_2)]}^{(M_3)}[q_3] \, \tilde{U}(q_1, q_2, q_3) \in C_{-1}(0, \infty). \tag{222}$$

Then the triple general fractional integral (triple GFI) is defined in the form

$$\hat{I}_W^{(M)}[q_1, q_2, q_3] \, U(q_1, q_2, q_3) :=$$

$$I_{[a,b]}^{(M_1)}[q_1] \, I_{[y_1(q_1),y_2(q_1)]}^{(M_2)}[q_2] \, I_{[z_1(q_1,q_2),z_2(q_1,q_2)]}^{(M_3)}[q_3] \, \tilde{U}(q_1, q_2, q_3).$$

A volume general fractional integral (volume GFI) of a scalar field is a triple general fractional integral for the $Z$-simple region $W \subset \mathbb{R}^3_{0,+}$.

**Example 11.1** Using the parallelepiped region
$$W := \{0 \leq q_1 \leq a_1, \quad 0 \leq q_2 \leq a_2, \quad 0 \leq q_3 \leq a_3\},$$
the volume general fractional integral can be written as
$$I_W^{(M)}[q_1, q_2, q_3] f(q_1, q_2, q_3) =$$

$$\int_0^{a_1} dq_1 \int_0^{a_2} dq_2 \int_0^{a_3} dq_3 \, M_1(a_1 - q_1) M_2(a_2 - q_2) M_3(a_3 - q_3) \tilde{f}(q_1, q_2, q_3). \quad (223)$$

For the power-law kernels, the general fractional flux (223) has the form
$$\hat{I}_W^{(M)}[q_1, q_2, q_3] f(q_1, q_2, q_3) =$$

$$\int_0^{a_1} dq_1 \int_0^{a_2} dq_2 \int_0^{a_3} dq_3 \, \frac{(a_1-q_1)^{\alpha_1-1}}{\Gamma(\alpha_1)} \frac{(a_2-q_2)^{\alpha_2-1}}{\Gamma(\alpha_2)} \frac{(a_3-q_3)^{\alpha_3-1}}{\Gamma(\alpha_3)} \tilde{f}(q_1, q_2, q_3). \quad (224)$$

For $\alpha_1 = \alpha_2 = \alpha_3 = 1$, equation (224) is gives as
$$\hat{I}_W^{(M)}[q_1, q_2, q_3] f(q_1, q_2, q_3) = \int_0^{a_1} dq_1 \int_0^{a_2} dq_2 \int_0^{a_3} dq_3 \, \tilde{f}(q_1, q_2, q_3), \quad (225)$$
which is the standard volume integral for the function $f(q_1, q_2, q_3)$.

### 11.7.2 General Fractional Divergence in OCC

In this section, we give the definition of general fractional divergence in OCC for $\mathbb{R}^3_{0,+}$.

Let us define sets of vector fields that will be used in the definition of theregional general fractional divergence.

**Definition 11.22** Let $F(q_1, q_2, q_3)$ be a vector field that satisfies the conditions
$$D_{(K_1)}^{q_1,*}[q'_1]\tilde{F}_1(q'_1, q_2, q_3) \in C_{-1}(\mathbb{R}^3_+),$$

$$D_{(K_2)}^{q_2,*}[q'_2]\tilde{F}_2(q_1, q'_2, q_3) \in C_{-1}(\mathbb{R}^3_+),$$

$$D_{(K_3)}^{q_3,*}[q'_3]\tilde{F}_3(q_1, q_2, q'_3) \in C_{-1}(\mathbb{R}^3_+),$$

where
$$\tilde{F}_k(q_1, q_2, q_3) = H_1(q_1, q_2, q_3) H_2(q_1, q_2, q_3) H_3(q_1, q_2, q_3) \frac{F_k(q_1,q_2,q_3)}{H_k(q_1,q_2,q_3)}.$$
Then the set of such vector fields will be denoted as $\mathbb{F}^1_{-1,Div,H}(\mathbb{R}^3_+)$.

We can also consider vector fields $F(q_1, q_2, q_3)$ that satisfy the condition $\hat{F}_k \in C^1_{-1}(\mathbb{R}^3_+)$.

**Definition 11.23** Let $F(q_1, q_2, q_3)$ be a vector field that satisfies the conditions
$$\tilde{F}_k(q_1, q_2, q_3) \in C^1_{-1,H}(\mathbb{R}^3_+)$$
for all $k = 1,2,3$.
Then the set of such vector fields will be denoted as $C^1_{-1,H}(\mathbb{R}^3_+)$.

Let us define the general fractional divergencein OCC.

**Definition 11.24** Let $F(q_1, q_2, q_3)$ be a vector field that belongs to the set $\mathbb{F}^1_{-1,Div,H}(\mathbb{R}^3_+)$ or $C^1_{-1,H}(\mathbb{R}^3_+)$.

Then the regional general fractional divergence $Div_W^{(K)}$ in OCC for the region $W = \mathbb{R}^3_{0,+}$ is defined as

$$Div_W^{(K)} F =$$

$$\frac{1}{H_1 H_2 H_3} (D_{(K_1)}^{q_1,*}[q'_1]\tilde{F}_1(q'_1, q_2, q_3) + D_{(K_2)}^{q_2,*}[q'_2]\tilde{F}_2(q_1, q'_2, q_3) + D_{(K_3)}^{q_3,*}[q'_3]\tilde{F}_3(q_1, q_2, q'_3)),$$
(226)

where $H_k = H_k(q_1, q_2, q_3)$.

Divergence (226) is the regional general fractional differential operator in OCC. We can also define the line and surface general fractional divergence in OCC.

**Remark 11.2** The equation, whcih defines the regional general fractional divergence in OCC, can be written in compact form. If the vector field
$$F(q_1, q_2, q_3) = \sum_{k=1}^{3} e_k F_k(q_1, q_2, q_3)$$
belongs to the function space $C^1_{-1}(\mathbb{R}^3_+)$, then the regional general fractional divergence in OCC for the region $W = \mathbb{R}^3_+$ can be defined as
$$(Div_W^{(K)} \tilde{F})(q_1, q_2, q_3) = (\widehat{D}_W^{(K)}, F) = \frac{1}{H_1 H_2 H_3} \sum_{k=1}^{3} D_{(K_k)}^{q_k,*}[q'_k]\tilde{F}_k,$$
where
$$\widehat{D}_W^{(K)} := e_1 D_{(K_1)}^{q_1,*}[q'_1] + e_2 D_{(K_2)}^{q_2,*}[q'_2] + e_3 D_{(K_3)}^{q_3,*}[q'_3].$$

**Remark 11.3** The general fractional divergence in OCC can be defined not only for $W = \mathbb{R}^3_{0,+}$, but also for regions $W \subset \mathbb{R}^3_{0,+}$, surfaces $S \subset \mathbb{R}^3_{0,+}$ and line $L \subset \mathbb{R}^3_{0,+}$. Note that the general fractional divergence in OCC for regions $W \subset \mathbb{R}^3_{0,+}$ is used in the general fractional Gausss theorem for OCC to be given and proved in the following sections.

## 11.8  General Fractional Gauss Theorem for OCC

Let us define a set of vector field, for which general fractional Gauss theorem will be formulated.

**Definition 11.25** Let $W$ in $\mathbb{R}^3_{0,+}$ be $Q_3$-simple region such that $W$ is a piecewise $Q_2$-simple and $Q_1$-simple region
$$W = \bigcup_{k=1}^{n} W_{Q1,k} = \bigcup_{j=1}^{m} W_{Q2,j},$$
where $W_{Q1,k}$ are the $Q!X$-simple regions that is described by $q_1 = x_{k,1}(q_2, q_3)$, $q_1 = x_{k,2}(q_2, q_3)$ for $(q_2, q_3) \in D_{12}$, and $W_{Q2,j}$ are the $Q_2$-simple regions that are described by the functions $q_2 = y_{j,1}(q_1, q_3)$, $q_2 = y_{j,2}(q_1, q_3)$ for $(q_1, q_3) \in D_{13}$.

Let vector field $F(q_1, q_2, q_3)$ satisfy the conditions
$$\tilde{F}_1(x_{k,2}(q_2, q_3), q_2, q_3), \tilde{F}_1(x_{k,1}(q_2, q_3), q_2, q_3) \in C_{-1}(\mathbb{R}^2_+),$$

$$\tilde{F}_2(q_1, y_{j,2}(q_1, q_3), q_3), \tilde{F}_2(q_1, y_{j,1}(q_1, q_3), q_3) \in C_{-1}(\mathbb{R}^2_+),$$

$$\tilde{F}_3(q_1, q_2, z_2(q_1, q_2)), \tilde{F}_3(q_1, q_2, z_1(q_1, q_2)) \in C_{-1}(\mathbb{R}^2_+),$$

where
$$\tilde{F}_k(q_1, q_2, q_3) := H_1(q_1, q_2, q_3) H_2(q_1, q_2, q_3) H_3(q_1, q_2, q_3) \frac{F_k(q_1, q_2, q_3)}{H_k(q_1, q_2, q_3)}.$$

Then the set of such vector fields $F$ will be denoted as $\tilde{\mathbb{F}}_{-1,\partial W}(\mathbb{R}^2_+)$.

**Definition 11.26** Let vector field $F(q_1, q_2, q_3)$ satisfy the conditions
$$f_1(q_1, q_2, q_3) := D^{q_1,*}_{(K_1)}[q'_1] \tilde{F}_1(q'_1, q_2, q_3) \in C_{-1}(\mathbb{R}^3_+),$$

$$f_2(q_1, q_2, q_3) := D^{q_2,*}_{(K_2)}[q'_2] \tilde{F}_2(q_1, q'_2, q_3) \in C_{-1}(\mathbb{R}^3_+),$$

$$f_3(q_1, q_2, q_3) := D^{q_3,*}_{(K_3)}[q'_3] \tilde{F}_3(q_1, q_2, q'_3) \in C_{-1}(\mathbb{R}^3_+).$$

Then the set of such vector fields $\mathbb{F}$ will be denoted as $\tilde{\mathbb{F}}^1_{-1}(\mathbb{R}^3_+)$.

**Theorem 11.4** *(General fractional Gauss theorem for OCC)*

Let $W$ in $\mathbb{R}^3_{0,+}$ be $Q_3$-simple region such that $W$ is a piecewise $Q_1$-simple and $Q_2$-simple region. Let $W$ is bounded above and below by smooth surfaces $S_{12b}$, $S_{12a}$, which is described by equations (219), and a lateral surface $S_3$, whose generatrices are parallel to the $Q_3$-axis.

Let the vector field $F(q_1, q_2, q_3)$ belongs to the sets $\tilde{\mathbb{F}}_{-1}(W)$ and $\tilde{\mathbb{F}}^1_{-1}(\mathbb{R}^3_+)$.

Then the general fractional Gauss equation has the from
$$\hat{I}^{(M)}_W[q_1, q_2, q_3] (Div^{(K)}_W F) = (\hat{I}_{\partial W}, F)$$

that can be written as
$$I^{(M)}_W[q_1, q_2, q_3] (D^{q_1,*}_{(K_1)}[q'_1] \tilde{F}_1(q'_1, q_2, q_3) + D^{q_2,*}_{(K_2)}[q'_2] \tilde{F}_2(q_1, q'_2, q_3) + D^{q_3,*}_{(K_3)}[q'_3] \tilde{F}_3(q_1, q_2, q'_3)) =$$

$$\hat{I}^{(M)}_S[q_2, q_3] F_1(q_1, q_2, q_3) + \hat{I}^{(M)}_S[q_1, q_2] F_2(q_1, q_2, q_3) + \hat{I}^{(M)}_S[q_1, q_2] F_3(q_1, q_2, q_3).$$

*Proof.* Let us consider the triple GFI
$$\hat{I}^{(M)}_W[q_1, q_2, q_3] U = I^{(M)}_W[q_1, q_2, q_3] \tilde{U} = I^{(M)}_{W_{12}}[q_1, q_2] I^{(M_3)}_{[z_1(q_1,q_2), z_2(q_1,q_2)]}[q_3] \tilde{U}(q_1, q_2, q_3)$$
(227)

for the scalar field
$$U(q_1, q_2, q_3) = \frac{1}{H_1(q_1,q_2,q_3) H_2(q_1,q_2,q_3) H_3(q_1,q_2,q_3)} D^{q_3,*}_{(K_3)}[q'_3] \tilde{F}_3(q_1, q_2, q'_3) \in C_{-1}(0, \infty), (228)$$
where $F_3(q_1, q_2, q'_3) \in C^1_{-1}(0, \infty)$ for each point $(q_1, q_2) \in D_{12}$.

Using equations (221), (227) and (228), we get

$$I_W^{(M)}[q_1,q_2,q_3]\widetilde{U} = I_{W_{12}}^{(M)}[q_1,q_2] \, D_{(K_3)}^{q_3,*}[q'_3] \, \tilde{F}_3(q_1,q_2,q'_3), \tag{229}$$

since $H_1 H_2 H_3 (H_1 H_2 H_3)^{-1} = 1$.

Using the second fundamental theorem of GFC, we get

$$I_{[z_1(q_1,q_2),z_2(q_1,q_2)]}^{(M_3)}[q_3] \, D_{(K_3)}^{q_3,*}[q'_3] \, \tilde{F}_3(q_1,q_2,q'_3) =$$

$$I_{(M_3)}^{z_2(q_1,q_2)}[q_3] \, D_{(K_3)}^{q_3,*}[q'_3] \, F_3(q_1,q_2,q'_3) - I_{(M_3)}^{z_1(q_1,q_2)}[q_3] \, D_{(K_3)}^{q_3,*}[q'_3] \, \tilde{F}_3(q_1,q_2,q'_3) =$$

$$(\tilde{F}_3(q_1,q_2,z_2(q_1,q_2)) - \tilde{F}_3(q_1,q_2,0)) - (\tilde{F}_3(q_1,q_2,z_2(q_1,q_2)) - \tilde{F}_3(q_1,q_2,0))) =$$

$$\tilde{F}_3(q_1,q_2,z_2(q_1,q_2)) - \tilde{F}_3(q_1,q_2,z_1(q_1,q_2)). \tag{230}$$

Substituting expression (230) into equation (229) and assuming that the condition

$$\tilde{F}_3(q_1,q_2,z_1(q_1,q_2)), \tilde{F}_3(q_1,q_2,z_2(q_1,q_2)) \in C_{-1}(\mathbb{R}^2_+),$$

is satisfied, we obtain the equation

$$I_W^{(M)}[q_1,q_2,q_3] \, D_{(K_3)}^{q_3,*}[q'_3] \, \tilde{F}_3(q_1,q_2,q'_3) = I_{W_{12}}^{(M)}[q_1,q_2] \, (\tilde{F}_3(q_1,q_2,z_2(q_1,q_2)) - \tilde{F}_3(q_1,q_2,z_1(q_1,q_2))).$$

Therefore, the triple GFI can be represented through the surface GFIs in the form

$$I_W^{(M)}[q_1,q_2,q_3] D_{(K_3)}^{q_3,*}[q'_3] \, \tilde{F}_3(q_1,q_2,q'_3) =$$

$$I_{S_{12b}}^{(M)}[q_1,q_2]\tilde{F}_3(q_1,q_2,q_3) + I_{S_{12a}}^{(M)}[q_1,q_2]\tilde{F}_3(q_1,q_2,q_3),$$

where the surface GFIs are represented by the definition in the form

$$I_{S_{12b}}^{(M)}[q_1,q_2] \, \tilde{F}_3(q_1,q_2,q_3) = I_{W_{12}}^{(M)}[q_1,q_2] \, \tilde{F}_3(q_1,q_2,z_2(q_1,q_2)),$$

$$I_{S_{12a}}^{(M)}[q_1,q_2]\tilde{F}_3(q_1,q_2,q_3) = I_{W_{12}}^{(M)}[q_1,q_2] \, \tilde{F}_3(q_1,q_2,z_1(q_1,q_2)).$$

Then, taking into account that the surface GFI oves the surface $S_z$ is equati to zero

$$I_{S_3}^{(M)}[q_1,q_2] \, \tilde{F}_3(q_1,q_2,q_3) = 0,$$

we obtain

$$I_W^{(M)}[q_1,q_2,q_3] \, D_{(K_3)}^{q_3,*}[z'] \, \tilde{F}_3(q_1,q_2,q'_3) = I_S^{(M)}[q_1,q_2] \, \tilde{F}_3(q_1,q_2,q_3) =$$

$$\hat{I}_S^{(M)}[q_1,q_2] \, F_3(q_1,q_2,q_3),$$

where $\partial W = S = S_{12a} \cup S_{12b} \cup S_3$ is closed surface that contains the region $W$ inside the surface $\partial W$.

Further, similarly with the proof of general fractional Gauss theorem for the Cartesian coordinate system and and by analogy with the transformation described above, we obtain the equation of the theorem to be proved.

**Example 11.2** As an example, let us give the general fractional Gauss theorem in cylindrical coordinates. Let us consider the simple domain in the form

$$W := \{(r,\varphi,z): \; 0 \leq r_1 \leq r \leq r_2, \;\; 0 \leq \varphi_1 \leq \varphi \leq \varphi_2 < 2\pi, \;\; 0 \leq z_1 \leq z \leq z_2\}. \tag{231}$$

The vertices of the domain are the following points

$$A_1(r_1,\varphi_1,z_1), \quad A_2(r_2,\varphi_1,z_1), \quad B_1(r_1,\varphi_2,z_1), \quad B_2(r_2,\varphi_2,z_1),$$

$$A'_1(r_1,\varphi_1,z_2), \quad A'_2(r_2,\varphi_1,z_2), \quad B'_1(r_1,\varphi_2,z_2), \quad B'_2(r_2,\varphi_2,z_2).$$

Let $F \in \mathbb{F}^1_{-1}(W)$, then

$$\widehat{I_W^{(M)}}[r,\varphi,z]\, Div_W^{(K)} F \;=\; (\widehat{I_{\partial W}^{(M)}}[r,\varphi,z], F(r,\varphi,z)).$$

For simple domain (231) the cylindrical general fractional Gauss equation has the form

$$\widehat{I_{[r_1,r_2]}^{(M_r)}}[r]\, \widehat{I_{[\varphi_1,\varphi_2]}^{(M_\varphi)}}[\varphi]\, \widehat{I_{[z_1,z_2]}^{(M_z)}}[z]\, (\tfrac{1}{r}\, D_{(K_r)}^{r,*}[r']\,(r'\,F_r(r',\varphi,z)) \;+$$

$$\tfrac{1}{r}\, D_{(K_\varphi)}^{\varphi,*}[\varphi']\, F_\varphi(r,\varphi',z) \;+\; D_{(K_z)}^{z',*}[z']\, F_z(r,\varphi,z')) \;=$$

$$\widehat{I_{\partial W_{\varphi,z}}^{(M)}}[\varphi,z]\, F_r(r,\varphi,z) \;+\; \widehat{I_{\partial W_{r,z}}^{(M)}}[r,z]\, F_\varphi(r,\varphi,z) \;+\; \widehat{I_{\partial W_{r,\varphi}}^{(M)}}[r,\varphi]\, F_z(r,\varphi,z).$$

## 12 Conclusion

The general fractional vector calculus (General FVC) is proposed as a generalization of the fractional vector calculus suggested in [30, 31]. The formulation of the General FVC is based on the results in the calculus of general fractional integrals and derivatives that is proposed in Luchko's work [50].

The formulation of General FVC is self-consistent form, i.e. definitions of fractional generalizations of differential and integral vector operators are consistent with each other, and generalizations of fundamental theorems of vector calculus were proved. In this paper, the definitions of general fractional integral vector operators: the general fractional circulation, general fractional flux and general fractional volume integral are proposed. Definitions of general fractional differential vector operators, including the regional and line general fractional gradients, the regional and surface general fractional curl operators, the general fractional divergence, are suggested. Fundamental theorems of epy General FVC, which are general analogs of the standard gradient, Green's, Stokes' and Gauss's theorems, are proved for simple and complex regions. Let us emphasize that the fractional vector analogs of fundamental theorems (such as the gradient, Stock's and Gauss theorems) are not fulfilled for all type (regional, surface and line) of the general fractional vector operators (the gradient, curl and divergence). In the general case, the general fractional (GF) gradient theorem should be considered for the line GF Gradient, the FG Stock's theorem should be considered for the surface GF Curl operator, and the GF Green theorem should be considered for the regional GF Divergence. This is due to violation of the chain rule for general fractional derivatives. The General FVC for orthogonal curvilinear coordinates, which includes general fractional vector operators for the spherical and cylindrical coordinates, is described.

The proposed General FVC can be used as a mathematical tool in general fractional dynamics (GFDynamics) [41, 42], in which non-locality in space is taken into account in general form.